\documentclass{article}
 \usepackage{amssymb}                   
 \usepackage{amsmath}                   
 \usepackage{amsthm}                    
 \usepackage{latexsym}                  

 \usepackage{amscd}
 \usepackage{amsmath}
 \usepackage[all]{xy}

 \newtheorem{theorem}{Theorem}[section]
 
 \newtheorem{corollary}[theorem]{Corollary}
 \newtheorem{lemma}[theorem]{Lemma}
 \newtheorem{proposition}[theorem]{Proposition}
 \newtheorem{claim}[theorem]{Claim}
 \newtheorem{problem}[theorem]{Problem} 
 \newtheorem{condition}[theorem]{Condition}

 \newtheorem*{lemmaA}{Lemma A}
 \newtheorem*{lemmaB}{Lemma B}
 \newtheorem*{keyproposition}{Key Proposition}

 \theoremstyle{definition}\newtheorem{definition}[theorem]{Definition}
 \newtheorem*{division}{Division into Cases}
 
 \theoremstyle{remark}\newtheorem{remark}[theorem]{Remark}
 
 \newtheorem{example}[theorem]{Example}

 \numberwithin{equation}{section}

 \newcommand{\ZZ}{\ensuremath{\mathbb{Z}}}

 \newcommand{\LL}{\ensuremath{\mathbf{L}}}
 \newcommand{\ltensor}{\overset{\LL}{\otimes}}
 \newcommand{\CC}{\ensuremath{\mathbb{C}}}
 \newcommand{\Aut}{\ensuremath{\operatorname{Aut}}}
 \newcommand{\Ext}{\ensuremath{\operatorname{Ext}}}
 \newcommand{\Hom}{\ensuremath{\operatorname{Hom}}}
 \newcommand{\ext}{\ensuremath{\mathcal{E}xt}}
 \renewcommand{\hom}{\ensuremath{\mathcal{H}om}}
 \newcommand{\Pic}{\ensuremath{\operatorname{Pic}}}
 \newcommand{\Spec}{\ensuremath{\operatorname{Spec}}}
 \newcommand{\Supp}{\ensuremath{\operatorname{Supp}}}
 \newcommand{\Ker}{\ensuremath{\operatorname{Ker}}}
 \newcommand{\Coker}{\ensuremath{\operatorname{Coker}}}
 \newcommand{\Image}{\ensuremath{\operatorname{Im}}}
 \newcommand{\PP}{\mathbb P}
 \newcommand{\Auteq}{\ensuremath{\operatorname{Auteq}}}
 \newcommand{\Exc}{\ensuremath{\operatorname{Exc}}}

 \newcommand{\GL}{\operatorname{GL}}

 \newcommand{\rank}{\ensuremath{\operatorname{rank}}}

 \newcommand{\ev}{\ensuremath{\operatorname{ev}}}
 \newcommand{\owe}{\ensuremath{\mathcal{O}}}
 \newcommand{\length}{\operatorname{length}} 
 \newcommand{\mc}{\mathcal}

  \newcommand{\Span}[1]{\left<#1\right>}

 \newcommand{\cohom}{\mathcal{H}}
 \newcommand{\coh}{\operatorname{Coh}}
 
 \newcommand{\tor}{\mbox{\small\it tor}}
 \newcommand{\eqcat}{D^G_{\{0\}}(\CC^2)}
 \newcommand{\GExt}{\operatorname{G-Ext}}
 \newcommand{\GHom}{\operatorname{G-Hom}}

\newcommand{\zero}{
\setlength{\unitlength}{1ex}
\begin{picture}(2, 2)(-1, -1)
\put(0,0){\circle{2}}
\put(0,0){\makebox(0,0){\tiny$0$}}
\end{picture}
}
\newcommand{\mone}{
\setlength{\unitlength}{1ex}
\begin{picture}(2, 2)(-1, -1)
\put(0,0){\circle{2}}
\put(0,0){\makebox(0,0){\tiny $-1$}}
\end{picture}
}

\newcommand{\no}{
\setlength{\unitlength}{1ex}
\begin{picture}(2, 2)(-1, -1)
\put(0,0){\circle{2}}
\end{picture}
}

\newcommand{\cia}{\setlength{\unitlength}{1ex}
\begin{picture}(2, 2)(-1, -1)
\put(0,0){\circle{2}}
\put(0,0){\makebox(0,0){\tiny$a$}}
\end{picture}
}

\title{Autoequivalences of derived categories on the minimal resolutions
of $A_n$-singularities on surfaces}

\author{Akira Ishii \and Hokuto Uehara}

\begin{document}
\maketitle

\begin{abstract}
In this article, we study the group of
autoequivalences of derived categories of coherent 
sheaves on the minimal resolution of $A_n$-singularities on surfaces. 
Our main result is to find generators of this group.
\end{abstract}

 \section{Introduction}\label{section:introduction}

Let $X$ be a smooth projective variety over $\CC$ and $D(X)(=D^b(\coh X))$ the bounded derived category of 
coherent sheaves on $X$.
$D(X)$ carries a lot of geometric information on $X$;
for instance, Bondal and Orlov show in \cite{Bondal:rvfd} that
if $K_X$ or $-K_X$ is ample, then $X$ can be entirely reconstructed from $D(X)$.
To the contrary, there are examples of mutually non-isomorphic varieties
$X$ and $Y$ having mutually equivalent derived categories.
Given a smooth projective variety $X$, it is an interesting problem
to find all the varieties $Y$ with $D(X) \cong D(Y)$.
In dimension $2$, the answer is given by Bridgeland and Maciocia in \cite{Bridgeland:csed},
and Kawamata \cite{Kawamata:deke} and in dimension $3$ some results are shown by Toda \cite{To03}. 
Moreover Orlov gives a satisfactory answer in \cite{Orlov:dcav} 
to this problem for the case where $X$ is an abelian variety.
The subject of this paper is related to another important problem:
\begin{problem}\label{problem2}
Given a smooth projective variety $X$, determine the group of isomorphism classes of autoequivalences of $D(X)$.
\end{problem}
\noindent
We denote this group by $\Auteq D(X)$.
We note that $\Auteq D(X)$ always contains the group $A(X):=(\Aut X \ltimes \Pic X) \times \ZZ$,
generated by functors of tensoring with invertible sheaves, automorphisms of $X$ and the shift functor.
When $K_X$ or $-K_X$ is ample, it is shown that $\Auteq D(X)\cong A(X)$ in \cite{Bondal:rvfd}. 
When $X$ is an abelian variety, Orlov solves Problem \ref{problem2} in \cite{Orlov:dcav}.
In this case, $\Auteq D(X)$ is strictly larger than $A(X)$.

The twist functors along spherical objects are autoequivalences of another kind that are not in $A(X)$.
Seidel and Thomas \cite{Seidel:bga} introduced them, expecting that they should correspond 
via Kontsevich's homological mirror conjecture to the generalized Dehn twists along Lagrangian spheres.
These functors play an essential role in this paper and we recall the
definition.

For an object $\mc{P}\in D(X\times Y)$, an \emph{integral functor} 
$$
\Phi ^{\mc{P}}_{X\to Y}:D(X)\to D(Y)
$$ is defined by 
$$
\Phi ^{\mc{P}}_{X\to Y}(-)=\mathbf{R}\pi_{Y*}(\mc{P}\ltensor \mathbf{L}\pi^*_X(-)),
$$
where $\pi_X:X\times Y\to X$ and $\pi_Y:X\times Y\to Y$ are the projections.
\begin{definition}[\cite{Seidel:bga}]
\begin{enumerate}
\item
We say that an object $\alpha \in D(X)$ is \emph{spherical} if 
we have $\alpha \otimes \omega_{X} \cong \alpha$ and
$$
\Hom^{k}_{D(X)}(\alpha,\alpha)\cong\begin{cases}  0 & k\ne 0,\dim X\\
                                                \CC & k=0,\dim X. 
\end{cases}
$$
\item
Let $\alpha\in D(X)$ be a spherical object. 
We consider the mapping cone 
$$
\mc{C}=Cone(\pi_1^*\alpha^\vee\ltensor \pi_2^*\alpha \to \mc{O}_{\Delta})
$$
of the natural evaluation $\pi_1^*\alpha^\vee\ltensor \pi_2^*\alpha \to \mc{O}_{\Delta}$,
where $\Delta\subset X\times X$ is the diagonal, and $\pi_i$ is the projection of
$X\times X$ to the $i$-th factor. Then the integral functor 
$T_{\alpha}:=\Phi^{\mc{C}}_{X\to X}$ defines an autoequivalence of $D(X)$, 
called the \emph{twist functor} along the spherical object $\alpha$. 
\end{enumerate}
\end{definition}

Consider the derived category $D(X)$ for a smooth surface $X$.
It is natural to ask how large the subgroup of $\Auteq D(X)$ generated by
$A(X)$ and the twists along spherical objects is.
An example of a spherical object in $D(X)$ is given by a line bundle
$\mc{R}$ on a chain of $-2$-curves on $X$,
considered as a sheaf on $X$.
In this paper, we consider a chain $Z$ of $-2$-curves
on a smooth surface $X$ and study the autoequivalences
of the derived category $D_Z(X)$ of coherent sheaves on $X$ supported by $Z$.

Note that the twist functor $T_{\alpha}$ can be defined
as long as the support of  $\alpha$ is projective,
even if $X$ is not projective.
Moreover, the category $D_Z(X)$ depends only on the formal
neighborhood of $Z$ in $X$.
Thus we can assume as follows:
$$
Y=\Spec \CC[[x, y, z]]/(x^2+y^2+z^{n+1})
$$
is the $A_n$-singularity,
$$
f:X \to Y
$$
its minimal resolution and
$$
Z=f^{-1}(P)=C_1 \cup \dots \cup C_n
$$
where $P\in Y$ is the closed point.

For an autoequivalence $\Phi \in \Auteq D_Z(X)$,
we don't know if it is always isomorphic to an integral functor.
Here, an integral functor from $D_Z(X)$ to $D_Z(X)$ is defined
by an object $\mc P \in D(X\times X)$ whose support is projective
over $X$ with respect to each projection.
If an autoequivalence is given as an integral functor,
we call it a {\it Fourier-Mukai transform} ({\it FM transform}).
Let
$$
\Auteq^{\text{FM}} D_Z(X) \subset \Auteq D_Z(X)
$$
be the subgroup consisting of FM transforms.
Remark that $\Aut X\cong\Aut Y$
and $\Pic X\cong \Pic (X/Y)$ act faithfully on $D_Z(X)$
in our setting; therefore 
we see $A(X)\subset\Auteq^{\text{FM}} D_Z(X)$. 

We also define a normal subgroup
$$
N(D_Z(X)) \subset \Auteq D_Z(X)
$$
consisting of $\Phi$ with $\Phi(\alpha) \cong \alpha$ for
every object $\alpha \in D_Z(X)$.
This group is trivial if every autoequivalence is an FM transform.
We denote the dualizing sheaf on $Z$ by $\omega_Z$ and put
$$
B=\Span{T_{\owe_{C_l}(-1)}, T_{\omega_Z}\bigm| 1\le l\le n } \subset \Auteq D_Z(X).
$$
The following is a main result of this article.


\begin{theorem}\label{main result:3}
We have
$$
\Auteq D_Z(X) = \Auteq^{\text{FM}} D_Z(X) \ltimes N(D_Z(X))
$$
and
$$
\Auteq^{\text{FM}} D_Z(X)=(\Span{B, \Pic X} \rtimes \Aut X) \times \ZZ.
$$
Here $\ZZ$ is the group generated by the shift $[1]$.
\end{theorem}


\begin{remark}[see Proposition \ref{subgroups} and Corollary \ref{corollary:allspherical}]\label{root&weight}
We know more about subgroups of $\Auteq^{\text{FM}} D_Z(X)$,
that is, we have the following:
\begin{itemize}
\item $B\cap\Pic X=\Span{\otimes\mc{O}_X(C_1),\ldots,\otimes\mc{O}_X(C_n)}$.
\item $\Span{B, \Pic X} \cong B \rtimes \ZZ/(n+1)\ZZ$.
\item $B=\Span{T_\alpha \bigm| \alpha \in D_Z(X), \textit{ spherical } }$.
\end{itemize}
\end{remark}

Put $\alpha_i:=\owe_{C_i}(-1)$ $(1\le i\le n)$
and $\alpha_0:=\alpha_{n+1}:=\omega_Z$, where we consider the suffix $i$ of $\alpha _i$ modulo $n+1$ 
(that is, $\alpha_{i}=\alpha_{n+1+i}$ for all $i\in \ZZ$).
$B$ is generated by all $T_{\alpha_i}$'s by definition.
We denote by $B_k$ the subgroup of $B$ 
generated by all $T_{\alpha_i}$'s except $T_{\alpha_k}$. 
The result in \cite{Seidel:bga} implies 
that the defining relation of the group $B_k$ is as follows:
$$
\begin{cases}
T_{\alpha_i}T_{\alpha_{i+1}}T_{\alpha_i}\cong T_{\alpha_{i+1}}T_{\alpha_i}T_{\alpha_{i+1}} 
& \textit{ if } 0\le i\le n,\quad i\ne k-1,k \\
T_{\alpha_i}T_{\alpha_j}\cong T_{\alpha_j}T_{\alpha_i} & \textit{ if } i-j\ne \pm 1,0.
\end{cases}
$$
In other words, $B_k$ is the Artin group of type $A_n$ 
(or the braid group on $n+1$ strands). 
Conjecturally our group $B$ is the Artin group of type $\tilde{A}_n$.

According to Orlov's theorem \cite{Orlov:edck3}, any autoequivalence $\Phi\in \Auteq D(S)$ 
for a smooth \emph{projective} variety $S$ is isomorphic to an integral functor $\Phi^{\mc P}_{S \to S}$
for some $\mc P \in D(S \times S)$.
Using this, we obtain another main result:


\begin{theorem}\label{main result:2}
Let $S$ be a smooth projective surface of general type
whose canonical model has $A_n$-singularities at worst.
Then we have
$$
\Auteq D(S)=\Span { T_{\mc{O}_C(a)}, A(S)\bigm | \mbox{$C: -2$-curve, $a\in\ZZ$} }.
$$
\end{theorem}

In the proofs of Theorems \ref{main result:3} and \ref{main result:2},
the following proposition is essential.


\begin{keyproposition}\label{keyproposition}
For any $\Phi \in \Auteq D_Z(X)$,
there exists an integer $i$ and $\Psi \in B$
such that $\Psi\circ \Phi$ sends every skyscraper sheaf $\owe_x$ with $x \in Z$
to $\owe_y[i]$
for some $y \in Z$.
\end{keyproposition}


\paragraph{Strategy for the proof of Key Proposition.}
Our main results, Theorems \ref{main result:3} and \ref{main result:2},
follow from Key Proposition and rather a formal argument.
Here we shall explain how to prove Key Proposition
because it is essential in this article.
For $\alpha \in D_Z(X)$, let us put 
$$
l(\alpha)=\sum _{i,p}\length _{\owe_{X,\eta _i}} \cohom ^p (\alpha)_{\eta_i},
$$
where $\eta_i$ is the generic point of $C_i$.
When $\alpha$ is spherical, we can see that every cohomology sheaf 
$\cohom ^p(\alpha)$ is a pure one-dimensional $\owe _Z$-module (Corollary \ref{rigid & pure}). 
Hence if $l(\alpha)=1$, we get $\alpha\cong\owe _{C_b}(a)[i]$ for some $a,b,i\in\ZZ$.
To show Key Proposition, we first prove that for a spherical $\alpha$ with $l(\alpha)>1$, there is
an autoequivalence $\Psi\in B$ such that $l(\alpha)>l(\Psi (\alpha))$. 
Then, since $\Psi(\alpha)$ is again spherical, induction on $l(\alpha)$ yields the following:


\begin{proposition}\label{proposition:step -1 of A_n}
Let $\alpha$ be a spherical object in $D_Z(X)$. 
Then there are integers $a$, $b$ $(1\le b\le n)$ and $i$,
and there is an autoequivalence $\Psi\in B$ such that
$$
\Psi(\alpha) \cong \owe_{C_b}(a)[i].
$$
\end{proposition}
\noindent
Next step to prove Key Proposition is to show: 


\begin{proposition}\label{proposition:step -2 of A_n}
Suppose that an autoequivalence $\Phi$ of $D_Z(X)$ is given.
Then, there are integers $a$, $b$ $(1\le b\le n)$ and $i$,
and there is an autoequivalence $\Psi\in B$ such that
$$
\Psi \circ \Phi(\owe_{C_1}) \cong \owe_{C_{b}}(a)[i]
$$
and
$$
\Psi \circ \Phi(\owe_{C_1}(-1)) \cong \owe_{C_{b}}(a-1)[i].
$$
In particular, 
for any point $x\in C_1$, we can find a point $y\in C_b$ 
with $\Psi\circ\Phi (\mc{O}_x)\cong \mc{O}_y[i]$.  
\end{proposition}
\noindent
Put $\alpha=\Phi(\owe_{C_1})$ and $\beta=\Phi(\owe_{C_1}(-1))$.
By Proposition \ref{proposition:step -1 of A_n}, we may assume that $l(\alpha)=1$. To prove 
Proposition \ref{proposition:step -2 of A_n}, we show the existence of $\Psi\in B$ 
such that $l(\Psi(\alpha))=1$ and $l(\beta)>l(\Psi(\beta))$.
Then we can complete the proof by induction on $l(\beta)$.

Once we get Proposition \ref{proposition:step -2 of A_n}, 
we can rather easily show Key Proposition by induction on $n$. 


\paragraph{Construction of this article.}
In \S \ref{section:the main result},
we first demonstrate that Proposition \ref{proposition:step -2 of A_n}
implies Key Proposition.
We then prove our main results, Theorem \ref{main result:3}
and Theorem \ref{main result:2}.
The rest of this paper is devoted to showing Proposition \ref{proposition:step -2 of A_n}.

In \S \ref{section:tool}, we study spherical objects and the twist functors
for a smooth surface $X$, which play the leads in our article.
We first observe that the isomorphism class of an object $\alpha \in D(X)$
is determined by the cohomology sheaves $\cohom^i(\alpha)$ and some connecting
data $e^i(\alpha) \in \Ext^2_X(\cohom^i(\alpha), \cohom^{i-1}(\alpha))$.
Then we give a necessary and sufficient condition for $\alpha$ to be spherical
in terms of $\cohom^i(\alpha)$ and $e^i(\alpha)$.
Especially, for a chain $Z$ of $-2$-curves on $X$ and a spherical object $\alpha\in D_Z(X)$,
we see that $\bigoplus_p\mc{H}^p(\alpha)$ is a rigid $\mc{O}_Z$-module, 
pure of dimension $1$ (Corollary \ref{rigid & pure}). 
This result, combined with Lemma \ref{sheaf on A_n} on pure sheaves on $Z$,
enables explicit computations in the latter sections.

In \S \ref{section:a_1}, as a first step, we consider the $A_1$ cases of 
Proposition \ref{proposition:step -1 of A_n} and Proposition \ref{proposition:step -2 of A_n}.
We show Proposition \ref{proposition:step -1 of A_n} in \S \ref{section:preliminary}
and Proposition \ref{proposition:step -2 of A_n} in \S \ref{section:a_n} respectively.
In \S \ref{section:preliminary}, we compute $l(\Psi(\alpha))-l(\alpha)$ for various $\Psi$'s in $B$
by using results from \S \ref{section:tool} and Lemma \ref{sheaf on A_n}.
We use similar methods in \S \ref{section:a_n} and find $\Psi$ 
in the statement of Proposition \ref{proposition:step -2 of A_n} 
via case-by-case arguments.


\paragraph{Notation and Convention.}
We work over the complex number field $\CC$. Let $X$ be an algebraic variety and 
$Z$ a closed subset of $X$.
$D_Z(X)$ denotes the full subcategory of $D(X)$ consisting 
of objects supported on $Z$.
Here, the support of an object of $D_Z(X)$ is, by definition, the union of the supports
of its cohomology sheaves.
It is known that $D_Z(X)$ is naturally equivalent 
to the bounded derived category of coherent sheaves on $X$, supported on $Z$
(see \cite[Proposition 1.7.11]{Kashiwara:som}).  
When we write $D_Z(X)$ for a closed subscheme $Z$ of $X$, 
we forget the scheme structure of $Z$ and regard it as a closed subset of $X$.  
Let $D_{c}(X)$ denote the derived category of ``compactly supported'' coherent sheaves on $X$,
i.e. coherent sheaves whose supports are proper over $\CC$.

Next let $Z=C_1\cup \cdots \cup C_n$ be a chain of $-2$-curves on a smooth surface $X$.
Namely, each $C_l$ is a smooth rational curve with $C_l ^2=-2$ and 
$$
C_l\cdot C_m = 
\begin{cases}
                1 & |l-m|=1 \\
                0 & |l-m|\ge 2. 
\end{cases}
$$
We regard $Z$ as a closed subscheme of $X$ with respect to the reduced induced structure.
For a coherent sheaf
$\mc{R}$ on $Z$, we denote by $\deg_{C_l}\mc{R}$
the degree of the restriction $\mc{R}|_{C_l}$ on $C_l\cong\PP^1$.
We denote by
$$
\mc{R}_0=\mc{O}_{C_1\cup \cdots \cup C_n}(a_1, \dots, a_n)
$$
the line bundle (or $\owe_Z$-invertible sheaf) on $Z$ such that $\deg_{C_l}\mc{R}_0=a_l$ for all $l$.
When we write $*$ instead of $a_l$, we don't specify the degree at $C_l$.
For instance, when we put 
$$
\mathcal{R}_1 = \owe_{C_1\cup C_2 \cup C_3}(a, b, *),
$$
this means that $\mc{R}_1$ is a line bundle on $C_1 \cup C_2 \cup C_3$
such that $\deg_{C_1}\mc{R}_1=a$, $\deg_{C_2}\mc{R}_2=b$ and $\deg_{C_3}\mc{R}_1$ arbitrary.
The expression
$$
\mc{R}_2=\mc{O}_{C_1\cup \cdots}(a,*)
$$
means that there exists $t \ge 2$ with $\mc{R}_2= \owe_{C_1 \cup C_2 \cup \cdots \cup C_t}(a,*, \dots, *)$.
Note that the support of $\mc{R}_2$ is strictly larger than $C_1$.
We often use figures
\[ 
\xymatrix@R=1ex@M=0ex{  & C_1 & C_2 & C_3 \\
   \mc{R}_1:  & 
\setlength{\unitlength}{1ex}
\begin{picture}(2, 2)(-1, -1)
\put(0,0){\circle{2}}
\put(0,0){\makebox(0,0){\tiny$a$}}
\end{picture} 
\ar@{-}[r] &
\setlength{\unitlength}{1ex}
\begin{picture}(2, 2)(-1, -1)
\put(0,0){\circle{2}}
\put(0,0){\makebox(0,0){\tiny$b$}}
\end{picture}
\ar@{-}[r] &{\no} & \\
\mc{R}_2:  &
\setlength{\unitlength}{1ex}
\begin{picture}(2, 2)(-1, -1)
\put(0,0){\circle{2}}
\put(0,0){\makebox(0,0){\tiny$a$}}
\end{picture}
\ar@{-}[r] & & & }
\]
\noindent
to define $\mc{R}_1, \mc{R}_2$ above.
We use a dotted line
\[ 
\xymatrix@R=1ex@M=0ex{  &C_1  & \\
          \mc{R}_3:    &\setlength{\unitlength}{1ex}
\begin{picture}(2, 2)(-1, -1)
\put(0,0){\circle{2}}
\put(0,0){\makebox(0,0){\tiny$a$}}
\end{picture}
\ar@{--}[r]& &}
\]
to indicate that $\mc{R}_3$ is either $\mc{O}_{C_1}(a)$ or $\mc{O}_{C_1\cup\cdots}(a,*)$.

For an object $\alpha \in D_Z(X)$, we put
$$
l(\alpha)=\sum _{i,p}\length _{\owe_{X,\eta _i}} \cohom ^p (\alpha)_{\eta_i},
$$
where $\owe_{X,\eta _i}$ is the local ring of $X$ at the generic point $\eta _i$ of $C_i$,
$\cohom ^p (\alpha)_{\eta_i}$ is the stalk over $\eta_i$
and $\length _{\owe_{X,\eta _i}}$ measures the length over $\owe_{X,\eta _i}$.

Throughout this paper, a \emph{point} on a variety always means a \emph{$\CC$-valued point} unless otherwise specified.   
For a point $x$ on a variety $X$, we denote the structure sheaf of $x$ by $\mc O_x$.
We regard it as a skyscraper sheaf on $X$.

\paragraph{Acknowledgements.} The authors thank Alastair Craw for his helpful comments and remarks.
The first author also thanks Yujiro Kawamata for valuable discussions.

 \section{Main results}\label{section:the main result}

In this section, we first show that Key Proposition
follows from Proposition \ref{proposition:step -2 of A_n}
that will be shown in \S \ref{section:a_n}.
As its application, we prove our main results, 
Theorem \ref{main result:3} and Theorem \ref{main result:2}.
In the proof of Theorem \ref{main result:3},
we use the facts that
$B\cap\Aut X=\{id\}$ and that $B$ is a normal subgroup of $\Span{B, A(X)}$,
which will be explained in Remark \ref{remark:normal sub'}.

 This section is logically the final part of this article. Therefore 
we do not use the results in \S 2 afterwards.


\subsection{Proof of Key Proposition.}
Let us first show the following claim. 
\begin{claim}\label{claim:silly}
Assume that $\Phi(\owe_{C_1})\cong \owe_{C_l}(a)$ and $\Phi(\owe_{C_1}(-1)) \cong \owe_{C_l}(a-1)$
for some $l$.
Then $l=1$ or $n$.
\end{claim}
\begin{proof}
The assumption implies that a closed point $x \in C_1$ corresponds bijectively to $y \in C_l$
such that $\Phi(\owe_x) \cong \owe_y$.
If $1<l<n$, there are points $y_0, y_1$ such that $C_l\cap C_{l+1}=\{y_0\}$ and $C_{l-1}\cap C_{l}=\{y_1\}$.
Let $x_0, x_1 \in C_1$ be the points with $\Phi(\owe_{x_0}) \cong \owe_{y_0}$
and $\Phi(\owe_{x_1}) \cong \owe_{y_1}$.
Then $x_0$ is contained in $\Supp \Phi^{-1}(\owe_{C_{l+1}}) \cap C_1$.
Since $\Supp \Phi^{-1}(\owe_{C_{l+1}})$ is connected and does not contain $C_1$,
$x_0$ is the intersection point of $C_1$ and $C_2$.
By the same argument, we obtain $x_0=x_1$, which is absurd.
\end{proof}
We want to show that there is an autoequivalence
$$
\Psi\in\Span{[i],B\bigm| i\in\ZZ}
$$
such that for any point $x\in Z$, we can find a point $y\in Z$ with $\Psi\circ\Phi (\mc{O}_x)\cong \mc{O}_y$.

The assertion for the case $n=1$ follows directly from Proposition \ref{proposition:step -2 of A_n}, 
and hence we may assume $n>1$.
Utilizing Proposition \ref{proposition:step -2 of A_n} and Claim \ref{claim:silly}, we obtain
an autoequivalence
$$
\Psi_1\in\Span{[i],B \bigm| i\in\ZZ}
$$
such that for any point $x\in C_1$, we have a point $y\in C_l$ with $\Psi _1\circ \Phi(\mc{O}_x)\cong\mc{O}_y$.
Here, $l=1$ or $n$ and we consider the case $l=n$, the other case is similar.
Put $Z_1=\sum_{k=2}^n C_k$ and $Z_2 =\sum_{k=1}^{n-1} C_k$.
Then we can see that $\Psi _1\circ \Phi$ induces an equivalence $D_{Z_1}(X) \cong D_{Z_2}(X)$.
By the induction hypothesis, there is
$$
\Psi_2 \in \Span{T_{\owe_{C_l}(a)}\bigm|  a\in \ZZ, 1\le l\le n-1 }
$$
such that $\Psi:=\Psi_2 \circ \Psi_1$ has the desired property, and we finish the proof of 
Key Proposition.
\qed\\

Let $\iota \in \Aut Y(\cong \Aut X)$ be an involution such that $\iota(C_i)=C_{n-i+1}$ for curves $C_i$.
The above proof also supplies the following:

 
\begin{corollary}\label{corollary:1}
For any $\Phi \in \Auteq D_Z(X)$, there is $\Psi \in \Span{B, \iota^*, [i] \bigm| i\in\ZZ}$
such that
$$
\Psi \circ \Phi (\mc{R}) \cong \mc{R}
$$
for every line bundle $\mc{R}$ on any subchain of $Z$. 
\end{corollary}


\subsection{Proof of Theorem \ref{main result:3}.}
First of all, we show the equality 
$$
\Auteq^{\text{FM}} D_Z(X)=(\Span{B, \Pic X} \rtimes \Aut X) \times \ZZ.
$$
Note that $B\cap\Aut X=\{id\}$ and $B$ is a normal subgroup of $\Span{B,A(X)}$ by Remark \ref{remark:normal sub'}.
Therefore it suffices to show that $\Phi$ belongs to $\Span{B,A(X)}$ for any $\Phi\in\Auteq^{\text{FM}} D_Z(X)$.
Key Proposition implies that there are $\Psi\in B$ and an integer $i$
such that for any point $x\in Z$, 
we have $\Psi\circ \Phi (\mc{O}_x)\cong\mc{O}_y[i]$ for some point $y\in Z$.
Then Lemma \ref{A(X)} assures that $\Psi\circ \Phi \in A(X)$,
and thus we get the conclusion.


\begin{lemma}\label{A(X)}(\cite[3.3]{Bridgeland:fmqv})
Suppose an autoequivalence $\Phi\in\Auteq ^{FM}D(X)$ for an algebraic variety
$X$ satisfies the following: for any point $x\in X$, there is a point $y\in X$ such that
$\Phi (\mc{O}_x)\cong \mc{O}_y$. Then $\Phi\in \Pic X \rtimes \Aut X$.
\end{lemma}

Next we prove
$$
\Auteq D_Z(X) = \Auteq^{\text{FM}} D_Z(X) \ltimes N(D_Z(X))
$$
by using the {\it McKay correspondence}.
Recall that $Y$ is isomorphic to (the germ of) a quotient singularity $\CC^2/G$,
where $G \subset \operatorname{SL}(2, \CC)$ is a finite subgroup;
the $A_n$-singularity corresponds to the case $G \cong \ZZ/(n+1)\ZZ$.
Let $\coh^G(\CC^2)$ be the abelian category of
$G$-equivariant coherent sheaves on $\CC^2$
and $D^G(\CC^2)$ its bounded derived category.
The McKay correspondence \cite{kapranov:kdh}
establishes an equivalence from the derived category of the minimal resolution
of $\CC^2/G$ to $D^G(\CC^2)$, which is an FM transform.
This induces an equivalence from $D_Z(X)$ to the full subcategory
$D^G_{\{0\}}(\CC^2)$ of objects supported on the set $\{0\}$.
Especially, it sends $\owe_{C_i}(-1) \in D_Z(X)$ to $\rho_i \otimes \owe_0 \in D^G_{\{0\}}(\CC^2)$,
where $\rho_1, \dots, \rho_n$ are the non-trivial irreducible representations of $G$.
Moreover, $\omega_Z$ corresponds to $\rho_0 \otimes \owe_0[-1]$ where $\rho_0$
is the the trivial representation of $G$.
Thus, an autoequivalence of $D_Z(X)$ which fixes $\omega_Z$ and $\owe_{C_i}(-1)$
for $i =1, \dots, n$
corresponds to an autoequivalence of $D^G_{\{0\}}(\CC^2)$ which fixes
$\rho_0 \otimes \owe_0, \dots, \rho_n \otimes \owe_0$. 
Recall that we have a natural isomorphism $\Aut X\cong\Aut Y$; via this isomorphism 
$\Aut X$ acts both on $D_Z(X)$ and on $\eqcat$ preserving the McKay correspondence.


\begin{proposition}\label{propostion:McKay}
Let $\Phi$ be an autoequivalence of $\eqcat$
satisfying $\Phi(\rho_i \otimes \owe_0) \cong \rho_i \otimes \owe_0$
for all irreducible representations $\rho_i$ of $G$.
Then there is an automorphism $\sigma\in \Aut Y$
such that
$$
\Phi(\alpha) \cong \sigma^* \alpha
$$
for all $\alpha \in \eqcat$.
\end{proposition}

\begin{proof}
Since any sheaf $\mc{F} \in \coh^G_{\{0\}}(\CC^2)$ is a successive extension of sheaves $\rho_i \otimes \owe_0$,
it follows from the assumption that $\Phi(\mc{F})$ is also a sheaf.
Moreover, $\Phi$ restricted to $\coh^G_{\{0\}}(\CC^2)$ is an exact functor of abelian categories.
Let $R$ be the affine coordinate ring of $\CC^2$
with maximal ideal $m$ of the origin.
We denote by $\widehat{R}$ the completion of $R$ with respect to $m$.

\begin{claim}\label{claim:R/m^l}
We have $\Phi(\rho_i \otimes R/m^l) \cong \rho_i \otimes R/m^l$ for
all irreducible representations $\rho_i$ and for all positive integers $l$.
\end{claim}

\begin{proof}
We prove the claim by induction on $l$.
The case $l=1$ is included in the assumption.
Assume $l>1$ and consider the short exact sequence
$$
0 \to \rho_i \otimes m^{l-1}/m^l \to \rho_i \otimes R/m^l \to \rho_i \otimes R/m^{l-1} \to 0.
$$
Since the equivalence $\Phi$ sends a sheaf to a sheaf,
the following is also an exact sequence of sheaves:
$$
0 \to \Phi(\rho_i \otimes m^{l-1}/m^l) \to \Phi(\rho_i \otimes R/m^l) \to \Phi(\rho_i \otimes R/m^{l-1}) \to 0.
$$
Here, we have $\Phi(\rho_i \otimes m^{l-1}/m^l) \cong \rho_i \otimes m^{l-1}/m^l$ since
$m^{l-1}/m^l$ is a direct sum of sheaves $\rho_j \otimes \owe_0$,
and $\Phi(\rho_i \otimes R/m^{l-1}) \cong \rho_i \otimes R/m^{l-1}$ by the induction hypothesis.
Therefore, the claim follows from the following lemma.
\end{proof}


\begin{lemma}
Let
$$
0 \to \rho_i \otimes m^{l-1}/m^l \to \mc{F} \to \rho_i \otimes R/m^{l-1} \to 0
$$
be the extension corresponding to a class $e \in \GExt^1_{\CC^2}(\rho_i \otimes R/m^{l-1}, \rho_i \otimes m^{l-1}/m^l)$.
Then, $\mc{F} \cong \rho_i \otimes R/m^l$ if and only if
$\phi \circ e \ne 0$ in $\GExt^1_{\CC^2}(\rho_i \otimes R/m^{l-1}, \rho_j \otimes \owe_0)$
for any $j$ and for any surjection $\phi: \rho_i \otimes m^{l-1}/m^l \to \rho_j \otimes \owe_0$.
\end{lemma}

\begin{proof}
The \lq only if' part is obvious.
Let $\mc{F}$ be an extension with the above property.
Lift $\rho_i\otimes 1 \subset \rho_i \otimes R/m^{l-1}$ to a $G$-invariant vector subspace $V \cong \rho_i$ of $\mc{F}$.
The assumption on $e$ implies that $V$ generates $\mc{F}$ as an $R$-module.
Therefore, $\mc{F}$ is of the form $\rho_i \otimes R/J$ for a $G$-invariant $R$-submodule $J$ of $\rho_i \otimes R$.
Since $\mc{F}$ fits into the above extension, $J$ must coincide with $\rho_i \otimes m^l$.
\end{proof}

\noindent
We denote by $j: \mc{C} \hookrightarrow \eqcat$ the full subcategory whose objects are
sheaves $\rho_i \otimes R/m^l$
where $i$ and $l$ vary.

\begin{claim}
There exists an automorphism $\sigma\in \Aut Y$ with an isomorphism $\phi: \sigma^* \circ j \cong \Phi \circ j$.
\end{claim}

\begin{proof}
$\Phi$ induces an isomorphism (of $\CC$-algebras)
$$
\sigma_l:\GHom_{\CC^2}(R/m^l, R/m^l) \cong \GHom_{\CC^2}(\Phi(R/m^l), \Phi(R/m^l)).
$$
\noindent
By Claim \ref{claim:R/m^l}, the right hand side is isomorphic to $(R/m^l)^G$
and this isomorphism does not depend on the choice of the isomorphism in Claim \ref{claim:R/m^l}.
Hence $\sigma_l$ is a $\CC$-algebra
automorphism of $(R/m^l)^G$.
Put
$$
\sigma = \varprojlim_{l} \sigma_l \in \Aut Y.
$$
By replacing $\Phi$ with $(\sigma^*)^{-1} \circ \Phi$, we may assume that $\sigma$ is the identity.
We choose isomorphisms $\phi_l^0: R/m^l \overset{\sim}{\to} \Phi(R/m^l)$ such that
$$
\xymatrix{
R/m^{l+1} \ar[d]^{\cong}_{\phi_{l+1}^0} \ar@{>>}[r]^{p_l} & R/m^l \ar[d]^{\cong}_{\phi^0_l} \\
\Phi(R/m^{l+1}) \ar@{>>}[r]^{\Phi(p_l)}& \Phi(R/m^l)
}
$$
commutes where $p_l$ is the projection.
We see that $\Phi(f)\circ \phi^0_l = \phi^0_l \circ f$ for any $G$-equivariant morphism $f: R/m^l \to R/m^l$
since $f$ is the multiplication by an element of $(R/m^l)^G$ and since $\sigma_l$ is the identity.

For $i \ne 0$, we first choose isomorphisms 
$\psi_l^i: \rho_i \otimes R/m^l \cong \Phi(\rho_i \otimes R/m^l)$ such that 
$\psi_l^i \circ (1_{\rho_i} \otimes p_l) = \Phi(1_{\rho_i} \otimes p_l) \circ \psi_{l+1}^i$.
For an element $a \in (\rho_i \otimes R/m^l)^G$,
denote by $m_a: R/m^l \to \rho_i \otimes R/m^l$ the multiplication by $a$.
Then $(\psi_l^i)^{-1} \circ \Phi(m_a) \circ \phi^0_l$ is also a morphism
from $R/m^l$ to $\rho_i \otimes R/m^l$ and hence is the multiplication
by an element $\xi_l(a)$ of $(\rho_i \otimes R/m^l)^G$.
Here, $\xi_l$ is an automorphism of $(\rho_i \otimes R/m^l)^G$
as an additive group.
Moreover, for any $b \in (R/m^l)^G$, the relation $m_{ba}=m_a \circ m_b$
implies that $\xi_l$ is $(R/m^l)^G$-linear.
Furthermore, $\xi_{l+1}$ induces $\xi_l$ on $(\rho_i \otimes R/m^l)^G$.
Therefore, we can define $\xi= \varprojlim_{l} \xi_l$
which is a $\widehat{R}^G$-module automorphism of $(\rho_i \otimes \widehat{R})^G$.
Since
$$
\rho_i \otimes \widehat{R} \cong \left( (\rho_i \otimes \widehat{R})^G \otimes_{(\widehat{R})^G}
\widehat{R}\right)^{\vee\vee}
$$
(\cite{Esnault:rmqs}, see also \cite[Theorem 12]{Riemenschneider:srtm}),
$\xi$ gives rise to automorphisms $\tilde \xi$ of $\rho_i \otimes \widehat{R}$
and therefore we obtain an automorphism $\tilde \xi_l$ of $\rho_i \otimes R/m^l$ for any $l$
which coincide with $\xi_l$ on $(\rho_i \otimes R/m^l)^G$.
Put
$$
\phi_l^i:=\psi_l^i \circ \tilde \xi_l.
$$
Then for any $a\in (\rho_i \otimes R/m^l)^G$, we have $\tilde \xi_l\circ m_a=m_{\xi_l(a)}$, hence the diagram
$$
\xymatrix{
R/m ^l \ar[r]^{m_a} \ar[d]^{\cong}_{\phi_l^0}    & \rho_i \otimes R/m^l \ar[r]^{\tilde \xi_l}
& \rho_i \otimes R/m^l \ar[d]^{\psi_l^i}_{\cong} \\
 \Phi(\rho_i \otimes R/m^l) \ar[rr]^{\Phi(m_a)}  &                                               &  \Phi(\rho_i \otimes R/m^l)
}
$$
is commutative. Then we obtain
\begin{equation}\label{equation:0icommutes}
\phi_l^i \circ m_a = \Phi(m_a) \circ \phi^0_l.
\end{equation}

Finally, we consider a $G$-equivariant morphism $f: \rho_i \otimes R/m^k \to \rho_j \otimes R/m^l$ for
arbitrary $i, j,k,l$ and show that
\begin{equation}\label{equation:ijcommutes}
\xymatrix{
\rho_i \otimes R/m^k \ar[r]^f \ar[d]^{\phi_k^i} & \rho_j \otimes R/m^l \ar[d]^{\phi_l^j} \\
\Phi(\rho_i \otimes R/m^k) \ar[r]^{\Phi(f)} & \Phi(\rho_j \otimes R/m^l)
}
\end{equation}
commutes.
When $k=l$, we write $\zeta^{ji}_l (f)=(\phi_l^j)^{-1} \circ \Phi(f) \circ \phi_l^i$
and put $\zeta^{ji}=\varprojlim_{l} \zeta^{ji}_l$.
Then $\zeta^{ji}$ is a $\widehat{R}^G$-automorphism of
$$
\left(\Hom_{\widehat{R}}(\rho_i \otimes \widehat{R}, \rho_j \otimes \widehat{R})\right)^G
\cong \Hom_{\widehat{R}^G}((\rho_i \otimes \widehat{R})^G, (\rho_j \otimes \widehat{R})^G).
$$
Take $f \in \Hom_{\widehat{R}^G}((\rho_i \otimes \widehat{R})^G, (\rho_j \otimes \widehat{R})^G)$ and
$g \in \Hom_{\widehat{R}^G}(\widehat{R}^G, (\rho_i \otimes \widehat{R})^G)$.
Then we have $\zeta^{j0}(f\circ g)=\zeta^{ji} (f) \circ \zeta^{i0} (g)$ by the definition of $\zeta^{ji}$'s.
\eqref{equation:0icommutes} shows that $\zeta^{i0}(g)=g$ and $\zeta^{j0}(f \circ g)=f\circ g$.
Since $g$ is arbitrary, these equalities imply that $\zeta^{ji}(f)=f$ and hence the commutativity
of \eqref{equation:ijcommutes} in the case $k=l$.
If $k>l$, then $f$ factors through $\rho_i \otimes R/m^l$
and if $k<l$ then $f$ can be composed with the surjection $\rho_i \otimes R/m^l \to \rho_i \otimes R/m^k$.
In this way, we obtain the commutativity of \eqref{equation:ijcommutes}.
\end{proof}


\begin{claim}\label{claim:sheaf}
Let $j': \coh^G_{\{0\}}(\CC^2) \hookrightarrow \eqcat$
be the natural embedding.
Then we have an isomorphism $\phi: \sigma^* \circ j' \cong \Phi \circ j'$.
Moreover, for $\mc{F} \in \coh^G_{\{0\}}(\CC^2)$, let us define
$\phi_{\mc{F}[n]}: \mc{F}[n] \to \Phi(\mc{F}[n])$ by $\phi_{\mc{F}[n]}=\phi_{\mc{F}}[n]$.
Then, these isomorphisms commute with $\Hom$'s between
 shifts of sheaves: $\mc{F}[n]$ and $\mc{G}[m]$.
\end{claim}

\begin{proof}
As in the proof of the previous claim, we may assume $\sigma$ is the identity.
For $\mc{F} \in \coh^G_{\{0\}}(\CC^2)$, we can take a presentation
$$
\mc{E}_1 \to \mc{E}_0 \to \mc{F} \to 0
$$
where $\mc{E}_0$ and $\mc{E}_1$ are direct sums of sheaves in $\mc{C}$.
Then, the proof is similar to that in \cite[2.16.1 -- 2.16.4]{Orlov:edck3}.
\end{proof}
\noindent
Now we give a proof of the proposition. We may assume $\sigma$ is the identity
by replacing $\Phi$ with $(\sigma^*)^{-1} \circ \Phi$.
Let $\alpha\ne 0$ be an object of $\eqcat$.
$\alpha=\alpha^{\bullet}$ is a bounded complex over $\coh^G_{\{0\}}(\CC^2)$.
Let $p$ and $q$ be the minimum and the maximum of $i$
with $\alpha^i \ne 0$,
and denote by $v$ the natural morphism $\alpha^{q}[-q] \to \alpha$.
We show by induction on $q-p$ that there is an isomorphism $\phi_{\alpha} : \alpha \to \Phi(\alpha)$
such that $\phi_{\alpha} \circ v = \Phi(v) \circ \phi_{\alpha^{q}}[-q]$.
Let $\beta=\beta^{\bullet}$ be an object such that
$$
\beta^i = \begin{cases} \alpha^i &(i\ne q) \\
0&(i=q) \end{cases}
$$
with the same differentials (except for $d^{q-1}:\beta^{q-1} \to \beta^q$) as $\alpha$.
Then $\beta$ fits into a distinguished triangle
$$
\alpha^{q}[-q] \to \alpha \to \beta \overset{t}{\to} \alpha^{q}[-q+1].
$$
By the induction hypothesis, we have an isomorphism
$\phi_{\beta}: \beta \to \Phi(\beta)$ such that
$\phi_{\beta} \circ u = \Phi(u) \circ \phi_{\beta^{q-1}}[-q+1]$
where $u: \beta^{q-1}[-q+1] \to \beta$ is the natural morphism.
For the existence of $\phi_{\alpha}$ with the prescribed property, it is enough to show
$\phi_{\alpha^q}[-q+1] \circ t = \Phi(t) \circ \phi_{\beta}$.
Consider the following diagram:
$$
\xymatrix{
\beta^{q-1}[-q+1] \ar[r]^(0.7){u} \ar[d]^{\phi_{\beta^{q-1}}[-q+1]} & \beta \ar[r]^(0.3){t} \ar[d]^{\phi_{\beta}}
 &\alpha^q[-q+1] \ar[d]^{\phi_{\alpha^q}[-q+1]} \\
 \Phi(\beta^{q-1})[-q+1] \ar[r]^(0.7){\Phi(u)} & \Phi(\beta) \ar[r]^(0.3){\Phi(t)} &
 \Phi(\alpha^q)[-q+1].
}
$$
Here the left square is commutative by virtue of the property of $\phi_{\beta}$
and the whole square is commutative by Claim \ref{claim:sheaf}.
Thus we obtain
\begin{equation}\label{equation:composeu}
\phi_{\alpha^q}[-q+1] \circ t \circ u= \Phi(t) \circ \phi_{\beta} \circ u.
\end{equation}
If we consider the object $\gamma$ in a distinguished triangle
$$
\beta^{q-1}[-q+1] \overset{u}{\to} \beta \to \gamma \to \beta^{q-1}[-q+2],
$$
then we see $\Hom(\gamma, \Phi(\alpha^q)[-q+1]) \cong \Hom(\gamma, \alpha^q[-q+1])=0$
and therefore
$$
u^*: \Hom(\beta, \Phi(\alpha^q)[-q+1]) \to \Hom(\beta^{q-1}[-q+1], \Phi(\alpha^q)[-q+1])
$$
is injective.
Thus we can remove \lq $\circ u$' from \eqref{equation:composeu} as desired.
\end{proof}

\noindent
We apply the above proposition to $D_Z(X)$ via the McKay correspondence.
Assume $\Phi \in \Auteq D_Z(X)$ is given.
>From Corollary \ref{corollary:1} and Proposition \ref{propostion:McKay}, 
we obtain an FM transform $\Psi \in \Span{B, \Aut X, [i] \bigm| i\in\ZZ}$ such that $\Psi \circ \Phi \in N(D_Z(X))$.

On the other hand, Lemma \ref{A(X)} implies that an autoequivalence
$$ 
\Phi \in \Auteq^{\text{FM}} D_Z(X) \cap N(D_Z(X))
$$
is induced by an automorphism $\sigma$ of $X$ such that
$\sigma(x)=x$ for all $x \in Z$.
Moreover, we have $\sigma^* \mc F \cong \mc F$ for any coherent sheaf $\mc F$
on $X$ supported by $\{x\} \subseteq Z$ and this implies that the 
automorphism
of the two-dimensional regular local ring $\owe_{X,x}$ induced by
$\sigma$ is the identity.
Consequently, $\sigma$ and hence $\Phi$ are the identity.
Now we obtain the splitting
$$
\Auteq D_Z(X)= \Auteq^{\text{FM}} D_Z(X) \ltimes N(D_Z(X)),
$$
which completes the proof of Theorem \ref{main result:3}.
\qed


\subsection{Proof of Theorem \ref{main result:2}.}
Let $f:S\to S_0$ be a composite of blowing-ups along a point and $S_0$ the minimal model of $S$.
\begin{claim}\label{KSC}
Let $C$ be an irreducible curve on $S$.
\begin{enumerate}
\item If $K_S\cdot C=0$, then $C$ is a $-2$-curve. Assume furthermore that $\Exc f\cap C\ne\emptyset$.
 Then $C\subset \Exc f$.
\item If $K_S\cdot C<0$, then $C$ is a $-1$-curve with $C\subset \Exc f$.
\end{enumerate}
\end{claim}
\begin{proof}
Put $K_S=f^*K_{S_0}+\sum a_iE_i$, where $E_i$'s are the components of $\Exc f$ and $a_i\in\ZZ _{>0}$.
Assume that $K_S\cdot C\le 0$. Then we have $0\ge K_S\cdot C\ge \sum a_iE_i\cdot C$, and hence
$C\cap \Exc f=\emptyset$ or $C=E_i$ for some $i$. In the former case, we get $K_S\cdot C=K_{S_0}\cdot f(C)=0$,
in particular $C$ is a $-2$-curve. If $K_S\cdot C<0$, then the latter case occurs
and we have $C^2<0$. Therefore we obtain $K_S\cdot C=-1$.
\end{proof}
Put $f=\varphi_1\circ\cdots\circ\varphi_n$ and $S=S_n$, where $\varphi_k:S_k\to S_{k-1}$ is the blow-up along
a point $x_{k-1}\in S_{k-1}$.
\begin{claim}\label{-1 curve}
Let $C$ be a $-1$-curve on $S$. If some $-2$-curve $C_1$ meets $C$, no other $-2$-curves meet $C$.
\end{claim}
\begin{proof}
We may assume that $\varphi_n$ contracts $C$, since $C\subset \Exc f$. For a contradiction, suppose that
there are two $-2$-curves $C_1, C_2$ such that both of them meet $C$. By Claim \ref{KSC},
$\varphi_n (C_1)$ and $\varphi_n (C_2)$ are $-1$-curves on $S_{n-1}$.
$\varphi_n (C_1)\cap \varphi_n (C_2)\ne\emptyset$ yields a contradiction with $C_1, C_2\subset \Exc f$.
\end{proof}
\begin{claim}\label{-2 curve}
Any connected component of the union of all $-2$-curves on $S_k$ $(0\le k\le n)$ forms a chain.
\end{claim}
\begin{proof}
We show the claim by induction on $k$. Note that the claim holds for $k=0$ by the assumption
of Theorem \ref{main result:2}.
Suppose that the claim is true for $S_k$. If there are no $-1$-curves passing through $x_k$,
we have $C\cap\Exc \varphi _{k+1}=\emptyset$ for any $-2$-curve $C$ on $S_{k+1}$.
Then the claim is true for $S_{k+1}$. If there is a $-1$-curve $C$ passing through $x_k$,
no other $-1$-curves pass through $x_k$ by Claim \ref{-1 curve} for $S_{k+1}$. Claim \ref{-1 curve}
for $S_{k}$ says that at most one $-2$-curve meets $C$. Now we get the conclusion by the induction assumption.
\end{proof}
Suppose that $\Phi\in \Auteq D(S)$ is given. Then Orlov's result \cite{Orlov:edck3} assures that
there is an object $\mc{P}\in D(S\times S)$ such that $\Phi\cong\Phi ^{\mc{P}}$.
By the proof of  \cite[Theorem 2.3]{Kawamata:deke}, we have a projective surface $Z\subset \Supp \mc{P}$ such that
$p_i|_Z:Z\to S$ $(i=1,2)$ is an isomorphism. Here $p_i$'s are the projections $S\times S\to S$.
Put  $q:=p_2|_Z\circ (p_1|_Z)^{-1}$.

Because $S$ is of general type and $\Phi (\mc{O}_x)\otimes \omega _S\cong \Phi (\mc{O}_x)$ 
for any $x\in S$ (see \cite[Theorem 2.7]{Bridgeland:csed}), we have $\dim \Phi (\mc{O}_x)\le 1$.
Assume that $\dim \Phi (\mc{O}_{x_0})= 1$ for some $x_0\in S$.
Then because $K_S\cdot C=0$ for any one-dimensional irreducible component $C$ of $\Supp \Phi (\mc{O}_{x_0})$
by the proof of \cite[Theorem 2.3]{Kawamata:deke}, Claim \ref{KSC} implies
that there is a $-2$-curve $C$ on $S$ such that $C\subset \Supp \Phi (\mc{O}_{x_0})$.
Since $q(x_0)\in C$, there is a $-2$-curve $C'$ such that $x_0\in C'$.
Therefore we can conclude that if a point $x\in S$ is not contained in any $-2$-curves,
we have $\Supp \Phi (\mc{O}_{x_0})=q(x)$. Moreover the proof of  \cite[Proposition 3.1]{Bridgeland:csed}
deduces that $\Phi (\mc{O}_{x})\cong\mc{O}_{q(x)}[i]$ for some $i\in\ZZ$.
Here the choice of $i$ is independent of the choice of $x$.

Let $\{ Z_j\}_j$ be the set of chains of $-2$-curves on $S$.
Take a point $x\in S\backslash \coprod_j Z_j$. 
Then we have $\Phi (\mc O_x)\cong \mc{O}_{q(x)}[i]$, and in particular
$q^*\circ\Phi (\mc O_x)\cong \mc O_x[i]$ for any $x\in S\backslash \coprod_j Z_j$. 
Therefore $q^*\circ\Phi$ preserves $D_{Z_j}(S)$ for each $j$.
Now Key Proposition and Lemma \ref{A(X)} complete the proof.
\qed

\section{Spherical objects and twist functors for the derived categories 
         of smooth surfaces}\label{section:tool}

This section provides technical tools used in the proofs of 
Proposition \ref{proposition:step -1 of A_n} and Proposition \ref{proposition:step -2 of A_n}.
In \S\ref{section:spectral}, we recall two kinds of spectral sequences;
their $d_2$-maps are determined by some {\it connecting data} $e^i(\alpha)$.
Then we see in \S\ref{section:reconstruction}
that the isomorphism class of an object $\alpha \in D(X)$, $X$ a smooth surface,
is determined by the cohomology sheaves $\cohom^i(\alpha)$ and the connecting
data $e^i(\alpha)$.
In \S\ref{section:spherical}, we give a necessary and sufficient condition for $\alpha$ to be spherical
in terms of $\cohom^i(\alpha)$ and $e^i(\alpha)$.
In \S\ref{section: Twist functors},
we summarize properties of twist functors and then do some computations.
We consider the group $B$ and its relation with $\Pic X$ in \S\ref{section:B}.


\subsection{Spectral sequences arising from the canonical filtration of a complex}\label{section:spectral}

In this subsection, we review some basic facts on spectral sequences.
See \cite[IV.2. Excercise 2]{Gelfand:mha} and the proof of \cite[III Proposition 4.4.6]{Verdier:cdca}
for details.

 Let $\mathcal{A}$ be an abelian category with enough injectives
and let $D(\mathcal{A})$ be
 the bounded derived category of it.
 For an object $\alpha \in D(\mathcal{A})$, we denote by
 $\cohom^i(\alpha) \in \mathcal{A}$ the $i$-th cohomology of the complex $\alpha$.
 For objects $\alpha, \beta \in D(\mathcal{A})$, there is a spectral sequence
 \begin{equation}\label{equation:spectral1}
 E_2^{p, q}=\bigoplus_i \Hom_{D(\mathcal{A})}^p(\cohom^i(\alpha), \cohom^{i+q}(\beta)) \Longrightarrow
 E^{p+q}=\Hom_{D(\mathcal A)}^{p+q}(\alpha, \beta).
 \end{equation}
 
 For a cohomological functor $F$ from $D(\mathcal{A})$ to
 an abelian category $\mathcal{B}$, we have another spectral sequence
 \begin{equation}\label{equation:spectral2}
 E_2^{p,q}= F^p(\cohom^q(\alpha)) \Longrightarrow E^{p+q}=F^{p+q}(\alpha).
 \end{equation}
We use \eqref{equation:spectral1} for a single spherical object $\alpha=\beta$
in the proof of Proposition \ref{proposition:step -1 of A_n},
and two spherical objects $\alpha=\Phi(\owe_{C_1})$ and $\beta=\Phi(\owe_{C_1}(-1))$ in the proof of Proposition
\ref{proposition:step -2 of A_n}.

In addition, we also use the description of the maps $d_2$ of the above spectral sequences.
 We denote by $\tau_{\le p}\alpha$ the following complex:
 $$
 (\tau_{\le p}\alpha)^n= \begin{cases}
                        \alpha^n & n < p \\
            \ker d^p & n = p \\
                     0 & n > p  
\end{cases}
 $$
We define
 $\tau_{> p}\alpha$($=\tau_{\ge p+1}\alpha$) so that it fits into a distinguished triangle
 $$
 \tau_{\le p}\alpha \to \alpha \to \tau_{>p}\alpha \to \tau_{\le p}\alpha[1]
 $$
 and we put
 $$
 \tau_{[p, q]}\alpha = \tau_{\ge p}\tau_{\le q}\alpha.
 $$
 Especially, we have an isomorphism $\tau_{[p, p]}\alpha \cong \cohom^p(\alpha)[-p]$ and
 a distinguished triangle
\begin{equation}\label{equation:e^2}
\cohom^{p-1}(\alpha)[-p+1] \to \tau_{[p-1, p]}\alpha \to \cohom^p(\alpha)[-p]
\to \cohom^{p-1}(\alpha)[-p+2].
\end{equation}
 The last morphism determines an element
 $$
 e^p(\alpha) \in \Hom_{D(\mathcal{A})}(\cohom^p(\alpha), \cohom^{p-1}(\alpha)[2])
 \cong \Ext^2_{\mathcal{A}}(\cohom^p(\alpha), \cohom^{p-1}(\alpha)).
 $$
 This class gives rise to the morphisms $d_2$ of the above spectral sequences: 
\begin{proposition}\label{differential map}
 The morphisms $d_2^{p, q} : E_2^{p, q} \to E_2^{p+2, q-1}$ in the spectral sequences
 in \eqref{equation:spectral1} and \eqref{equation:spectral2} are determined as follows.
 \begin{itemize}
 \item[\eqref{equation:spectral1}:]
 For $\oplus_i f_i\in \bigoplus_i \Hom_{D(\mc{A})}^p(\cohom^i(\alpha), \cohom^{i+q}(\beta))$,
 $$
 d_2^{p,q}(\oplus_i f_i)= \bigoplus_i ((-1)^{p+q}f_{i-1}\circ e^i(\alpha)-e^{i+q}(\beta)\circ f_i).
 $$
 \item[\eqref{equation:spectral2}:]
 $d_2^{p, q}$ is the morphism
 $F(e^q(\alpha)[p]): F(\cohom^q(\alpha)[p]) \to F(\cohom^{q-1}(\alpha)[p+2])$.
 \end{itemize}
 \end{proposition}


\subsection{Reconstruction of objects of the derived category of a smooth surface}\label{section:reconstruction}
 
 Let $X$ be a smooth surface.
We denote by $D(X)=D^b(\coh X)$ the bounded derived category of coherent sheaves on $X$.
The following proposition shows that an object $\alpha$ of $D(X)$ is determined by
its cohomology sheaves $\cohom^i(\alpha)$ and the classes $e^i(\alpha)$,
up to (non-canonical) isomorphisms.

\begin{proposition}\label{e^i.2}
Suppose we are given coherent sheaves $\mc{G}^i$ on $X$ and elements
$$
e^i\in \Ext^2_X(\mc{G}^i, \mc{G}^{i-1})
$$
for all $i\in\ZZ$
such that $\mc{G}^i$'s are zero except for finitely many $i$'s.  
Then there is an object $\alpha\in D(X)$ and isomorphisms $\mu_i:\mc{H}^i(\alpha)\cong\mc{G}^i$
such that $\mu_{i-1}[2]\circ e^i(\alpha)=e^i \circ \mu_i$. 
This $\alpha$ is uniquely determined up to isomorphisms.
\end{proposition}

\begin{proof}
Define $q_0=\max \bigl\{q \bigm| \mc{G}^q\ne 0 \bigr\}$ and $q_1=\min \bigl\{q \bigm| \mc{G}^q\ne 0 \bigr\}$.
We use induction on the non-negative integer $q_0-q_1$. When $q_0-q_1=0$,
we just define $\alpha$ to be $\mc{G}^{q_0}[-q_0]$.
Let us consider the case $q_0-q_1>0$. By the induction hypothesis, we can find $\beta\in D(X)$ and isomorphisms
$$
\nu_i:\mc{H}^i(\beta)\cong\begin{cases} \mc{G}^i & i\ne q_0\\
                              0 & i=q_0
\end{cases}
$$
such that $\nu_{i-1}[2]\circ e^i(\beta)=e^i \circ \nu_i$
if $i \ne q_0$.

Let us consider the spectral sequence \eqref{equation:spectral1}
$$
E_2^{p,q}=\Ext^p_X(\mathcal{G}^{q_0}, \mathcal{H}^{q_0+q}(\beta)) 
\Longrightarrow \Hom_{D(X)}^{p+q}(\mathcal{G}^{q_0}[-q_0], \beta). 
$$ 
Then, since $E_2^{p,q}=0$ for $q\ge 0$ or $p\not\in [0,2]$, we have an isomorphism
\[ 
\xymatrix{
f :\Ext^2_X(\mc{G}^{q_0},\mc{G}^{q_0-1})\ar[r]^(.6){(\nu_{q_0-1})_*}&E_2^{2,-1}\cong E^1.} 
\]
>From the morphism $-f (e^{q_0})[-1]$,
we obtain an object $\alpha\in D(X)$ and a distinguished triangle 
\[ 
\xymatrix{
\mc{G}^{q_0}[-q_0-1]\ar[rr]^(.6){-f (e^{q_0})[-1]}&& \beta\ar[r] &\alpha\ar[r] & \mc{G}^{q_0}[-q_0].
}
\]
We denote the last morphism by $\varphi$.
Then we have an isomorphism $\xi: \tau_{<q_0} \alpha \cong \beta$ and a morphism of distinguished triangles:
\[ 
\xymatrix{
\beta \ar[r]\ar[d] & \alpha \ar[r]^(.4){\varphi}\ar[d]& \mc{G}^{q_0}[-q_0] \ar[r]^(.6){f(e^{q_0})}\ar@{=}[d] &\beta[1] \ar[d]\\
\cohom^{q_0-1}(\beta)[-q_0+1] \ar[r]& \tau_{[q_0-1,q_0]}\alpha \ar[r] & \mc{G}^{q_0}[-q_0] \ar[r]&
\cohom^{q_0-1}(\beta)[-q_0+2]
}
\]
Here, the triangle in the second row is isomorphic to the one in \eqref{equation:e^2}.
Thus, putting  $\mu_i = \nu_i \circ \cohom^i(\xi)$ for $i \ne q_0$ and
$\mu_{q_0}=\cohom^{q_0}(\varphi)$,
we have $\mu_i: \cohom^i(\alpha) \cong \mc{G}^i$ and $\mu_{i-1}[2]\circ e^i(\alpha)=e^i \circ \mu_i$. 

For the uniqueness, let $\alpha$ and  $\beta$ be objects of $D(X)$
with isomorphisms $\xi_i: \cohom^i(\alpha) \cong \cohom^i(\beta)$
satisfying $\xi_{i-1}[2]\circ e^i(\alpha)=e^i(\beta) \circ \xi_i$.
Then $\oplus_i \xi_i$ lies in $E_2^{0,0}$ in the spectral sequence \eqref{equation:spectral1}
and the condition on $\xi_i$ implies that $d_2^{0,0}(\oplus_i \xi_i)=0$.
Since $X$ is non-singular of dimension $2$, $E_2^{p, q}$ vanishes unless $0 \le p \le 2$
and hence \eqref{equation:spectral1} is $E_3$-degenerate.
Therefore, $\oplus_i \xi_i$ survives at the infinity and there exists
$\xi \in \Hom_{D(X)}(\alpha, \beta)$ which induces $\xi_i$ on the cohomology sheaves.
Since each $\xi_i$ is an isomorphism, we see that $\xi$ is an isomorphism.
\end{proof}

In the light of Proposition \ref{e^i.2}, we obtain the following.


\begin{lemma}\label{decomposable}
Let $\alpha$ be an object of $D(X)$ which satisfies
$\cohom^i(\alpha)=\mc{G}_1^i\oplus \mc{G}_2^i$ 
for some coherent sheaves $\mc{G}_1^i, \mc{G}_2^i$. 
For the class $e^i(\alpha) \in \Ext^2_X(\cohom^i(\alpha),\cohom^{i-1}(\alpha))$,
we write
$$
e^i(\alpha)=\begin{pmatrix} a_i & b_i \\ c_i & d_i \end{pmatrix}
$$
so that
\begin{align*}
a_i &\in \Ext^2_X(\mc{G}_1^i,\mc{G}_1^{i-1}), \\
b_i &\in \Ext^2_X(\mc{G}_2^i,\mc{G}_1^{i-1}), \\
c_i &\in \Ext^2_X(\mc{G}_1^i,\mc{G}_2^{i-1}), \\
d_i &\in \Ext^2_X(\mc{G}_2^i,\mc{G}_2^{i-1})
\end{align*}
respectively.
If all $b_i$ and $c_i$ are zero, then we have objects $\alpha_1, \alpha_2\in D(X)$
such that $\alpha \cong \alpha_1 \oplus \alpha_2$,
$\cohom^i(\alpha_k) \cong \mc{G}^i_k$, $e^i(\alpha_1)=a_i$ and $e^i(\alpha_2)=d_i$.
\end{lemma}


\subsection{Spherical objects}\label{section:spherical}
The definition of a spherical object on an $n$-dimensional smooth quasi-projective variety $X$
is given by Seidel and Thomas:

\begin{definition}[\cite{Seidel:bga}]
We say that an object $\alpha \in D_c(X)$ is \emph{spherical} if 
we have $\alpha \otimes \omega_{X} \cong \alpha$ and
$$
\Hom^{k}_{D(X)}(\alpha,\alpha)\cong\begin{cases}  0 & k\ne 0,n\\
                                                \CC & k=0,n. 
\end{cases}
$$
\end{definition}
\noindent
Here suppose that $\dim X=2$ and take an object $\alpha$ of $D_c(X)$.
We shall give conditions for $\alpha$ to be spherical 
under the assumption $\alpha\otimes \omega_X\cong \alpha$.


\begin{proposition}\label{TFAE for spherical}
Assume that $\alpha\otimes \omega_X\cong \alpha$. The following are equivalent.
\begin{enumerate}
\item $\alpha$ is spherical.
\item In the spectral sequence \eqref{equation:spectral1} (for $\alpha=\beta$)
$$
E_2^{p,q}=\bigoplus _i\Hom^p_X(\mathcal{H}^i(\alpha), \mathcal{H}^{i+q}(\alpha)) 
\Longrightarrow \Hom_{D(X)}^{p+q}(\alpha, \alpha), 
$$
we have the following:
\begin{itemize}
\item $d_2^{0,q}$ is injective for all $q\ne 0$.
\item $\Ker d_2^{0,0}$ is a one-dimensional $\CC $-vector space generated by the element $\oplus_i id_i\in E_2^{0,0}$. 
\item $E_2^{1,q}=0$ for all $q$, i.e.,
$\Ext^1_X(\cohom^i(\alpha), \cohom^j(\alpha))=0$ for all $i, j$.
\end{itemize}
\end{enumerate}
\end{proposition}

\begin{proof}
Notice that the spectral sequence in (ii) degenerates at the $E_3$-level, since $X$ is two-dimensional.
We have 
\begin{equation}\label{duality}
E_2^{0,q}\cong (E_2^{2,-q})^\vee
\end{equation}
for all $q$ by the Grothendieck--Serre duality.

Let us first give the proof of the implication from (i) to (ii).
Notice that 
$$
\dim\Ker d_2^{0,0}=\dim E_3^{0,0}\le \dim E^0=1.
$$
Since $\oplus_i id_i\in \Ker d_2^{0,0}$, we obtain the second condition in (ii)
and $E_2^{1,-1}=E_3^{2,-2}=0$. Especially, we get $E_2^{1,1}=0$ by \eqref{duality}.
Since $\dim E_2^{1,q}\le \dim E^{1+q}=0$ for all $q\ne -1,1$, we have $E_2^{1,q}=0$ for all $q$, as desired.
Now let us show the first condition in (ii).
Obviously, the condition (i) implies that $d_2^{0,q}$ is injective for $q\ne 0,2$. 
On the other hand, we know that $d_2^{0,2}$ is surjective by $E^3=0$ and $d_2^{0,-1}$ is isomorphic by $E_3^{2,-2}=0$.
In particular, we see
$$
\dim\Ker d_2^{0,2}=\dim E_2^{2,1}-\dim E_2^{0,2}=\dim E_2^{0,-1}-\dim E_2^{2,-2}=0,
$$
which implies the conclusion.

Conversely, assume that (ii) holds. We have
$$
\dim E_2^{0,q}-\dim E_2^{2,q-1}\le\dim\Ker d_2^{0,q}=
\begin{cases}
 1 & q=0 \\
 0 & q\ne 0.
\end{cases}
$$
Combining this and \eqref{duality} together, we get 
$$
\dim E_2^{0,q}=\dim E_2^{2,q-1}
$$
for $q\ne 0,1$.
Since $d_2^{0,q}$ is injective for $q\ne 0$,
we know that $d_2^{0,q}$ is isomorphic for $q\ne 0,1$, in particular, 
$E_3^{2,-2}=\Coker d_2^{0,-1}=0$. This equality and (ii) imply
$$
\Coker d_2^{0,1}\cong E_3^{2,0}\cong\Hom^2_{D(X)}(\alpha,\alpha )
$$
and
$$
\Hom^0_{D(X)}(\alpha,\alpha )\cong E_3^{0,0}\cong\Ker d_2^{0,0}\cong\CC .
$$
Hence it follows from the duality that 
$$
\Coker d_2^{0,1}\cong \Hom^0_{D(X)}(\alpha,\alpha )^\vee\cong \CC .
$$
Therefore we have 
$$
\dim E_2^{0,0}-\dim E_2^{2,-1}=\dim E_2^{2,0}-\dim E_2^{0,1}=\dim\Coker d_2^{0,1}=1.
$$
Especially, we get the surjectivity of $d_2^{0,0}$ and
$$
\dim \Hom^1_{D(X)}(\alpha,\alpha )=\dim E_3^{2,-1}=\dim\Coker d_2^{0,0}=0.
$$
This completes the proof.
\end{proof}


\begin{remark}
Via Proposition \ref{differential map}, Proposition \ref{TFAE for spherical} (ii)
is regarded as a condition on $\cohom^i(\alpha)$ and $e^i(\alpha)$.
Consequently the condition for $\alpha\in D(X)$ to be spherical is entirely expressed
in terms of $\cohom^i(\alpha)$ and $e^i(\alpha)$.
\end{remark}


\begin{example}\label{A_5}
Let $X$ be a smooth surface.
\begin{enumerate}
\item
Let $Z$ be a chain of $-2$-curves on $X$ and  $\mc{L}$ a line bundle on $Z$. Then $\mc{L}$ is 
a spherical object of $D(X)$.
\item
We give a rather non-trivial example of a spherical object $\alpha\in D(X)$, supported on $C_1\cup\cdots\cup C_5$,
a union of $-2$-curves in $A_5$-configuration on $X$.
First we define the cohomology sheaves of $\alpha$ as follows:
\[ 
\xymatrix@R=1ex@M=0ex{     &   & C_1 & C_2 & C_3 & C_4 & C_5  \\
          \mc{H}^2(\alpha):&   &{\zero}\ar@{-}[r] &{\mone}\ar@{-}[r] &{\zero} & & &   \\
          \mc{R}_1        :&   &{\mone}\ar@{-}[r]&{\zero}\ar@{-}[r]&{\zero}\ar@{-}[r]&{\zero }&          &   \\    
          \mc{R}_2        :&   &{\zero}\ar@{-}[r]&{\zero}\ar@{-}[r]&{\mone}          &        &          &   \\   
          \mc{H}^0(\alpha):&   &{\mone}\ar@{-}[r]&{\zero}\ar@{-}[r]&{\zero}\ar@{-}[r]&{\zero}\ar@{-}[r]&{\zero}& }      
\]
with $\mc{H}^1(\alpha)=\mc{R}_1\oplus\mc{R}_2$. Notice that
$$
\Ext_X^2(\mc{H}^2(\alpha),\mc{H}^1(\alpha))\cong
\Ext_X^2(\mc{H}^2(\alpha),\mc{R}_1)\oplus\Ext_X^2(\mc{H}^2(\alpha),\mc{R}_2)\cong\CC\oplus\CC  
$$
and
$$
\Ext_X^2(\mc{H}^1(\alpha),\mc{H}^0(\alpha))\cong
\Ext_X^2(\mc{R}_1,\mc{H}^0(\alpha))\oplus\Ext_X^2(\mc{R}_2,\mc{H}^0(\alpha))\cong\CC\oplus\CC .  
$$
Keep these isomorphisms in mind, and take 
$$
e^2(\alpha)=(e^2_1, e^2_2)\in\Ext_X^2(\mc{H}^2(\alpha),\mc{H}^1(\alpha))
$$
and 
$$
e^1(\alpha)=(0, e^1_2)\in\Ext^2_X(\mc{H}^1(\alpha),\mc{H}^0(\alpha))
$$
with $e^2_1, e^2_2, e^1_2\in\CC ^*$. 
The data $\mc{H}^i(\alpha)$ and $e^i(\alpha)\in \Ext_X^i(\mc{H}^i(\alpha),\mc{H}^{i-1}(\alpha))$
determine an object $\alpha\in D(X)$ by Proposition \ref{e^i.2}.
We can see that $\alpha$ is spherical by checking the conditions in Proposition \ref{TFAE for spherical} (ii).
\end{enumerate}
\end{example}

Proposition \ref{TFAE for spherical} holds for any compactly supported object on a smooth surface $X$.
In the situation of our problem, we can say more about
the cohomology sheaves of a spherical object.


\begin{lemma}\label{cohomology sheaf}
Let $f:X\to Y$ be a surjective morphism from a smooth variety $X$ to a variety $Y$,
and let $Z=f^{-1}(y)$ be the scheme-theoretic fiber of a closed point $y \in Y$.
If $\alpha \in D_Z(X)$ satisfies $\Hom_{D(X)}(\alpha, \alpha) \cong \CC$,
then every cohomology sheaf $\mathcal{H}^i(\alpha)$ is an $\mathcal{O}_Z$-module.
\end{lemma} 

\begin{proof}
Take an affine open neighborhood $U:=\Spec R$ of $y$
and denote by $m_y \subset R$ the maximal ideal of $y$ in $U$.
Then the spectral sequence \eqref{equation:spectral1} is a spectral sequence of $R$-modules
and we have $E^0=\Hom_{D(X)}(\alpha, \alpha)\cong R/m_y$.
On the other hand, this spectral sequence satisfies
$$
E_2^{0,0} \supset E_3^{0,0} \supset\dots \supset E_\infty^{0,0}
$$
and the image of $E^0$ in $E_\infty^{0,0} \hookrightarrow E_2^{0,0}=\bigoplus _i\Hom_X(\mathcal{H}^i(\alpha), \mathcal{H}^{i}(\alpha))$
contains $\oplus_i id_i$.
Thus, for each identity map $id_i$ on $\mathcal{H}^i(\alpha)$, we have
$m_y\cdot id_i=0$ and $\mathcal{I}_Z\cdot\mathcal{H}^i(\alpha)=0$.  
This completes the proof. 
\end{proof}


Recall that a coherent sheaf $\mc{F}$ on a variety $X$ is {\it rigid}
if $\Ext^1_X(\mc{F}, \mc{F})=0$.

\begin{lemma}\label{pure}
Let $\mathcal{F}$ be a one-dimensional rigid coherent sheaf on a smooth surface $X$.
Then $\mathcal{F}$ is purely one-dimensional, that is, every non-zero subsheaf of $\mc{F}$
is one-dimensional. 
\end{lemma} 

\begin{proof}
Let $\mathcal{F}_{\tor}$ be the `torsion' part of $\mathcal{F}$, namely 
the maximal zero-dimensional subsheaf of $\mathcal{F}$. 
Our aim is to show $\mathcal{F}_{\tor}=0$.
Take a surjection $\mathcal{E}\to\mathcal{F}$ from a locally free sheaf $\mathcal{E}$ and 
denote the kernel of it by $\mathcal{G}$.
We consider the following commutative diagram with exact rows.

\noindent
\[ \xymatrix{
   0 \ar[r] & \mathcal{G} \ar[d]^g \ar[r] & \mathcal{E} \ar[r] \ar@{=}[d] & \mathcal{F} \ar[r] \ar[d]^f & 0  \\
   0 \ar[r] & \mathcal{G}^{\vee\vee} \ar[r] &  \mathcal{E} \ar[r]           &  \mathcal{F}' \ar[r]      & 0
} \]

\noindent
Here, $\mathcal{G}^{\vee\vee}$ is the double dual
of $\mathcal{G}$ and $\mathcal{F}' = \mathcal{F}/\mc{F}_{\tor}$.
Note that $\mathcal{F}_{\tor}\cong \mathcal{G}^{\vee\vee}/\mathcal{G}$
by the snake lemma.
Let us consider the composite of the natural maps
$$
\hom_X(\mc{G}, \mc{F}_{\tor}) \hookrightarrow \hom_X(\mc{G}, \mc{F}) \rightarrow \ext_X^1(\mc{F}, \mc{F})
$$
and denote it by $\varphi$ .
Since
$\hom_X(\mathcal{G},\mathcal{F}_{\tor})$ is a zero-dimensional sheaf,
the vanishing of $H^0(\ext^1_X(\mathcal{F},\mathcal{F}))$ implies that
$\varphi$ is the zero map.
This means that in the exact sequence
$$
\hom_X(\mathcal{E},\mathcal{F})\to
  \hom_X(\mathcal{G},\mathcal{F})\to
  \ext^1_X(\mathcal{F},\mathcal{F})\to 0,
$$
we can extend a (local) map
$\psi\in \hom_X(\mathcal{G},\mathcal{F}_{\tor})$ to a (local) map $\bar{\psi}\in\hom_X(\mathcal{E},\mathcal{F})$.
$\bar{\psi}$ sends $(\mathcal{G}^{\vee\vee})$ into $\mathcal{F}_{\tor}$, 
since $\mathcal{G}^{\vee\vee}/\mathcal{G}(\cong\mathcal{F}_{\tor})$ is zero-dimensional. 
Therefore, in the exact sequence
$$
\hom_X(\mathcal{G}^{\vee\vee},\mathcal{F}_{\tor})\to
  \hom_X(\mathcal{G},\mathcal{F}_{\tor})\to
  \ext^1_X(\mathcal{F}_{\tor},\mathcal{F}_{\tor})\to 0,
$$
the first map is surjective. It follows that 
$\ext^1_X(\mathcal{F}_{\tor},\mathcal{F}_{\tor})=0$.
Since $\mc{F}_{\tor}$ is zero-dimensional and rigid, 
we obtain $\mc{F}_{\tor}=0$.
\end{proof}

We summarize Proposition \ref{TFAE for spherical}, Lemmas \ref{cohomology sheaf} and
\ref{pure} in our situation as follows.


\begin{corollary}\label{rigid & pure}
Let $\{C_i\}$ be a collection of $-2$-curves in an $ADE$-configuration on
a smooth surface $X$ and let $Z$ be the fundamental cycle of $\bigcup _iC_i$.
If $\alpha\in D_Z(X)$ is a spherical object,
then the sheaf $\bigoplus _p\mathcal{H}^p(\alpha)$ is a rigid 
$\mc{O}_Z$-module, pure of dimension $1$.
\end{corollary}

Recall we defined $l(\alpha)$ for an object $\alpha \in D_Z(X)$ in Introduction.
The following is a basic tool in the proofs of Propositions \ref{proposition:step -1 of A_n}
and \ref{proposition:step -2 of A_n}.


\begin{lemma}\label{strategy} 
Under the notation in Corollary \ref{rigid & pure},
we have 
\begin{equation}\label{equation:le}
l(\Phi(\alpha))\le\sum_q l(\Phi (\cohom ^q(\alpha)))
\end{equation}
for any $\Phi \in \Auteq D_Z(X)$.
The equality in \eqref{equation:le} implies
the vanishing $d_2^{p,q}=0$ for all $p,q$
in the spectral sequence
\begin{equation}\label{equation:spe}
E_2^{p,q}= \mc{H}^p(\Phi (\cohom^q(\alpha))) \Longrightarrow E^{p+q}=\mc{H}^{p+q}(\Phi(\alpha)),
\end{equation}
if every $E_2^{p,q}$ is purely one-dimensional. 
\end{lemma}

\begin{proof}
In \eqref{equation:spe}, we see that
\begin{align*}
l(\Phi(\alpha))=&\sum_n l(E^n)=\sum_{p,q}l(E_{\infty}^{p,q})\le\cdots\\
\le&\sum_{p,q}l(E_3^{p,q})\le\sum_{p,q}l(E_2^{p,q})=\sum_q l(\Phi(\mc{H}^q(\alpha))),
\end{align*}
which implies \eqref{equation:le}.
If the equality holds in \eqref{equation:le}, then 
$\sum_{p,q}l(E_3^{p,q})=\sum_{p,q}l(E_2^{p,q})$.
This ensures $l(\Image(d_2^{p,q}))=0$,
and consequently $\dim \Image(d_2^{p,q}) \le 0$.
Since $\Image(d_2^{p,q})$ is a subsheaf of $E_2^{p+2,q-1}$
which is pure of dimension $1$, it must be zero.
\end{proof}

\begin{remark}\label{remark:always pure}
If $Z$ forms an $A_n$-configuration in Lemma \ref{strategy}, 
we can actually show that every $E_2^{p,q}$ is always purely one-dimensional 
by Corollary \ref{rigid & pure} and Lemma \ref{sheaf on A_n}. 
\end{remark}


\subsection{Twist functors}\label{section: Twist functors}
Let $X$ be an $n$-dimensional smooth quasi-projective variety.
The following definition is due to Seidel and Thomas.
\begin{definition}[\cite{Seidel:bga}]
Let $\alpha\in D_{c}(X)$ be a spherical object and
consider the mapping cone 
$$
\mc{C}=Cone(\pi_1^*\alpha^\vee\ltensor \pi_2^*\alpha \to \mc{O}_{\Delta})
$$
of the natural evaluation $\pi_1^*\alpha^\vee\ltensor \pi_2^*\alpha \to \mc{O}_{\Delta}$,
where $\Delta\subset X\times X$ is the diagonal and $\pi_i$ is the $i$-th projection
$\pi_{i}: X\times X\to X$. 
Then $T_{\alpha}:=\Phi_{X\to X}^{\mc{C}}$ defines an autoequivalence, called
the \emph{twist functor} along a spherical object $\alpha$. 
The object $T_{\alpha}(\beta)$ fits into a distinguished triangle
\[
\mathbf{R}\Hom_{\owe_{X}}(\alpha, \beta) \ltensor_{\CC} \alpha \overset{\ev}{\longrightarrow} \beta \longrightarrow T_{\alpha}(\beta)
\]
for any $\beta \in D(X)$, where ev is the evaluation morphism. 
For the inverse \(T'_{\alpha}\) of $T_{\alpha}$, we have a distinguished triangle
\[
T'_{\alpha}(\beta) \longrightarrow \beta  \overset{\text{ev}}{\longrightarrow} \mathbf{R}\Hom_{\owe_{X}}(\beta,\alpha)^{\vee} 
\ltensor_{\CC} \alpha
\] 
for any $\beta \in D(X)$.
\end{definition}

We list several lemmas on twist functors that will be used later. 


\begin{lemma}\label{easy}
\begin{enumerate}
\item
Let $\alpha \in D(X)$ be a spherical object.
For an FM transform $\Phi: D(X)\to D(X)$ with quasi-inverse $\Phi ^{-1}$, we have 
$$
\Phi\circ T_\alpha \circ\Phi ^{-1}\cong T_{\Phi(\alpha)}.
$$ 
For an integer $i$, we also have
$$
T_\alpha \cong T_{\alpha[i]}.
$$
\item 
Let $Z\subsetneq X$ be a closed subscheme of $X$ which is proper over $\CC$.
Then we have 
$$
\Span{T_\alpha \bigm| \alpha \in D_Z(X), \textit{ spherical } }\cap \Aut X=\{id\}.
$$
\end{enumerate}
\end{lemma}

\begin{proof}
(i) is readily verified by definition. 
The kernel $\mc{P}$ of an integral functor $\Phi^{\mc{P}}$ in the left hand side of (ii) satisfies that 
$\mc{P}|_{(X\backslash Z)\times (X\backslash Z)}\cong \mc{O}_{\Delta}|_{(X\backslash Z)\times (X\backslash Z)}$,
where 
$\Delta \subset X\times X$ is the diagonal. This leads us to the equality in (ii).
\end{proof}


\begin{lemma}\label{calculation}
Let $X$ be a smooth surface.
\begin{enumerate}
\item
For a $-2$-curve $C$ on $X$ 
and an integer $a$, we have the following:
\begin{enumerate}
\item
$$
T_{\mc{O}_C(a)}(\mc{O}_C(a))=\mc{O}_C(a)[-1]
$$
and
$$
T_{\mc{O}_{C}(a-1)}(\mc{O}_C(a))=\mc{O}_C(a-2)[1].
$$
\item
$$
T_{\mc{O}_C(a-1)}\circ T_{\mc{O}_C(a)}\cong \otimes\mc{O}_X(C).
$$
\end{enumerate}

\item 
Let $Z=\sum _{l=1}^nC_l$ be a chain of $-2$-curves $C_l$ on $X$ with $n>1$ 
and put $\alpha=\mc{O}_Z(a_1,a_2,\dots,a_n)$ for some $a_l\in\ZZ$.
Then we have the following:
\begin{enumerate}
\item
$$
\mc{H}^p(T_{\mc{O}_{C_1}(a_1)}(\alpha))=
\begin{cases}
\alpha & p=0 \\
\mc{O}_{C_1}(a_1) & p=1 \\
0 & p\ne 0,1.
\end{cases}
$$ 

\item
$$
T_{\mc{O}_{C_1}(a_1-1)}(\alpha)=\mc{O}_{C_2\cup \cdots\cup C_n}(a_2,\dots,a_n).
$$

\item
$$
\mc{H}^p(T_{\mc{O}_{C_1}(a_1-2)}(\alpha))=
\begin{cases}
\mc{O}_{C_1}(a_1-3) & p=-1 \\
\mc{O}_Z(b_1,\dots,b_n) & p=0\\
0 & p\ne -1,0.
\end{cases}
$$
Here 
$$
b_l=
\begin{cases}
a_1-2 & l=1 \\
a_2+1 & l=2 \\
a_l & l\ne 1,2.
\end{cases}
$$

\item
$$
T_{\mc{O}_{C_k}(a_k-1)}(\alpha)=\alpha
$$
for all $k$ $(1<k<n)$.

\item
$$
T_{\mc{O}_{C_k}(a_k-2)}(\alpha)
=\mc{O}_Z(b_1,\dots,b_n)
$$
for all $k$ $(1<k<n)$.

Here
$$
b_l=
\begin{cases}
a_l & i\ne k-1,k,k+1\\
a_l+1 & l= k-1,k+1\\
a_k-2 & l= k.
\end{cases}
$$
\end{enumerate}
\end{enumerate}
\end{lemma}

\begin{proof}
(i.1) and (ii) are easy calculations.
It follows from (i.1) that
$T_{\mc{O}_C(a-1)}\circ T_{\mc{O}_C(a)}$
sends $\owe_C(a)$ to $\owe_C(a-2)$ and $\owe_C(a+1)$ to $\owe_C(a-1)$.
Hence, for any point $x\in X$,
$$
T_{\mc{O}_C(a-1)}\circ T_{\mc{O}_C(a)}(\mc{O}_x)\cong \mc{O}_y
$$
for some $y \in X$.
Thus Lemma \ref{A(X)} implies that 
$T_{\mc{O}_C(a-1)}\circ T_{\mc{O}_C(a)}$ is an element of $A(X)$.
Lemma \ref{easy} (ii) then yields it must be $\otimes L$ for some line bundle $L$.
Since $\owe_C(a) \otimes L \cong \owe_C(a-2)$, we see $L\cong \owe_X(C)$.
\end{proof}


\subsection{On the group $B$}\label{section:B}
Let $Z=C_1\cup\cdots\cup C_n \subset X$ be as in Introduction.
Recall we defined
$$
B=\Span{T_{\owe_{C_l}(-1)}, T_{\omega_Z}\bigm| 1\le l\le n } \subset \Auteq D_Z(X),
$$
where $\omega_Z$ denotes the dualizing sheaf on $Z$.
Put
$$
B'=\Span{T_{\owe_{C_l}(a)}\bigm|  a\in \ZZ, 1\le l\le n }.
$$
Then we have


\begin{lemma}\label{B=B'}
$B=B'$.
\end{lemma}
\begin{proof}
The proof is by induction on $n$.
When $n=1$, we write $C=C_1$.
In this case, $B=\Span{T_{\owe_C(-2)}, T_{\owe_C(-1)}} \subset B'$ by definition.
Then, Lemma \ref{calculation} (i.2) shows $\otimes \owe_X(C) \in B$.
Thus, we obtain from Lemma \ref{easy} (i)
$$
T_{\mc{O}_C(2a-2)}\cong\otimes\mc{O}_X(-aC)\circ T_{\mc{O}_C(-2)}\circ \otimes\mc{O}_X(aC)
\in B
$$
and
$$
T_{\mc{O}_C(2a-1)}\cong\otimes\mc{O}_X(-aC)\circ T_{\mc{O}_C(-1)}\circ \otimes\mc{O}_X(aC)
\in B.
$$ 

Let us consider the case $n>1$.
By the induction hypothesis, we have
\begin{equation}\label{equation:Z_n}
\Span{T_{\owe_{C_l}(a)}\bigm|  a\in \ZZ, 2 \le l \le n}
=\Span{T_{\owe_{C_l}(-1)}, T_{\omega_{Z_1}}\bigm|  2 \le l\le n}
\end{equation}
and
\begin{equation}\label{equation:Z_1}
\Span{T_{\owe_{C_l}(a)}\bigm|  a\in \ZZ, 1\le l\le n-1 }
=\Span{T_{\owe_{C_l}(-1)}, T_{\omega_{Z_n}}\bigm| 1\le l\le n-1 },
\end{equation}
where $Z_1=\sum ^n_{l=2}C_l$ and $Z_n=\sum ^{n-1}_{l=1}C_l$.
Since we have
$$
T_{\omega_{Z_1}}\cong T_{\mc{O}_{C_1}(-1)}'\circ T_{\omega_Z}\circ T_{\mc{O}_{C_1}(-1)}\in B
$$ 
and 
$$
T_{\omega_{Z_n}}\cong T_{\mc{O}_{C_n}(-1)}'\circ T_{\omega_Z}\circ T_{\mc{O}_{C_n}(-1)}\in B
$$ 
by Lemmas \ref{easy} (i) and \ref{calculation} (ii), 
\eqref{equation:Z_n} and \eqref{equation:Z_1} show that $T_{\owe_{C_l}(a)}\in B$
for all $l$ $(1\le l\le n)$, that is, $B' \subset B$.
Conversely, we see from Lemmas \ref{easy} (i) and \ref{calculation} (ii) that
$T_{\omega_Z} \in B'$. Thus we obtain $B = B'$.
\end{proof}
\noindent
We further see in Corollary \ref{corollary:allspherical} that
$T_{\alpha} \in B$ for every spherical object $\alpha \in D_Z(X)$.

\begin{remark}\label{remark:normal sub'}
We see from Lemma \ref{easy} (i) and Lemma \ref{B=B'} that $B$ is a normal subgroup of $\Span{A(X),B}$.
It also follows from Lemma \ref{easy} (ii) that $B \cap \Aut X = \{\text{id}\}$.
\end{remark}

Next we consider the relation between $B$ and $\Pic X$ in $\Auteq D_Z(X)$.


\begin{proposition}\label{subgroups}
We have the following.
\begin{enumerate}
\item
$B\cap\Pic X=\Span{\otimes\mc{O}_X(C_1),\ldots,\otimes\mc{O}_X(C_n)}$.
\item
$\Span{B, \Pic X} \cong B \rtimes \ZZ/(n+1)\ZZ$.
\end{enumerate}
\end{proposition}
\begin{proof}
(i) Lemma \ref{calculation} (i.2) implies
that the right hand side is contained in the left hand side.
Let $i:X\backslash Z \to X$ be the open immersion. For a spherical object $\alpha\in D_Z(X)$,
we have $(i^*\circ T_{\alpha})(\mc{O}_X)\cong \mc{O}_{X\backslash Z}$.
Hence for an autoequivalence  $\otimes \mc{L}\in B\cap\Pic X$,
we have $i^*\mc{L}\cong \mc{O}_{X\backslash Z}$. Thus $\mc{L}$ belongs to the right hand side.

(ii) Note that the natural map
$$
\deg:\Pic X\to\ZZ ^{\oplus n}\quad \mc{L}\longmapsto (\deg \mc{L}|_{C_l})_l
$$
is isomorphic \cite{Artin:ira}.
We denote by $\mc{O}_X(a_1,\ldots,a_n)$ the element of $\Pic X$ which goes to $(a_1,\ldots,a_n)\in \ZZ ^{\oplus n}$.

By (i), $B \cap \Pic X$ can be regarded as the root lattice;
then $\Pic X$ is the weight lattice of it.
As is well-known (see \cite[\S13, Exercise 4]{Humphreys:LART}),
the weight lattice modulo the root lattice of type $A_n$ is
isomorphic to $\ZZ/(n+1)\ZZ$.
Thus, we have
$$
\Span{B,\Pic X}/B\cong \Pic X/(B\cap\Pic X)\cong \ZZ/(n+1)\ZZ.
$$
Put
$$
\Phi_0:=T_{\mc{O}_{C_1}(-1)}\circ\cdots \circ T_{\mc{O}_{C_n}(-1)}\circ \otimes \mc{O}_X(0,\ldots,0,1),
$$
and $\alpha_l:=\mc{O}_{C_l}(-1)$ for $l=1,\ldots,n$, and $\alpha_0:=\alpha_{n+1}:=\omega_Z[1]$.
Then we can show by direct computation that $\Phi_0(\alpha_l)\cong \alpha _{l+1}$ for $l=0,\ldots,n$.
Thus we have ${\Phi_0}^{n+1}(\alpha_l)\cong \alpha_l$ for all $l$ $(0\le l\le n)$,
which implies that for any point $x\in C_l$,
we obtain ${\Phi_0}^{n+1}(\mc{O}_{x})\cong \mc{O}_y$ for some $y\in C_l$. Then we get
${\Phi_0}^{n+1}\in \Aut X\cap B$, and
therefore ${\Phi_0}^{n+1}\cong id$ by Lemma \ref{easy} (ii)
and $\Span{B, \Pic X} \cong B \rtimes \Span{\Phi _0}$.
\end{proof}

\begin{remark}
Consider the McKay correspondence $D_Z(X) \cong D_{\{0\}}(\CC^2)$.
Then it is easy to find an autoequivalence
of $D_{\{0\}}(\CC^2)$ of order $n+1$.
In fact, tensoring by a one-dimensional representation of $G$ is such an equivalence
and this lies in our subgroup.
\end{remark}

Finally, we state a fact which we frequently use in the proofs of Propositions 1.6 and 1.7.

\begin{lemma}\label{lemma:normal sub}
Let $\alpha$ be an object of $D_Z(X)$.
If there is $\Psi_0\in\Span{B, A(X)}$ 
such that $l(\Psi_0(\alpha))<l(\alpha)$, then there is
$\Psi \in B$ with the same property.
\end{lemma}
\begin{proof}
We know by Remark \ref{remark:normal sub'} that $B$ is a normal subgroup of $\Span{B, A(X)}$,
and by definition that $l(\alpha)=l(\Psi(\alpha))$ holds
for $\Psi\in A(X)$ and for a spherical object $\alpha\in D_Z(X)$. 
The assertion follows from this.
\end{proof}

 \section{The $A_1$ cases of Propositions \ref{proposition:step -1 of A_n} 
and \ref{proposition:step -2 of A_n}}\label{section:a_1}

In this section, we consider the $A_1$ cases of Propositions \ref{proposition:step -1 of A_n} 
and \ref{proposition:step -2 of A_n};
thus we are given a single $-2$-curve $C=Z$.
Let $\alpha \in D_Z(X)$ be a spherical object.
By Corollary \ref{rigid & pure},
we may assume that there is an integer $a$
such that
$$
\cohom^p(\alpha) \cong \owe_{C}(a-1)^{\oplus r_p} \oplus \owe_{C}(a)^{\oplus s_p}
$$
for all $p$, where $r_p$ and $s_p$ are non-negative integers.
In this case, $l(\alpha)$ is written as
$$
l(\alpha)=\sum_p (r_p + s_p).
$$


\begin{proposition}[The $A_1$ case of Proposition \ref{proposition:step -1 of A_n}]\label{proposition:spherical}
Let $\alpha\in D_Z(X)$ be a spherical object.
Then, there are integers $a, i$ and a functor $\Psi \in B$ such that
$$
\Psi(\alpha) \cong\owe_C(a)[i]
$$
\end{proposition}

\begin{proof}
Since we have $B=B'$ by Lemma \ref{B=B'},
it suffices to show the following:
\begin{claim}\label{claim:a1}
If $l(\alpha)>1$, then $l(T_{\owe_C(a-1)}(\alpha)) < l(\alpha)$.
\end{claim}
\noindent
The class $e^q(\alpha) \in \Ext^2_X(\cohom^q(\alpha), \cohom^{q-1}(\alpha))$
is of the form
$$
e^q(\alpha)=\begin{pmatrix} a_q & b_q \\ c_q & d_q \end{pmatrix},
$$
where
\begin{align*}
a_q &\in \Ext^2_X(\owe_C(a-1)^{\oplus r_q}, \owe_C(a-1)^{\oplus r_{q-1}}), \\
b_q &\in \Ext^2_X(\owe_C(a)^{\oplus s_q}, \owe_C(a-1)^{\oplus r_{q-1}}), \\
c_q &\in \Ext^2_X(\owe_C(a-1)^{\oplus r_q}, \owe_C(a)^{\oplus s_{q-1}})=0, \\
d_q &\in \Ext^2_X(\owe_C(a)^{\oplus s_q}, \owe_C(a)^{\oplus s_{q-1}})
\end{align*}
respectively.

Consider the spectral sequence \eqref{equation:spectral2}:
$$
E_2^{p, q}=\cohom ^p(T_{\owe_C(a-1)}(\cohom^q(\alpha))) \Longrightarrow
\cohom ^{p+q}(T_{\owe_C(a-1)}(\alpha)).
$$
In this spectral sequence, we have
\begin{align*}
E_2^{-1, q}&=\cohom ^{-1}( T_{\owe_C(a-1)}(\cohom^q(\alpha)))\cong \owe_C(a-2)^{\oplus s_q},
\\
E_2^{1, q}&=\cohom ^{1}(T_{\owe_C(a-1)}(\cohom^q(\alpha)))\cong \owe_C(a-1)^{\oplus r_q}
\end{align*}
and $E_2^{p, q}=0$ for $p \ne \pm 1$ by Lemma \ref{calculation}.
Especially, Lemma \ref{strategy} implies $l(T_{\owe_C(a-1)}(\alpha)) \le l(\alpha)$; if the equality holds,
then $d_2^{-1, q}=0$ for all $q$.
Assume, by contradiction, that $d_2^{-1, q}=0$ for all $q$.
Then we see by Proposition \ref{differential map}
that $b_q=0$ for all $q$.
Therefore we have
$$
e^q(\alpha)=\begin{pmatrix} a_q & 0 \\ 0 & d_q \end{pmatrix}.
$$
Lemma \ref{decomposable} implies that there are objects $\alpha_1,\alpha_2\in D_Z(X)$ 
such that $\alpha \cong \alpha_1 \oplus \alpha_2$ with $\cohom^q(\alpha_1) \cong \owe_C(a-1)^{\oplus r_q}$
and $\cohom^q(\alpha_2) \cong \owe_C(a)^{\oplus s_q}$.
Since $\alpha$ is spherical, either $\alpha_1$ or $\alpha_2$ must be zero.
Let $q_0, q_1$ be the maximum and the minimum of the integers $q$
with $\cohom^q(\alpha) \ne 0$.
Since $\alpha\cong \alpha_1$ or $\alpha_2$,
we have $\Hom_X(\cohom^{q_0}(\alpha), \cohom^{q_1}(\alpha)) \ne 0$.
If $q_0 > q_1$, then the spectral sequence \eqref{equation:spectral1} for $\alpha=\beta$
implies that $\Hom_{D(X)}^{q_1-q_0}(\alpha, \alpha) \ne 0$
contradicting the assumption that $\alpha$ is spherical.
Thus we have $q_0=q_1$.
Then, since $\dim \Hom_{D(X)}(\alpha, \alpha)=1$, $l(\alpha)$ must be $1$.
\end{proof}


\begin{proposition}[The $A_1$ case of Proposition \ref{proposition:step -2 of A_n}]\label{proposition:A_1}
Let $\Phi$ be an autoequivalence of $D_Z(X)$.
Then, there are integers $a$ and $i$, and 
there is an autoequivalence $\Psi\in B$ such that
$$
\Psi\circ \Phi(\owe_C) \cong \owe_C(a)[i]
$$
and
$$
\Psi\circ \Phi(\owe_C(-1)) \cong \owe_C(a-1)[i].
$$
In particular, 
for any point $x\in C$, we can find a point $y\in C$ 
with 
$$
\Psi\circ\Phi (\mc{O}_x)\cong \mc{O}_y[i].
$$  
\end{proposition}

\begin{proof}
Put $\alpha=\Phi(\owe_C)$ and $\beta=\Phi(\owe_C(-1))$.
By Proposition \ref{proposition:spherical}, we may assume $l(\alpha)=1$.
We can further assume $\alpha \cong \owe_C$ by Lemma \ref{lemma:normal sub}.
Note that we have
\begin{equation}\label{eq:abc}
\Hom_{D(X)}^q(\beta, \owe_C)\cong\begin{cases} \CC ^2& q=0\\
                                                 0& q\ne 0.
\end{cases}
\end{equation}
We prove the first statement in the proposition by induction on $l(\beta)$;
the second follows from the first.
We first consider the case $l(\beta)=1$.
Then \eqref{eq:abc} implies that
$\beta$ is isomorphic to either $\owe_C(-1)$ or $\owe_C(1)[-2]$.
In the latter case, $T_{\owe_C}(\alpha) \cong \owe_C[-1]$
and $T_{\owe_C}(\beta) \cong \owe_C(-1)[-1]$ as desired.

Next assume $l(\beta)>1$.
As before, there is an integer $a$ such that
$$
\cohom^q(\beta) \cong \owe_{C}(a-1)^{\oplus r_q} \oplus \owe_{C}(a)^{\oplus s_q}.
$$
Let $q_0, q_1$ be the maximum and the minimum of the integers $q$
with $\cohom^q(\beta) \ne 0$.
If $q_0=q_1$, then $l(\beta)$ must be $1$ since $\beta$ is spherical.
If  $q_0>q_1$, then we have $\Hom_X^0( \cohom^{q_0}(\beta), \cohom^{q_1}(\beta))= 0$
and hence $r_{q_0}=s_{q_1}=0$.
Then we can see that either $\Hom_X^2( \cohom^{q_1}(\beta),\owe_C)$ or
$\Hom_X^0(\cohom^{q_0}(\beta),\owe_C)$ is non-zero.  
It follows from \eqref{eq:abc} and the spectral sequence
$$
E_2^{p,q}=\Hom_X^p(\cohom^{-q}(\beta),\owe_{C}) \Longrightarrow \Hom_{D(X)}^{p+q}(\beta,\owe_C)
$$
that $q_0=0$ or $q_1=2$, and in particular that $\cohom^1(\beta)=0$. 
Consequently we have 
$\Hom_X^1(\cohom^{-q}(\beta),\owe_{C})=0$ for all $q$
and hence that $a=0$ or $1$.
Therefore, we have $l(T_{\owe_C(a-1)}(\alpha))=1$.
On the other hand,  Claim \ref{claim:a1} implies
$l(T_{\owe_C(a-1)}(\beta))< l(\beta)$ and
we complete the proof by induction on $l(\beta)$.
\end{proof}

\section{Proof of Proposition \ref{proposition:step -1 of A_n}}\label{section:preliminary}

Our main purpose in this section is to show Proposition \ref{proposition:step -1 of A_n}.
As explained in Introduction, the essential part is to find $\Psi \in B$ such that
$l(\Psi(\alpha)) < l(\alpha)$ for a spherical object $\alpha \in D_Z(X)$ with $l(\alpha)>1$.
In Lemma \ref{sheaf on A_n} of \S \ref{subsection:Grothendieck}, we clarify the structure of an $\owe_Z$-module of pure dimension $1$,
generalizing a well-known theorem of  Grothendieck.
This gives an expression of cohomology sheaves of a spherical object $\alpha \in D_Z(X)$
in a computable way.
Then using results in \S \ref{section:tool} and in \S \ref{subsection:Grothendieck},
we show Lemma A in \S \ref{subsection:lemmaA} 
and Lemma B in \S \ref{subsection:lemmaB};
these lemmas provide sufficient conditions for the existence of $\Psi \in B$ as above.
Finally, we show in \S \ref{subsection:step 1 of A_n}
that we can always apply Lemma A or B, and thus obtain Proposition
\ref{proposition:step -1 of A_n}.


\subsection{Generalization of a theorem of Grothendieck}\label{subsection:Grothendieck}
Grothendieck proved that every vector bundle on a smooth rational curve 
decomposes into a direct sum of line bundles. 
We generalize this result in the case of a chain of smooth rational curves.

We first introduce some notation that we use in the statement and in the proof.
Let $Z=\bigcup _{i=1}^nC_i$ be a chain of smooth rational curves $C_i$.
We denote by $\Sigma(Z)$ the set of the isomorphism classes of
sheaves $\owe_{C_s \cup \dots \cup C_t}(a_s,\dots,a_t)$,
where $1 \le s \le t \le n$ and $a_s, \dots, a_t \in \ZZ$.
$\Sigma_{C_1}(Z) \subset \Sigma(Z)$ is the subset consisting of
$\mc{R} \in \Sigma(Z)$ with $\Supp \mc{R} \supset C_1$.
We define the {\it lexicographic order} on $\Sigma_{C_1}(Z)$ by setting
$$
\owe_{C_1 \cup \dots \cup C_s}(a_1, \dots, a_s)
>\owe_{C_1 \cup \dots \cup C_t}(b_1, \dots, b_t)
$$
if either of the following holds.
\begin{itemize}
\item For some integer $k$ ($1 \le k \le s, t$),
we have $a_i = b_i$ ($1 \le i \le k-1$)
and $a_k > b_k$.
\item We have $s<t$ and $a_i=b_i$($1 \le i \le s$).
\end{itemize}
Let $x \in C_1 \setminus (C_1 \cap C_2)$ be a point.
Then we can see that for $\mc{R}, \mc{S} \in \Sigma_{C_1}(Z)$,
the inequality $\mc{R} \le \mc{S}$ holds if and only if the restriction map
$$
\Hom_Z(\mc{R}, \mc{S}) \to \Hom_{\CC}(\mc{R}|_x, \mc{S}|_x)
$$
is non-zero.


\begin{lemma}\label{sheaf on A_n} 
Let $Z=\bigcup _{i=1}^nC_i$ be a chain of smooth rational curves $C_i$ and
let $\mathcal{E}$ be a coherent $\mathcal{O}_Z$-module, pure of dimension $1$.
Then $\mathcal{E}$ decomposes into a direct sum of sheaves in $\Sigma(Z)$.
Moreover, such a decomposition is unique up to isomorphism.
\end{lemma}

\begin{proof}
The case $n=1$ is due to Grothendieck, so we consider the case $n\ge 2$.
We define
$$
l(\mathcal{E})=\rank \mathcal{E}|_{C_1}+\cdots+ \rank \mathcal{E}|_{C_n}
$$
and use induction on $l(\mc{E})$. 
We may assume that $\Supp\mathcal{E}$ contains $C_1$.
Replacing $\mathcal{E}$ with 
$\mathcal{E}\otimes\mathcal{L}$ for some line bundle $\mc{L}$ on $Z$,
we may also assume that 
$\Hom _Z^0(\mathcal{E}, \mathcal{O}_{C_1})\ne 0$ and $\Hom _Z^0(\mathcal{E}, \mathcal{O}_{C_1}(-1))= 0$.
Then there exists an exact sequence
$$
0\to\mathcal{E}' \to \mathcal{E}\to\mathcal{O}_{C_1}\to 0,
$$
where $\mc{E}'$ is an $\owe_Z$-module of pure dimension $1$.
By the induction hypothesis, we can decompose $\mc{E}'$
into sheaves in $\Sigma(Z)$.
We write
$$
\mathcal{E}'=\bigoplus _i\mathcal{E}_i\oplus\bigoplus _i\mathcal{F}_i \oplus\bigoplus _i\mathcal{G}_i,
$$ 
where $\mc{E}_i \in \Sigma_{C_1}(Z)$, $\mc{F}_i \in \Sigma_{C_2}(C_2 \cup \dots \cup C_n)$ and $\mc{G}_i \in \Sigma (C_3 \cup \dots \cup C_n)$.
It follows from $\mc{E}_i \in \Sigma_{C_1}(Z)$ that
$$
\Ext _Z^1(\mathcal{O}_{C_1},\mathcal{E}_i)
\cong H^1(\hom_Z(\mathcal{O}_{C_1},\mathcal{E}_i)),
$$
which is zero by $\Hom _Z^0(\mathcal{E}', \mathcal{O}_{C_1}(-1))=0$.
Therefore, we have
$$
\mc{E} \cong \mc{K} \oplus \bigoplus_i \mc{E}_i \oplus \bigoplus_i \mc{G}_i,
$$
where $\mc{K}$ is given by an extension
\begin{equation}\label{equation:extensionK}
0 \to \bigoplus_i \mc{F}_i \to \mc{K} \to \owe_{C_1} \to 0.
\end{equation}
Let $e=\oplus e_i \in \Ext^1_Z(\owe_{C_1}, \bigoplus _i \mc{F}_i)$
be the class corresponding to this extension.
If $e=0$, then \eqref{equation:extensionK} splits and consequently
$\mc{E}$ has a desired decomposition.
Thus we may assume $e \ne 0$.
We reorder the indices $i$ of $\mc{F}_i$ so that if $i>j$,
then $\mc{F}_i \ge \mc{F}_j$ holds with respect to the lexicographic order in
$\Sigma_{C_2}(C_2 \cup \dots \cup C_n)$.
Then the image of the restriction map
$$
\Aut_{\owe_Z} \left(\bigoplus _i^r\mathcal{F}_i\right)\to
\Aut_{\CC} \left(\bigoplus _i^r\mathcal{F}_i|_y\right) \cong \GL (r,\CC),
$$
at the point $y\in C_1\cap C_2$,
contains every lower triangular matrix in
$\GL(r, \CC)$.
Since $\Aut_{\owe_Z} (\bigoplus _i^r\mathcal{F}_i)$ acts on
$$
\Ext _Z^1(\mathcal{O}_{C_1},\bigoplus^r_i\mathcal{F}_i)\cong \CC ^r
$$
thorough the natural action of 
$\GL(r, \CC)$,
there is an element $g \in \Aut_{\owe_Z}(\bigoplus_i^r \mc{F}_i)$ such that
if we put $g\cdot e = \bigoplus_i e'_i$, then $e'_i=0$ except for one index $i=i_0$.
Let $\mc{F}'_{i_0}$ be the unique non-trivial extension of $\owe_{C_1}$ by $\mc{F}_{i_0}$.
Then $\mc{F}_{i_0}$ belongs to $\Sigma(Z)$ and there is an isomorphism
$$
\mc{K} \cong \bigoplus_{i \ne i_0} \mc{F}_i \oplus \mc{F}'_{i_0},
$$
which proves the existence part of the lemma.

For the uniqueness, fix a point $x \in C_1 \setminus (C_1 \cap C_2)$
and let $\mc{R} \in \Sigma_{C_1}(Z)$ be the maximum element
that has the property that the restriction map
$$
\eta: \Hom_Z(\mc{R}, \mc{E}) \to \Hom_{\CC}(\mc{R}|_x, \mc{E}|_x)
$$
is non-zero.
We denote by $r$ the rank of the linear map $\eta$.
Then, in any decomposition of $\mc{E}$ as in the lemma,
$\mc{E}$ contains exactly $r$ copies of $\mc{R}$ as direct summands.
We fix such a decomposition and write $\mc{E} = \mc{E}_1 \oplus \mc{E}_2$
with $\mc{E}_1 = \CC^r \otimes \mc{R}$.
For another such decomposition $\mc{E}=\mc{E}'_1 \oplus \mc{E}'_2$,
$V:= \Hom_Z(\mc{R}, \mc{E}'_1) \subset \Hom_Z(\mc{R}, \mc{E})$
is an $r$-dimensional subspace
such that the restriction $\eta|_V$ is an isomorphism
to the image of $\eta$.
Then the composite of the evaluation map
$\ev_V : V \otimes \mc{R} \to \mc{E}$ and the projection $\mc{E} \to \mc{E}_1$
is an isomorphism.
Since the image of $\ev_V$ is $\mc{E}'_1$,
this proves $\mc{E}_2 \cong \mc{E}_2'$
and completes the proof by induction on $l(\mc{E})$.
\end{proof}

Lemma \ref{sheaf on A_n} provides an explicit form of an $\mc{O}_Z$-module
of pure dimension $1$. Our proofs of Propositions \ref{proposition:step -1 of A_n}
and \ref{proposition:step -2 of A_n} heavily use this explicit form. When $Z$ forms 
a $D_n$- or an $E_n$-configuration, we cannot directly generalize Lemma \ref{sheaf on A_n};
a purely one-dimensional sheaf on $Z$ (even with respect to the reduced induced structure) is
not necessarily a direct sum of line bundles on its subtrees.

Till the end of this section,
$Z$ and $X$ denote the varieties as in Introduction, namely,  
$X$ is the minimal resolution of an $A_n$-singularity
$$
Y=\Spec \CC[[x, y, z]]/(x^2+y^2+z^{n+1})
$$
and $Z$ is the exceptional locus of it with reduced induced structure.
  
Suppose that a spherical object $\alpha\in D_Z(X)$ is given. 
Then Corollary \ref{rigid & pure} and Lemma \ref{sheaf on A_n} say that
every cohomology sheaf $\mathcal{H}^p(\alpha)$ can be written as 
$$
\mathcal{H}^p(\alpha)=\mathcal{R}^p_1\oplus\cdots\oplus\mathcal{R}^p_{k_p},
$$
where every $\mathcal{R}^p_l$ $(1\le l\le k_p)$ belongs to $\Sigma(Z)$.
Note that 
\begin{equation}\label{Ext^1}
\Ext^1_X(\mc{R}^p_{l}, \mc{R}^{p'}_{m})=0
\end{equation}
for all $p,p',l,m$ by Corollary \ref{rigid & pure}. 
For example, \eqref{Ext^1} yields  
$$
|\deg_C \mc{R}^p_{l} -\deg_C \mc{R}^{p'}_{m}|\le 1
$$
for any $-2$-curve $C\subset\Supp \mc{R}^p_{l}\cap \Supp\mc{R}^{p'}_{m}$.
We have another application of \eqref{Ext^1}, which is useful later.
In the expression $\bigoplus_p \cohom^p(\alpha) = \bigoplus_j \mc{R}_j$
with $\mc{R}_j \in \Sigma(Z)$, we always assume that $\mc{R}_j$
is a direct summand of $\cohom^p(\alpha)$ for some $p$.

\begin{lemma}\label{lemma:split}
Let $\alpha \in D_Z(X)$ be a spherical object. Suppose that we have a decomposition
$$
\bigoplus_p \mc{H}^p(\alpha)=\bigoplus_j^{r_1}\mc{R}_{1,j}\oplus\bigoplus_j^{r_2}\mc{R}_{2,j}
$$
with $\mc{R}_{k,j} \in \Sigma(Z)$
such that
$$
\chi(\mc{R}_{1,i}, \mc{R}_{2,j})=0
$$
for all $i,j$.
Then either $r_1$ or $r_2$ is zero.
\end{lemma}
\begin{proof}
The vanishing of $\chi(\mc{R}_{1,i}, \mc{R}_{2,j})$ and \eqref{Ext^1} implies the
vanishing of $\Ext^p_X(\mc{R}_{1,i}, \mc{R}_{2,j})$ for all $p$.
Especially, we have
$$
\Ext^2_X(\mc{R}_{1,i}, \mc{R}_{2,j})=\Ext^2_X(\mc{R}_{2,j}, \mc{R}_{1,i})=0
$$
for all $i,j$.
Then, $\alpha$ splits as in Lemma \ref{decomposable}.
Since $\alpha$ is spherical, we obtain the assertion.
\end{proof}

To obtain Proposition \ref{proposition:step -1 of A_n}, as we explain in Introduction, 
we find an autoequivalence $\Psi\in B$ such that $l(\alpha)>l(\Psi (\alpha))$,
assuming $l(\alpha)>1$. 
For this purpose it suffices to find $\Psi\in B$  such that $\sum_p l(\Psi (\cohom ^p(\alpha)))<l(\alpha)$
by Lemma \ref{strategy}. 


\subsection{Lemma A: a case where we can reduce $l(\alpha)$}\label{subsection:lemmaA}
As a first candidate for $\Psi \in B$ with $l(\Psi(\alpha)) < l(\alpha)$,
we consider functors of the form $T_{\owe_{C_i}(a)}$.
We start with an easy but fundamental case.


\begin{lemma}\label{lemma:fundamental}
Let $\alpha \in D_Z(X)$ be a spherical object and $C \subset Z$ a $-2$-curve.
Assume that for every $p$ we have a decomposition
$$
\mathcal{H}^p(\alpha)=\bigoplus_j^{r_1^p}\mathcal{R}_{1,j}^p
\oplus\bigoplus_j^{r_2^p}\mathcal{R}_{2,j}^p
\oplus\bigoplus_j^{r_3^p}\mathcal{R}_{3,j}^p
\oplus\bigoplus_j^{r_4^p}\mathcal{R}_{4,j}^p
\oplus \mc{S}^p,
$$
where $\mathcal{R}^p_{k,j}$'s are sheaves of the forms
\[ 
\xymatrix@R=1ex@M=.3ex{
                           & & C &  \\ 
          \mathcal{R}_{1,j}^p:& \ar@{-}[r]&{\zero}\ar@{-}[r]&\\
          \mathcal{R}_{2,j}^p:& &{\zero}\ar@{-}[r] &\\
          \mathcal{R}_{3,j}^p:& &{\mone}\ar@{-}[r] &\\
          \mathcal{R}_{4,j}^p:& &{\mone} &
         }
\]
and where $\Supp \mc{S}^p \cap C = \emptyset$.
In this situation, we have the following:
\begin{enumerate}
\item If $\sum_p r_2^p > \sum_p r_3^p$,
then $l(T_{\owe_C(-1)}(\alpha)) < l(\alpha)$.
\item If $\sum_p r_2^p < \sum_p r_3^p$,
then $l(T_{\owe_C(-2)}(\alpha)) < l(\alpha)$.
\end{enumerate}
\end{lemma}

\begin{proof}
Combining the assumption of (i) with Lemma \ref{calculation}, 
we deduce that 
$$
\sum_p l(T_{\owe_C(-1)}(\cohom ^p(\alpha)))<\sum_p l(\cohom ^p(\alpha)),
$$
and then obtain the conclusion from
Lemma \ref{strategy}.
(ii) can be seen in a similar way.
\end{proof}
We cannot always find $C$ as above with $\sum_p r_2^p \ne \sum_p r_3^p$ (see Example \ref{A_5} (ii))
and it is important to consider the case $\sum_p r_2^p = \sum_p r_3^p$.


\begin{lemma}\label{lemma:equality}
Let $\alpha \in D_Z(X)$ be a spherical object and $C \subset Z$
a $-2$-curve.
Assume that for every $p$ we have 
$$
\cohom^p(\alpha) \cong \bigoplus_j^{r_2^p} \mc{R}_{2,j}^p \oplus
\bigoplus_j^{r_3^p} \mc{R}_{3,j}^p \oplus \mc{S}^p
$$
with the properties
\begin{itemize}
\item $\mathcal{R}_{2,j}^p$ and $\mathcal{R}_{3,j}^p$ are as in the previous lemma, and
\item $\mc{S}^p$'s are sheaves satisfying that the composition maps
\begin{align*}
\Hom_X(\owe_C(-1), \mc{R}_{2,j}^p) \times \Hom_X(\mc{R}_{2,j}^p, \mc{S}^q) &\to
\Hom_X(\owe_C(-1), \mc{S}^q)\\
\Hom_X(\mc{S}^p, \mc{R}_{3,j}^q) \times \Hom_X(\mc{R}_{3,j}^q, \owe_C(-1)) &\to
\Hom_X(\mc{S}^p, \owe_C(-1))
\end{align*}
are zero for all $p,q,j$.
\end{itemize}
Then, we have either $r_2^p \le r_3^{p-1}$ for all $p$ or
$r_2^p \ge r_3^{p-1}$ for all $p$.
Especially, if $\sum_p r_2^p = \sum_p r_3^p$, then
the equality $r_2^p = r_3^{p-1}$ holds for every $p$.
\end{lemma}
\begin{proof}
Put $\mc{R}_2^p=\bigoplus_j \mc{R}_{2,j}^p$ and $\mc{R}_3^p=\bigoplus_j \mc{R}_{3,j}^p$.
Let $e^p(\alpha) \in \Ext_X^2(\cohom^p(\alpha), \cohom^{p-1}(\alpha))$ be
the class determined by $\alpha$ as in \S \ref{section:spectral}.
According to the decomposition
$$
\cohom^p(\alpha)= \mc{R}_2^p \oplus \mc{R}_3^p \oplus \mc{S}^p,
$$
$e^p(\alpha)$ also decomposes and determines classes
\begin{align*}
\eta^p &\in \Ext_X^2(\mc{R}_2^p, \mc{R}_3^{p-1}),\\
\xi^p &\in \Ext_X^2(\mc{R}_3^p, \mc{R}_3^{p-1}), \\
\psi^p &\in \Ext_X^2(\mc{S}^p, \mc{R}_3^{p-1}).
\end{align*}
We denote by
$\bar\eta^p \in \Ext_X^2(\owe_C(-1)^{\oplus r_2^p}, \owe_C(-1)^{\oplus r_3^{p-1}})$
the following composite:
\begin{equation}\label{equation:etabar}
\xymatrix{
{\owe_C(-1)^{\oplus r_2^p}} \ar[rrr]^{\bar\eta^p} \ar[d]^{\cong}
&&&
{\owe_C(-1)^{\oplus r_3^{p-1}}[2]} \ar[d]^{\cong}
\\
{\hom_X(\owe_C, \mc{R}_2^p)} \ar@{^{(}->}[r]
&
{\mc{R}_2^p} \ar[r]^{\eta^p}
&
\mc{R}_3^{p-1}[2] \ar@{->>}[r]
&
\mc{R}_3^{p-1}|_C[2]
}
\end{equation}
Assume that the first assertion does not hold.
Then, there are $i, j$ with $r_2^i< r_3^{i-1}$ and $r_2^j > r_3^{j-1}$.
It follows from $r_2^i< r_3^{i-1}$ that there is a surjection
$\gamma: \mc{R}_3^{i-1}|_C \to \owe_C(-1)$ with
$\gamma \circ \bar\eta^i=0$.
Similarly, we have an injection $\delta:\owe_C(-1) \hookrightarrow \hom_X(\owe_C, \mc{R}_2^j)$
with $\bar\eta^j \circ \delta =0$.
Let $f: \cohom^{i-1}(\alpha) \to \cohom^j(\alpha)$ be the following composite:
$$
\cohom^{i-1}(\alpha) \twoheadrightarrow \mc{R}_3^{i-1}|_C \overset{\gamma}{\twoheadrightarrow}
\owe_C(-1) \overset{\delta}{\hookrightarrow}
\hom_X(\owe_C, \mc{R}_2^j)
\hookrightarrow \cohom^j(\alpha)
$$
We claim that
$f \circ e^i(\alpha)=0$ in $\Ext_X^2(\cohom^i(\alpha), \cohom^j(\alpha))$.
Let ${\bar f}$ and $p$ be as follows:
$$
\xymatrix{
{\mc{R}_2^i} \ar[dr]^(.3){\eta^i}^(.6){[2]} & {\cohom^{i-1}(\alpha)}\ar@{->>}[d] \ar[rr]^f && {\cohom^{j}(\alpha)} \\
{\mc{R}_3^i} \ar[r]^(.3){\xi^i}^(.6){[2]} & {\mc{R}_3^{i-1}} \ar@{->>}[d]^p \ar[rr]^{\bar f} && \mc{R}_2^j \ar@{^{(}->}[u] \\
{\mc{S}^i} \ar[ur]^(.3){\psi^i}^(.6){[2]} & {\mc{R}_3^{i-1}|_C} \ar@{->>}[r]^{\gamma} & {\owe_C(-1)}
\ar@{^{(}->}[r]^{\delta} & {\hom_X(\owe_C, \mc{R}_2^j)} \ar@{^{(}->}[u]
}
$$
It suffices to show $\bar f \circ \eta^i$, $\bar f \circ \xi^i$ and $\bar f \circ \psi^i$ are all zero.
Since $\gamma \circ \bar\eta^i =0$, we have $\gamma \circ p \circ \eta^i=0$
and therefore $\bar f \circ \eta^i=0$.
$\bar f \circ \xi^i$ factors through
$\gamma \circ p \circ \xi^i \in \Ext_X^2(\mc{R}_3^i, \owe_C(-1))=0$ and hence is zero.
Finally, $\bar f \circ \psi^i \in \Ext_X^2(\mc{S}^i, \mc{R}_2^j)$ is in the image of the
composition map
$$
\Ext_X^2(\mc{S}^i, \owe_C(-1)) \times \Hom_X(\owe_C(-1), \mc{R}_2^j) \to
\Ext_X^2(\mc{S}^i, \mc{R}_2^j)
$$
which is zero by the assumption and the Serre duality.
Thus we showed the claim.
Similarly, we have $e^j(\alpha) \circ f=0$.

Therefore, in the spectral sequence \eqref{equation:spectral1} (for $\alpha=\beta$),
$$
f \in \Hom_X(\cohom^i(\alpha), \cohom^j(\alpha)) \subset E_2^{0, j-i}
$$
lies in the kernel of $d_2^{0, j-i}$.
This contradicts Proposition \ref{TFAE for spherical}.
\end{proof}
The above proof is actually showing a slightly stronger statement:


\begin{lemma}\label{lemma:invertible}
Under the assumption of the above lemma, write $\bar\eta^p=M_p \otimes_{\CC} e$, where
$\bar\eta^p$ is defined in \eqref{equation:etabar}, $M_p$ is an $r_3^{p-1} \times r_2^p$ matrix
and $e \in \Ext_X^2(\owe_C(-1), \owe_C(-1)) \cong \CC$ is a fixed basis.
Then, we have either $\rank M_p = r_2^p$ for all $p$ or $\rank M_p = r_3^{p-1}$ for all $p$.
Especially, if $\sum_p r_2^p = \sum_p r_3^p$, then all $M_p$ are invertible.
\end{lemma}

Now we go back to the situation in Lemma \ref{lemma:fundamental}.


\begin{lemma}\label{lemma:r_4=0}
Under the assumptions of Lemma \ref{lemma:fundamental},
assume the equality $\sum_p r_2^p = \sum_p r_3^p \ne 0$ holds.
Then $\mc{R}_4^p=0$ for all $p$.
\end{lemma}

\begin{proof}
Put $\mc{R}_k^p=\bigoplus_j \mc{R}_{k,j}^p$ and
write $e^p(\alpha)=(e_{ij}^p)$, where $e_{ij}^p \in \Ext_X^2(\mc{R}_j^p, \mc{R}_i^{p-1})$.
Among these entries, $e_{24}^p, e_{43}^p, e_{41}^p, e_{14}^p$ are zero
because the corresponding Ext groups vanish.
If, in addition, $e_{34}^p$ and $e_{42}^p$ are zero, we have objects $\alpha_1, \alpha_2$ such that 
$\alpha \cong \alpha_1 \oplus \alpha_2$ with $\cohom^p(\alpha_1) \cong \mc{R}_1^p \oplus \mc{R}_2^p \oplus \mc{R}_3^p$
and $\cohom^p(\alpha_2) \cong \mc{R}_4^p$ by Lemma \ref{decomposable}.
Since $\alpha$ is spherical, either $\alpha_1$ or $\alpha_2$ must be zero and we are done.
Thus it is enough to show that $e_{34}^p$ and $e_{42}^p$ become zero if we change the decomposition
$$
\cohom^p(\alpha)=\mc{R}_1^p \oplus \mc{R}_2^p \oplus \mc{R}_3^p \oplus \mc{R}_4^p
$$
by suitable automorphisms of $\cohom^p(\alpha)$.
$e_{34}^p$ lies in 
$$
\Ext_X^2(\mc{R}_4^p, \mc{R}_3^{p-1}) \cong
\Hom_{\CC}(\CC^{r_4^p}, \CC^{r_3^{p-1}}) \otimes_{\CC} \Ext_X^2(\owe_C(-1), \owe_C(-1))
$$
and hence is of the form $A_p \otimes e$ for an $r_3^{p-1} \times r_4^p$ matrix $A_p$ and
the same $e$ as
in Lemma \ref{lemma:invertible}.
Lemma \ref{lemma:invertible} applied to
$\mc{S}^p= \mc{R}_1^p \oplus \mc{R}_4^p$
says that $e_{32}^p = \eta^p$ determines
$\bar\eta^p = M_p \otimes e$ with $M_p$ an invertible matrix.
We determine an automorphism $g^p=(g^p_{ij})$ of $\cohom^p(\alpha)$
by
$$
g_{24}^p=-M_p^{-1}A_p \in \Hom_X(\mc{R}_4^p, \mc{R}_2^p) \cong \Hom_{\CC}(\CC^{r_4^p}, \CC^{r_2^p})
$$
and $g_{ij}^p=\delta_{ij}I_{\mc{R}_i^p}$ for the other $(i,j)$.
If we replace $e^p(\alpha)$ by $(g^{p-1})^{-1} e^p(\alpha) g^p$,
then $e_{34}^p$ becomes zero and $e_{42}^p$ does not change.
$e_{42}^p$ is also of the form $B_p \otimes e$ for a matrix $B_p$ and in a similar way
we can find automorphisms that eliminate $e_{42}^p$ without changing $e_{34}^p$.
\end{proof}


\begin{lemmaA}\label{lemmaA}
Let $\alpha \in D_Z(X)$ be a spherical object and let $C \subset Z$ be a $-2$-curve.
Assume that we can write
$$
\bigoplus_p\mathcal{H}^p(\alpha)=
\bigoplus_j^{r_1}\mathcal{R}_{1,j}
\oplus\bigoplus_j^{r_2}\mathcal{R}_{2,j}
\oplus\bigoplus_j^{r_3}\mathcal{R}_{3,j}
\oplus\bigoplus_j^{r_4}\mathcal{R}_{4,j}
\oplus \mc{S}
$$
where $\mathcal{R}_{k,j}$'s are sheaves of the forms
\[ 
\xymatrix@R=1ex@M=.3ex{     & & C &  \\
          \mathcal{R}_{1,j}:& \ar@{-}[r]&{\no}\ar@{-}[r]&\\
          \mathcal{R}_{2,j}:& &{\no}\ar@{-}[r] &\\
          \mathcal{R}_{3,j}:& &{\no} &  \\ 
          \mathcal{R}_{4,j}:&  \ar@{-}[r]&{\no} &    
          }
\]
and where $\Supp \mc{S} \cap C = \emptyset$. 
Suppose that either $r_3\ne 0$ or $r_2\cdot r_4\ne 0$ holds, 
and suppose furthermore that $\Supp \alpha\ne C$.
Then, there is an integer $a$ such that $l(T_{\mc{O}_C(a)} (\alpha))<l(\alpha)$.
\end{lemmaA}
\begin{proof} 
We can freely replace $\alpha$ with $\alpha\otimes \mc{L}$ for some $\mc{L}\in\Pic X$ 
by Lemma \ref{lemma:normal sub}.
Hence we may assume that $\max_{k,j}\deg_C\mc{R}_{k,j}=0$,
and then we have $\deg_C\mc{R}_{k,j}\in\{-1,0\}$ for all $k,j$ by \eqref{Ext^1}.
Note that we have 
$$
\chi(\mathcal{R}_{1,j},\mathcal{R}_{3,i})= \chi(\mc{S}, \mc{R}_{3,i})=0
$$
for any $i,j$. 
Hence if $r_2=r_4=0$ (which implies $r_3\ne 0$ by our assumption), 
then we get $\bigoplus_j^{r_1}\mathcal{R}_{1,j}\oplus \mc{S}=0$ by Lemma \ref{lemma:split}.
This contradicts our assumption that $\Supp\alpha\ne C$.
Therefore, because the condition is symmetric, we may assume $r_2\ne 0$.

When $r_2\cdot r_4\ne 0$ holds, we see from \eqref{Ext^1} 
that 
$$
\deg_C\mc{R}_{2,i}=\deg_C\mc{R}_{4,k}=a
$$
for a fixed $a \in \{-1, 0\}$ and for all $i,k$,
and that $\deg_C \mc{R}_{3,j}$ is $a$ or $a-1$.
Then 
$l(T_{\mc{O}_C(a-1)}(\alpha))<l(\alpha)$ holds as desired. 

Next consider the case $r_2\cdot r_3\ne 0$ and $r_4=0$.
If $\deg_C\mc{R}_{3,j}=-1$ for all $j$, 
Lemma \ref{lemma:fundamental} and Lemma \ref{lemma:r_4=0} imply the conclusion. 
Hence suppose $\deg_C\mc{R}_{3,j}=0$ for some $j$. 
Then $\deg_C\mc{R}_{2,j}=0$ for all $j$ by \eqref{Ext^1}, 
and so $l(T_{\mc{O}_C(-1)}(\alpha))<l(\alpha)$ holds, as required.
\end{proof}


\subsection{Lemma B: another case where we can reduce $l(\alpha)$
}\label{subsection:lemmaB}


\begin{lemma}\label{C_s&C_t}
Let $\alpha \in D_Z(X)$ be a spherical object and $W=C_s \cup \dots \cup C_t \subset Z$
a chain of $-2$-curves with $s<t$.
Assume that for every $p$ we have 
$$
\mathcal{H}^p(\alpha)=
\bigoplus_j^{r_1^p}\mathcal{R}_{1,j}^p
\oplus\bigoplus_j^{r_2^p}\mathcal{R}_{2,j}^p
\oplus\bigoplus_j^{r_3^p}\mathcal{R}_{3,j}^p
\oplus\bigoplus_j^{r_4^p}\mathcal{R}_{4,j}^p
\oplus\bigoplus_j^{r_5^p}\mathcal{R}_{5,j}^p
\oplus \mc{S}^p,
$$
where $\mathcal{R}_{k,j}^p$'s are sheaves of the forms
\[ 
\xymatrix@R=1ex@M=.3ex{
& & C_s & C_{s+1} & & C_{t-1} & C_t & \\
\mathcal{R}_{1,j}^p:& \ar@{-}[r] &{\zero}\ar@{-}[r]&{\zero}\ar@{-}[r]&\cdots \ar@{-}[r]&{\zero}\ar@{-}[r]&{\zero}\ar@{-}[r]& \\
\mathcal{R}_{2,j}^p:& &{\zero}\ar@{-}[r]&{\zero}\ar@{-}[r]&\cdots \ar@{-}[r]&{\zero}\ar@{-}[r]&{\zero}\ar@{-}[r]& \\
\mathcal{R}_{3,j}^p:& &{\mone}\ar@{-}[r]&{\zero}\ar@{-}[r]&\cdots \ar@{-}[r]&{\zero}\ar@{-}[r]&{\zero}\ar@{-}[r]& \\
\mathcal{R}_{4,j}^p:& &{\zero}\ar@{-}[r]&{\zero}\ar@{-}[r]&\cdots \ar@{-}[r]&{\zero}\ar@{-}[r]&{\mone} & \\
\mathcal{R}_{5,j}^p:& &{\mone}\ar@{-}[r]&{\zero}\ar@{-}[r]&\cdots \ar@{-}[r]&{\zero}\ar@{-}[r]&{\zero}& 
         }
\]
and where $\Supp \mc{S}^p \cap W = \emptyset$.
Under these assumptions, either of the following holds:
\begin{enumerate}
\item At least one of $l(T_{\owe_{C_s}(-1)}(\alpha))$, $l(T_{\owe_{C_s}(-2)}(\alpha))$,
$l(T_{\owe_{C_t}(-1)}(\alpha))$ or $l(T_{\owe_{C_t}(-2)}(\alpha))$ is smaller than $l(\alpha)$
or
\item $r_4^p = r_5^p=0$ for all $p$.
\end{enumerate}
\end{lemma}
\begin{proof}
Assume that (i) does not hold.
Then Lemma \ref{lemma:fundamental} applied to $C=C_s$ and $C=C_t$  imply
$\sum_p(r_2^p+r_4^p)= \sum_p(r_3^p + r_5^p)$
and $\sum_p r_5^p = \sum_pr_4^p$ respectively.
These equalities also deduce $\sum_p r_2^p = \sum_p r_3^p$.
Then, applying Lemma \ref{lemma:equality} in three ways, we obtain
$$
\begin{aligned}
r_2^p + r_4^p &= r_3^{p-1} + r_5^{p-1} \\
r_5^p &= r_4^{p-1} \\
r_2^p &= r_3^{p-1}
\end{aligned}
$$
for all $p$.
Especially, we have both $r_5^p = r_4^{p-1}$ and $r_4^p = r_5^{p-1}$.
Since $\cohom^p(\alpha)=0$ except for finitely many integers $p$,
this means that all $r_4^p$ and $r_5^p$ are zero.
\end{proof}


\begin{lemma}\label{middle}
Let $\alpha\in D_Z(X)$ be a spherical object and fix 
positive integers $s,t$ with $s<t$. Assume 
$$
\bigoplus_p\mathcal{H}^p(\alpha)=
\bigoplus_j^{r_1}\mathcal{R}_{1,j}\oplus\bigoplus_j^{r_2}\mathcal{R}_{2,j}
\oplus\bigoplus_j^{r_3}\mathcal{R}_{3,j}\oplus\bigoplus_j^{r_4}\mathcal{R}_{4,j}
\oplus\bigoplus_j^{r_5}\mathcal{R}_{5,j}\oplus\bigoplus_j^{r_6}\mathcal{R}_{6,j}
\oplus\mathcal{S},
$$
where $\mathcal{R}_{k,j}$'s are sheaves of the forms 
\[ 
\xymatrix@R=1ex@M=.3ex{
                           & & C_s & C_{s+1} & & C_{t-1} & C_t & \\
          \mathcal{R}_{1,j}:&\ar@{-}[r]&{\zero}\ar@{-}[r]&{\zero}\ar@{-}[r]&\cdots \ar@{-}[r]&{\zero}\ar@{-}[r]&{\zero}\ar@{-}[r]& \\
          \mathcal{R}_{2,j}:&\ar@{-}[r]&{\zero}\ar@{-}[r]&{\zero}\ar@{-}[r]&\cdots \ar@{-}[r]&{\zero}\ar@{-}[r]&{\mone}\ar@{-}[r]& \\
          \mathcal{R}_{3,j}:&          &{\zero}\ar@{-}[r]&{\zero}\ar@{-}[r]&\cdots \ar@{-}[r]&{\zero}\ar@{-}[r]&{\zero}\ar@{-}[r]& \\
          \mathcal{R}_{4,j}:&          &{\zero}\ar@{-}[r]&{\zero}\ar@{-}[r]&\cdots \ar@{-}[r]&{\zero}\ar@{-}[r]&{\mone}\ar@{-}[r]& \\
          \mathcal{R}_{5,j}:&          &{\mone}\ar@{-}[r]&{\zero}\ar@{-}[r]&\cdots \ar@{-}[r]&{\zero}\ar@{-}[r]&{\zero}\ar@{-}[r]& \\
          \mathcal{R}_{6,j}:&          &{\mone}\ar@{-}[r]&{\zero}\ar@{-}[r]&\cdots \ar@{-}[r]&{\zero}\ar@{-}[r]&{\mone}\ar@{-}[r]&  
         }
\]
and where $\Supp\mathcal{S}\cap (C_s\cup\cdots\cup C_t)=\emptyset$.
Suppose that 
\begin{equation}\label{equation:assumption1}
\mbox{$l(\alpha)\le l(\Phi(\alpha))$ for all $\Phi\in \Span{T_{\owe_{C_l}(a)}\bigm|  a\in \ZZ, s\le l\le t }$}
\end{equation} 
and $r_3+r_4+r_5+r_6\ne 0$. Then we have $r_1=r_3=r_5=0$ or $r_2=r_4=r_6=0$.
In particular, $\deg_{C_t}\mathcal{R}_{k,j}$ does not depend on $j$ and $k$.
\end{lemma}

\begin{proof}
First note that $r_3\cdot r_6=0$ by \eqref{Ext^1}. We prove the following:
\begin{itemize}
\item If $r_3=0$, then we have $r_1=r_5=0$.
\item If $r_6=0$, then we have $r_2=r_4=0$.
\end{itemize}
First assume that $r_3=0$.
We apply Lemma \ref{lemma:fundamental} for $C=C_s$ and then obtain
\begin{equation}\label{C=C_s}
r_4=r_5+r_6
\end{equation} 
from the assumption \eqref{equation:assumption1}.
Put 
$$
\Phi=T_{\owe_{C_{s+1}}(-1)}\circ \dots \circ T_{\owe_{C_{t-1}}(-1)}\circ T_{\owe_{C_t}(-2)} 
$$
if $t>s+1$, and 
$$
\Phi=T_{\owe_{C_{s+1}}(-2)}
$$
if $t=s+1$.
Then $\Phi(\mathcal{R}_{i,j})$ are sheaves,
and we have $\deg_{C_s}\Phi(\mathcal{R}_{4,j})=\deg_{C_s}\Phi(\mathcal{R}_{5,j})=0$
and $\deg_{C_s}\Phi(\mathcal{R}_{6,j})=-1$.
If $r_5\ne 0$, then we see from \eqref{C=C_s} that $r_4+r_5>r_6$ and then from Lemma \ref{lemma:fundamental} that
$l(T_{\owe_{C_s}(-1)}\circ\Phi(\alpha))<l(\Phi(\alpha))=l(\alpha)$,
a contradiction to \eqref{equation:assumption1}.
If $r_1\ne 0$, we have $r_6=0$ by \eqref{Ext^1} and again $r_4+r_5 > r_6$. 
This contradicts \eqref{equation:assumption1} as above. 

In the case $r_6=0$, we get the assertion by a similar argument, using 
$$
\Psi=T'_{\owe_{C_{s+1}}(-1)}\circ \dots \circ T'_{\owe_{C_{t-1}}(-1)}\circ T'_{\owe_{C_t}(-1)}, 
$$
instead of $\Phi$.
\end{proof}

The above proof teaches us how to reduce $l(\alpha)$ for the spherical object $\alpha$ in Example \ref{A_5} (ii);
we can see that 
$$
l(T_{\mc{O}_{C_1}(-1)}\circ T_{\mc{O}_{C_2}(-2)}(\alpha))<l (\alpha).
$$ 
On the other hand, note that
$$
l(T_{\mc{O}_{C_l}(a)}(\alpha))\ge l(\alpha),\quad l(T'_{\mc{O}_{C_l}(a)}(\alpha))\ge l(\alpha)
$$ 
for any $a,l\in \ZZ$ $(1\le l\le 5)$ in the same example.


\begin{lemmaB}\label{lemmaB}
Let $\alpha\in D_Z(X)$ be a spherical object and fix positive integers $s,t$ with $s<t$. 
Assume that we can write
$$
\bigoplus_p\mathcal{H}^p(\alpha)=\bigoplus_j^{r_1}\mathcal{R}_{1,j}\oplus\bigoplus_j^{r_2}\mathcal{R}_{2,j}
\oplus\bigoplus_j^{r_3}\mathcal{R}_{3,j}\oplus\bigoplus_j^{r_4}\mathcal{R}_{4,j}\oplus\mathcal{S},
$$
where $\mathcal{R}_{k,j}$'s are sheaves of the forms 
\[ 
\xymatrix@R=1ex@M=.3ex{     &          & C_s & C_{s+1} & & C_{t-1} & C_t & \\
          \mathcal{R}_{1,j}:&\ar@{-}[r]&{\no}\ar@{-}[r]&{\no}\ar@{-}[r]&\cdots \ar@{-}[r]&{\no}\ar@{-}[r]&{\no}\ar@{-}[r]& \\
          \mathcal{R}_{2,j}:&          &{\no}\ar@{-}[r]&{\no}\ar@{-}[r]&\cdots \ar@{-}[r]&{\no}\ar@{-}[r]&{\no}\ar@{-}[r]& \\
          \mathcal{R}_{3,j}:&          &{\no}\ar@{-}[r]&{\no}\ar@{-}[r]&\cdots \ar@{-}[r]&{\no}\ar@{-}[r]&{\no}& \\
          \mathcal{R}_{4,j}:&\ar@{-}[r]&{\no}\ar@{-}[r]&{\no}\ar@{-}[r]&\cdots \ar@{-}[r]&{\no}\ar@{-}[r]&{\no}& }                   
\]
and where $\Supp\mathcal{S}\cap (C_s\cup\cdots\cup C_t)=\emptyset$.
Suppose that either $r_3\ne 0$ or $r_2\cdot r_4\ne 0$ holds.
Then there is  
$$
\Phi\in  \Span{T_{\owe_{C_l}(a)}\bigm|  a\in \ZZ, s\le l\le t }
$$
such that $ l(\Phi(\alpha))<l(\alpha)$.
\end{lemmaB}

\begin{proof}
For a contradiction, we assume
\begin{equation}\label{equation:assumption2}
\mbox{$l(\alpha)\le l(\Phi(\alpha))$ for all $\Phi\in  \Span{T_{\owe_{C_l}(a)}\bigm|  a\in \ZZ, s\le l\le t }$.}
\end{equation} 
Then it is enough to check $r_3=r_4=0$ or $r_2=r_3=0$. 
By Lemma \ref{middle} and by tensoring with a suitable line bundle on $X$ (cf. Lemma \ref{lemma:normal sub}), 
we may assume that $\deg_{C_l}\mathcal{R}_{k,j}=0$ for all $l$ $(s< l < t)$, $k$ and $j$.
Moreover, we assume
$$
\max_{k,j}\deg_{C_s}\mc{R}_{k,j}=\max_{k,j}\deg_{C_t}\mc{R}_{k,j}=0.
$$
Then we see that $\deg_{C_s}\mc{R}_{k,j},\deg_{C_t}\mc{R}_{k,j}\in\{-1,0\}$ 
for all $k,j$ by \eqref{Ext^1}.
We further claim
$$
\deg_{C_s} \mc{R}_{1,j} = \deg_{C_s} \mc{R}_{4,j}=0.
$$
Otherwise, \eqref{Ext^1} implies that
$\deg_{C_s} \mc{R}_{2,j} = \deg_{C_s} \mc{R}_{3,j}=-1$
and hence that $l(T_{\owe_{C_s}(-2)} (\alpha)) \le l(\alpha) -r_2-r_3$;
\eqref{equation:assumption2} shows $r_2=r_3=0$ as desired.
Similarly, we have
$$
\deg_{C_t} \mc{R}_{1,j} = \deg_{C_t} \mc{R}_{2,j}=0.
$$
Thus we can write 
\begin{eqnarray}
\bigoplus_j^{r_2}\mathcal{R}_{2,j}
&=&\bigoplus_j^{s_1}\mathcal{S}_{1,j}\oplus\bigoplus_j^{s_2}\mathcal{S}_{2,j}\nonumber\\
\bigoplus_j^{r_3}\mathcal{R}_{3,j}
&=&\bigoplus_j^{s_3}\mathcal{S}_{3,j}\oplus\bigoplus_j^{s_4}\mathcal{S}_{4,j}                         
\oplus\bigoplus_j^{s_5}\mathcal{S}_{5,j}\oplus\bigoplus_j^{s_6}\mathcal{S}_{6,j},\nonumber\\
\bigoplus_j^{r_4}\mathcal{R}_{4,j}
&=&\bigoplus_j^{s_7}\mathcal{S}_{7,j}\oplus\bigoplus_j^{s_8}\mathcal{S}_{8,j},\nonumber
\end{eqnarray}
where $\mathcal{S}_{k,j}$'s are sheaves of the forms 
in the following figure.
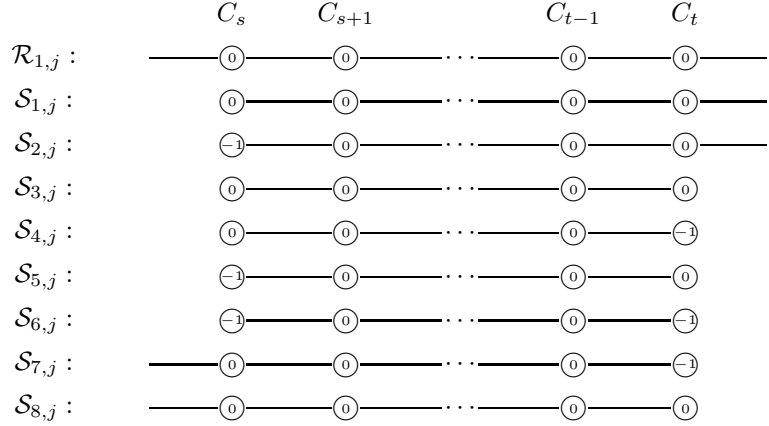
\[ 
\xymatrix@R=1ex@M=.3ex{     &          & C_s & C_{s+1} & & C_{t-1} & C_t & \\
          \mathcal{R}_{1,j}:&\ar@{-}[r]&{\zero}\ar@{-}[r]&{\zero}\ar@{-}[r]&\cdots \ar@{-}[r]&{\zero}\ar@{-}[r]&{\zero}\ar@{-}[r]& \\
          \mathcal{S}_{1,j}:&          &{\zero}\ar@{-}[r]&{\zero}\ar@{-}[r]&\cdots \ar@{-}[r]&{\zero}\ar@{-}[r]&{\zero}\ar@{-}[r]& \\
          \mathcal{S}_{2,j}:&          &{\mone}\ar@{-}[r]&{\zero}\ar@{-}[r]&\cdots \ar@{-}[r]&{\zero}\ar@{-}[r]&{\zero}\ar@{-}[r]& \\ 
          \mathcal{S}_{3,j}:&          &{\zero}\ar@{-}[r]&{\zero}\ar@{-}[r]&\cdots \ar@{-}[r]&{\zero}\ar@{-}[r]&{\zero}&\\
          \mathcal{S}_{4,j}:&          &{\zero}\ar@{-}[r]&{\zero}\ar@{-}[r]&\cdots \ar@{-}[r]&{\zero}\ar@{-}[r]&{\mone}& \\
          \mathcal{S}_{5,j}:&          &{\mone}\ar@{-}[r]&{\zero}\ar@{-}[r]&\cdots \ar@{-}[r]&{\zero}\ar@{-}[r]&{\zero}& \\
          \mathcal{S}_{6,j}:&         &{\mone}\ar@{-}[r]&{\zero}\ar@{-}[r]&\cdots \ar@{-}[r]&{\zero}\ar@{-}[r]&{\mone}& \\ 
          \mathcal{S}_{7,j}:& \ar@{-}[r]&{\zero}\ar@{-}[r]&{\zero}\ar@{-}[r]&\cdots \ar@{-}[r]&{\zero}\ar@{-}[r]&{\mone}& \\
          \mathcal{S}_{8,j}:& \ar@{-}[r]&{\zero}\ar@{-}[r]&{\zero}\ar@{-}[r]&\cdots \ar@{-}[r]&{\zero}\ar@{-}[r]&{\zero}& }
\]
Now, applying Lemma \ref{lemma:fundamental} for $C=C_s$ and $C_t$, we obtain from \eqref{equation:assumption2}
\begin{equation}\label{C_s}
s_1+s_3+s_4=s_2+s_5+s_6
\end{equation} 
and
\begin{equation}\label{C_t}
s_3+s_5+s_8=s_4+s_6+s_7
\end{equation} 
respectively. 

If $s_3\ne 0$, we have $s_2=s_6=s_7=0$ by \eqref{Ext^1}.
Substituting it into (\ref{C_s}) and (\ref{C_t}), we get $s_1+s_3+s_4=s_5$ and $s_3+s_5+s_8=s_4$, which is absurd.
By a similar argument, we also arrive at a contradiction when assuming $s_6\ne 0$. Therefore we obtain $s_3=s_6=0$.

Suppose that $s_1\ne 0$ and $s_8\ne 0$. In this case, we know $s_2=s_7=0$ by \eqref{Ext^1}. 
Then \eqref{C_s} and \eqref{C_t} become $s_1+s_4=s_5$ and $s_5+s_8=s_4$, but this is impossible. 
Next assume that $s_1=s_8=0$. Then (\ref{C_s}) and (\ref{C_t}) imply
that $s_2=s_7=0$ and $s_4=s_5$.
We have seen $r_2=r_4=0$ and thus we
apply Lemma \ref{C_s&C_t} to deduce $s_4=s_5=0$ from \eqref{equation:assumption2},
as desired. 
Finally suppose that precisely one of $s_1$ and $s_8$ is zero.
Because the conditions are symmetric, we may assume that $s_1\ne 0$ and $s_8=0$. Recall that we are in the case
$s_3=s_6=s_7=s_8=0$. Again Lemma \ref{C_s&C_t} and \eqref{equation:assumption2} imply that $s_4=s_5=0$.
\end{proof}


\subsection{Proposition \ref{proposition:step -1 of A_n}
: The main result of \S \ref{section:preliminary}}\label{subsection:step 1 of A_n}
For a spherical object $\alpha\in D_Z(X)$,
let us denote by
$$
\Sigma (\alpha)( \subset \Sigma(Z))
$$
the set of all the indecomposable direct summands of $\bigoplus_i\mathcal{H}^i(\alpha)$ obtained in Lemma \ref{sheaf on A_n}.

Now we are in a position to prove Proposition \ref{proposition:step -1 of A_n}.
In the proof, we freely use the equality  
$$
B=\Span{T_{\owe_{C_l}(a)}\bigm|  a\in \ZZ, 1\le l\le n }
$$
proved in Lemma \ref{B=B'}.\\


\noindent
{\it Proof of Proposition \ref{proposition:step -1 of A_n}.}
Notice that if we show the existence of an autoequivalence $\Phi\in B$
such that $l(\alpha)>l(\Phi(\alpha))$, 
then we can prove the statement by induction on $l(\alpha)$.
We assume $\Supp\alpha=Z=C_1\cup\cdots\cup C_n$.
Recall that the proof is already done for the case $n=1$ (and in particular the case $l(\alpha)=1$) by 
Proposition \ref{proposition:spherical}. Hence we consider the case $n\ge2$. 
Put
$$
l_i(\alpha):= \sum_p \length_{\owe_{X, \eta_i}} \cohom^p(\alpha)_{\eta_i}
$$
for each curve $C_i$ (see Introduction for the notation).
To simplify the argument, we also put $l_0(\alpha)=l_{n+1}(\alpha)=0$.

For $\mathcal{R}\in\Sigma (\alpha)$ with $\Supp \mc{R}=C_k \cup\dots\cup C_l$,
we define $s(\mc{R}):=k$ and  $t(\mc{R}):=l$.
Note that \eqref{Ext^1} guarantees that for $\mc{R} \in \Sigma(\alpha)$, there are no elements $\mc{S}\in\Sigma (\alpha)$ such that $t(\mc{S})=s(\mc{R})-1$ or $s(\mc{S})=t(\mc{R})+1$.  
Thus we have
$$
l_{s(\mc{R})-1}(\alpha) < l_{s(\mc{R})}(\alpha)
\text{ and } l_{t(\mc{R})}(\alpha) > l_{t(\mc{R})+1}(\alpha).
$$
Let $s \le t$ be integers such that $l_{s-1}(\alpha) < l_s(\alpha)=\dots=l_t(\alpha)>l_{t+1}(\alpha)$.
Then we are in the situation of Lemma A (if $s=t$) or Lemma B (if $s<t$).
\qed\\


\begin{remark}\label{remark:R_0}
Take an arbitrary element $\mc{R} \in \Sigma(\alpha)$.
Then, in the proof above, we can 
find $s, t$ such that $s(\mc{R}) \le s \le t \le t(\mc{R})$.
Thus Lemma A or B provides
$$
\Phi\in\Span{T_{\mc{O}_{C_l}(a)} \bigm|  a\in\ZZ, C_l\subset\Supp\mc{R} }
$$
such that $l(\alpha)>l(\Phi(\alpha))$.
We shall use this remark in \S \ref{section:a_n}.
\end{remark}


\begin{corollary}\label{corollary:allspherical}
$B=\Span{T_\alpha \bigm| \alpha \in D_Z(X), \textit{ spherical } }$.
\end{corollary}
\begin{proof}
$B$ is obviously contained in the right hand side.
For a spherical object $\alpha$, Proposition \ref{proposition:step -1 of A_n}
provides $\Psi \in B$ such that $\Psi(\alpha) \cong \owe_{C_b}(a)[i]$
for some $b$, $a$ and $i$.
Then Lemma \ref{easy} (i) shows
$$
T_{\alpha} \cong \Psi^{-1} \circ T_{\owe_{C_b}(a)} \circ \Psi,
$$
which is in $B$.
\end{proof}

\section{Proof of Proposition \ref{proposition:step -2 of A_n}}\label{section:a_n}

The aim of this section is to show Proposition \ref{proposition:step -2 of A_n}.
In the situation of Proposition \ref{proposition:step -2 of A_n},
put $\alpha=\Phi(\owe_{C_1})$ and $\beta=\Phi(\owe_{C_1}(-1))$. 
By Proposition \ref{proposition:step -1 of A_n}, we may assume $l(\alpha)=1$,
and hence $\Supp\alpha=C_b$ for an integer $b$ $(1\le b\le n)$.
The main part of the proof is the following.
\begin{claim}\label{claim:lbeta}
In this situation, suppose $l(\beta)>1$.
Then, there is an autoequivalence $\Psi\in B$ such that 
\begin{equation*}
l(\Psi(\alpha))=1 \text{ and } l(\beta)>l(\Psi(\beta)).
\end{equation*}
\end{claim}
In fact, Proposition \ref{proposition:step -2 of A_n} easily follows from this:\\

\noindent
{\it Proof of Proposition \ref{proposition:step -2 of A_n}.}
By Claim \ref{claim:lbeta}, we can reduce the problem to the case $l(\alpha)=l(\beta)=1$.
In this case, the supports of $\alpha$ and $\beta$ must be the same,
since $\chi(\alpha, \beta)=2$.
Therefore, we get the conclusion from the $A_1$ case.
\qed\\

Thus, the rest of this section is devoted to showing Claim \ref{claim:lbeta}.
In \S \ref{subsection:conditions}, we list conditions on $\alpha$ and $\beta$;
our arguments in the subsequent subsections are based on these conditions.
We divide the proof of Claim \ref{claim:lbeta} into three cases in \S \ref{subsection:division}.
We find $\Psi$ in the three cases in the remaining three subsections.

\subsection{Conditions on $\alpha$ and $\beta$}\label{subsection:conditions}
Before doing computation,
we list conditions that we assume for simplicity or
that our situation imposes on the spherical objects $\alpha$ and $\beta$.

We use the shift functor and a line bundle to simplify the computation
as in Lemma \ref{lemma:normal sub}.
First, using the shift functor $[i]$ $(i\in\ZZ)$,
we may assume that $\alpha$ is a sheaf on $X$ and therefore
$$
\alpha\cong \owe_{C_{b}}(a)
$$
for some $a\in\ZZ$.
\noindent
Secondly, we take a tensor product with a suitable line bundle to assume:

\begin{condition}\label{condition:(1)}
$\max \{\deg_{C_b}\mc{R}\,|\,\mc{R}\in\Sigma(\beta), \Supp \mc{R} \supset C_b \}=0$.
Especially, $\deg_{C_b}\mc{R}=0$ or $-1$ 
for all $\mc{R}\in\Sigma(\beta)$ with $\Supp \mc{R} \supset C_b$ by \eqref{Ext^1}.
\end{condition}
\noindent
Sometimes we also put conditions on the degrees on other curves, depending on the cases.

Relations between $\owe_{C_1}$ and $\owe_{C_1}(-1)$ impose conditions on $a$ and $\beta$.
>From the spectral sequence 
\begin{equation}\label{equation:spectral}
E_2^{p, q}=\Hom_X^p(\cohom^{-q}(\beta), \mc{O}_{C_b}(a)) \Longrightarrow
\Hom_{D(X)}^{p+q}(\beta, \alpha)= \begin{cases}
                                             \CC^2 & p+q=0 \\
                                             0 & p+q\ne 0
                                             \end{cases}
\end{equation}
we obtain
\begin{condition}\label{condition:(2)}
$E_2^{1, q}=0$ for $q\ne -1$
\end{condition}
\noindent
and
\begin{condition}\label{condition:(3)}
$d_2^{0,-1}:E_2^{0,-1}\to E_2^{2,-2}$ is injective, 
$d_2^{0, 0}:E_2^{0, 0}\to E_2^{2,-1}$ is surjective,
and $d_2^{0, q}:E_2^{0, q}\to E_2^{2,q+1}$ are isomorphic for all $q\ne 0,-1$.
\end{condition}

\noindent
In addition to Conditions \ref{condition:(2)} and \ref{condition:(3)}, \eqref{equation:spectral} implies
\begin{equation}\label{equation:addCondition(3)}
\dim \Coker d_2^{0, -1}+\dim \Ker d_2^{0,0}+\dim E_2^{1, -1} =2.
\end{equation}

Moreover, note that the following holds.

\begin{condition}\label{condition:chern class}
$c_1(\alpha) = c_1(\beta) (=C_b)$
holds in the Chow group of curves on $X$.
\end{condition}
\begin{proof}
Let us denote the Grothendieck group of $D_Z(X)$ by $K_Z(X)$ and 
the Euler form on it by $\chi(-,-):K_Z(X)\times K_Z(X)\to \ZZ$. 
Then for a point $x$ in $Z$, we have 
$$
\ZZ[\owe _x]=\bigl\{a\in K_Z(X)\bigm | \chi(a,b)=0 \mbox{ for all } b\in K_Z(X)\bigr \},
$$
since $\chi(-,-)$ is non-degenerate 
on $K_Z(X)/\ZZ[\owe _x]\cong \bigoplus _{i=1} ^n \ZZ[\owe _{C_i}]$. Now
$\Phi$ induces an isometry $\varphi$ on $K_Z(X)$ and it preserves $\ZZ[\mathcal{O}_x]$ by the above equality.
Because $[\owe_{C_1}]-[\owe_{C_1}(-1)]=[\owe _x]$ and $[\alpha]-[\beta]=[\Phi (\owe _x)]$, we get the result.
\end{proof}

\subsection{More on $a, \beta$ and the division into cases}\label{subsection:division}
>From now on, we don't use $\Phi$ in the argument.
In fact, it is sufficient to suppose
that we are given $\alpha \cong \owe_{C_b}(a)$ and a spherical object $\beta$ satisfying the conditions listed above.


\begin{claim}\label{claim:a>=-1}
We have $a\ge -1$.
\end{claim}

\begin{proof}
First note that since $c_1(\beta)=C_b$, there is an integer $q\ne 1$ such that $\mc{H}^q(\beta)\ne 0$.
Assume that $a\le -2$ and let $\mc{R} \in \Sigma(\beta) $ be a direct summand of $\bigoplus_{q\ne 1}\mc{H}^q(\beta)$.
Then it follows from
Conditions \ref{condition:(1)} and \ref{condition:(2)} that $\deg_{C_b}\mc{R}=-1$ and $a=-2$.
Therefore Condition \ref{condition:(1)} implies that there is a direct summand $\mc{R}' \in \Sigma(\beta)$
of $\mc{H}^1(\beta)$ such that $\Supp \mc{R}' \supset C_b$ and $\deg_{C_b}\mc{R}'=0$.
Especially, we have
$$
\Hom_X^0(\mc{O}_{C_b}(-2),\cohom^1(\beta)) \ne 0.
$$
On the other hand, Condition \ref{condition:(1)} also implies
$$
E^{0,0}_2= \Hom_X^0(\mc{H}^0(\beta),\mc{O}_{C_b}(-2))=0
$$
in \eqref{equation:spectral} and accordingly we obtain
$$
\Hom_X^0(\mc{O}_{C_b}(-2),\mc{H}^1(\beta))^\vee\cong E^{2,-1}_2=0
$$
by Condition \ref{condition:(3)}, a contradiction to the non-vanishing above.
\end{proof}
We sometimes use the following useful fact in the latter subsections.
\begin{claim}\label{claim:useful}
Fix $q \ne 0$.
If $E^{2,-q-1}_2= 0$ in \eqref{equation:spectral}, then we have $\deg_{C_b}\mc{R}>a$ for all direct summands $\mc{R} \in \Sigma(\beta)$ of $\mc{H}^q(\beta)$ with $\Supp \mc{R} \supset C_b$.
If, in addition, we suppose that $a\ge 0$, then we get $C_b\not\subset\Supp\mc{H}^q(\beta)$.
\end{claim}

\begin{proof}
The assumption and Condition \ref{condition:(3)} show that 
$$
\Hom_X^0(\mc{H}^q(\beta),\mc{O}_{C_b}(a))=E_2^{0,-q}=0,
$$ 
which implies the first statement.
Then the second statement follows from Condition \ref{condition:(1)}.
\end{proof}

Now we divide the proof into cases.
If there is an element $\mc{R}\in\Sigma (\beta)$ with $\Supp\mc{R}\cap C_{b}=\emptyset $,
then we can find $\Psi\in \Span{T_{\mc{O}_{C_l}(a)} \bigm| a\in\ZZ, C_l \subset \Supp \mc{R} }$
 such that $\Psi(\alpha)\cong \alpha$ and $l(\beta)>l(\Psi(\beta))$
by Remark \ref{remark:R_0}.
Therefore we may assume that
$$
\Supp\mc{R}\cap C_b\ne \emptyset
$$
for all $\mc{R}\in\Sigma (\beta)$ and
we have only to consider the three cases:
\begin{division}
We divide the proof of Claim \ref{claim:lbeta} into the following cases.
\begin{enumerate}
\item[(i)] $C_b\subset\Supp\mc{R}$ for all $\mc{R}\in\Sigma(\beta)$,
\item[(ii)] 
there is $\mc{R} \in \Sigma(\beta)$ with $\Supp \mc{R} \cap C_b = C_{b+1} \cap C_b$
but there is not $\mc{R}' \in \Sigma(\beta)$ with $\Supp \mc{R}' \cap C_b = C_{b-1} \cap C_b$,
\item[(iii)]
there are $\mc{R}, \mc{R}' \in \Sigma(\beta)$ with $\Supp \mc{R} \cap C_b = C_{b+1} \cap C_b$
and $\Supp \mc{R}' \cap C_b = C_{b-1} \cap C_b$.
\end{enumerate}
We subdivide the Case (i) according to the value of $a$:
(i.1) $a \ge 1$, (i.2) $a=0$, and (i.3) $a=-1$.
We also subdivide Case(ii) into
(ii.1) $a =0$ and (ii.2) $a=-1$, after showing $a \le 0$.
We further subdivide (ii.1) and (ii.2) into two cases respectively.
\end{division}


\subsection{Case (i)
}

\paragraph{Case (i.1): $a\ge 1$.}
In this case, it follows from Condition \ref{condition:(1)} that
$$
E^{2,-2}_2\cong \Hom_X^0(\mc{O}_{C_b}(a),\mc{H}^2(\beta))^\vee=0
$$
in \eqref{equation:spectral}.
Hence Claim \ref{claim:useful} and the case assumption show that $\mc{H}^1(\beta)=0$
and consequently that $E_2^{1,q}=0$ for all $q$ in Condition \ref{condition:(2)}.
Then Condition \ref{condition:(1)} implies that $a=1$ and $\Sigma(\beta)=\{\mc{O}_{C_b}\}$. This case has been already treated 
in Proposition \ref{proposition:A_1}.  

\paragraph{Case (i.2): $a=0$.}
\begin{claim}\label{(i.2)}
$\mc{O}_{\cdots\cup C_b}(*,-1),\mc{O}_{C_b\cup\cdots}(-1,*), \mc{O}_{\cdots\cup C_b\cup\cdots}(*,-1,*)\not\in\Sigma(\beta)$.
\end{claim}
\begin{proof}
Note that any sheaf $\mc{R}$ in the assertion satisfies $\Hom^1_X(\owe_{C_b}, \mc{R}) \ne 0$.
Thus, if $\owe_{C_b} \in \Sigma(\beta)$, then the assertion follows from \eqref{Ext^1}.
Therefore we may assume that $\mc{O}_{C_b}\not\in\Sigma(\beta)$.
Under this assumption, the same argument as in Case (i.1) shows that $E_2^{1,q}=0$ for all $q$ in \eqref{equation:spectral}.
It follows that the sheaves in the assertion cannot be in $\Sigma(\beta)$.
\end{proof}
\noindent
By Claim \ref{(i.2)}, we see that
$l(\mc{R}) \ge l(T_{\mc{O}_{C_b}(-1)}(\mc{R}))$ for all $\mc{R} \in \Sigma(\beta)$
and that the inequality is strict if $\mc{R}=\mc{O}_{\cdots\cup C_b}(*,0)$ or $\mc{O}_{C_b\cup\cdots}(0,*)$.
Hence, if $l(\beta)=l(T_{\mc{O}_{C_b}(-1)}(\beta))$, then $\Sigma(\beta)$ consists only of
$\owe_{C_b}(*)$ and $\mc{O}_{\cdots\cup C_b\cup\cdots}(*,0,*)$.
Now we know $c_1(\beta)=C_b$ from Condition \ref{condition:chern class} and therefore $\Sigma(\beta)$ must contain $\owe_{C_b}(*)$.
Then, Lemma \ref{lemma:split} shows $\Supp \beta = C_b$ and Proposition 
\ref{proposition:A_1}
completes the proof for the case (i.2).
\qed
\paragraph{Case (i.3): $a=-1$.}
We put
$$
\bigoplus_p\mathcal{H}^p(\beta)=
\bigoplus_j^{r_1}\mathcal{R}_{1,j}
\oplus\bigoplus_j^{r_2}\mathcal{R}_{2,j}
\oplus\bigoplus_j^{r_3}\mathcal{R}_{3,j}
\oplus\bigoplus_j^{r_4}\mathcal{R}_{4,j}
$$
where $\mathcal{R}_{k,j}$'s are sheaves as follows:
\[ 
\xymatrix@R=1ex@M=.3ex{     & & C_b &  \\
          \mathcal{R}_{1,j}:& \ar@{-}[r]&{\no}\ar@{-}[r]&\\
          \mathcal{R}_{2,j}:& &{\no}\ar@{-}[r] &\\
          \mathcal{R}_{3,j}:& &{\no} &  \\ 
          \mathcal{R}_{4,j}:&  \ar@{-}[r]&{\no} &    
          }
\]
 When $\Supp\beta=C_b$, we can apply Proposition \ref{proposition:A_1}, and hence we may assume that 
 $\Supp\beta\ne C_b$.
On the other hand, since $c_1(\beta)=C_b$, we can see either $r_3\ne 0$ or $r_2\cdot r_4\ne 0$ holds.  
Therefore the proof of Lemma A in \S\ref{subsection:lemmaA} implies 
$l(\beta)>l(\Psi(\beta))$ for $\Psi=T_{\mc{O}_{C_b}(-1)}$ or $T_{\mc{O}_{C_b}(-2)}$.
In each case, we can see $l(\Psi(\alpha))=1$.


\subsection{Case (ii)
}

The existence of $\mc{R} \in \Sigma(\beta)$ with $\Supp \mc{R} \cap C_b=C_b \cap C_{b+1}$
and \eqref{Ext^1} imply the non-existence of $\mc{S} \in \Sigma(\beta)$ with
$\Supp \mc{S} \cap C_{b+1} = C_b\cap C_{b+1}$.
Thus we have
$$
\Sigma(\beta)\subset\bigl\{ \mc{O}_{C_b\cup\cdots}(a',*), \mc{O}_{C_{b+1}}(*),
\mc{O}_{C_{b+1}\cup\cdots}(*), \mc{O}_{\cdots\cup C_b\cup\cdots}(*,a',*)\bigm|a'=-1,0 \}.
$$
By Condition \ref{condition:(2)}, $\mc{R}$ as above exists only in $\cohom^1(\beta)$.
Moreover, because of the condition $c_1(\beta)=C_b$, $\cohom^1(\beta)$ has
precisely one such direct summand $\mc{R}$.
It also follows from Condition \ref{condition:(2)} and Claim \ref{claim:a>=-1} that $a=-1$ or $0$.

\paragraph{Case (ii.1): $a=0$.}
In this case,
$$
E^{2,q}_2\cong \Hom_X^0(\mc{O}_{C_b},\mc{H}^{-q}(\beta))^\vee=0
$$
holds for all $q$ in \eqref{equation:spectral}.
Therefore, Claim \ref{claim:useful} implies that
$\mc{H}^q(\beta)=0$ for $q\ne 0,1$ and that $\Supp \cohom^1(\beta) \not\supset C_b$.
Then, from Condition \ref{condition:(1)} and the condition $c_1(\beta)=C_b$,
we can see
\begin{equation*}
(\mc{H}^0(\beta),\mc{H}^1(\beta))=(\mc{O}_{C_b\cup\cdots\cup C_{b''}}(0,*),\mc{O}_{C_{b+1}\cup\cdots\cup C_{b''}}(*))
\end{equation*}
with $b+1\le b''$.
Applying Lemma \ref{middle} (n.b. $C_s$ in Lemma \ref{middle} is $C_{b''}$ here),
we may assume that $\deg_{C_l}\mc{H}^q(\beta)=0$ for all $l$ $(b+1<l<b'')$ and all $q$.
Now we can classify spherical objects with such cohomology sheaves.
Note that by virtue of Lemma \ref{decomposable}, we have
\begin{equation}\label{equation:nonsplit}
\Ext^2_X(\cohom^1(\beta),\cohom^0(\beta))\ne 0.
\end{equation}
We divide the proof into two cases:
\subparagraph{Case (ii.1.a): $b+1<b''$.}
In this case, we may assume $\deg_{C_{b"}} \cohom^0(\beta) \ne \deg_{C_{b"}} \cohom^1(\beta)$
by Lemma \ref{lemma:fundamental}.
Then, by virtue of \eqref{Ext^1} and the conditions listed above,
the cohomology sheaves of $\beta$ must be of the following forms,
up to tensoring a line bundle:
\[ 
\xymatrix@R=1ex@M=.3ex{ & C_b& C_{b+1} & C_{b+2} & & C_{b''-1} & C_{b''}\\
        \mc{H}^0(\beta):& {\zero}\ar@{-}[r]&{\zero}\ar@{-}[r]&{\zero}\ar@{-}[r]&\cdots \ar@{-}[r]&{\zero}\ar@{-}[r]&{\mone} \\
        \mc{H}^1(\beta):&                     &{\mone}\ar@{-}[r]&{\zero}\ar@{-}[r]&\cdots \ar@{-}[r]&{\zero}\ar@{-}[r]&{\zero}}
\]
In this case, $\Psi:=T_{\mc{O}_{C_b}}\circ T_{\mc{O}_{C_{b+1}}(-2)}$ satisfies the conditions
$l(\beta)>l(\Psi (\beta))$ and $l(\Psi (\alpha))=1$ as desired.
\subparagraph{Case (ii.1.b): $b+1=b''$.}
In this case, \eqref{Ext^1}, \eqref{equation:nonsplit} and Condition \ref{condition:(1)} show
$$
(\cohom^0(\beta),\cohom^1(\beta))=(\mc{O}_{C_b\cup C_{b+1}},\mc{O}_{C_{b+1}}),
$$
up to tensoring a line bundle.
Then we can see that 
$T_{\mc{O}_{C_{b+1}}(-1)}(\beta)=\mc{O}_{C_b\cup C_{b+1}}(1,-2)$ and 
$T_{\mc{O}_{C_{b+1}}(-1)}(\alpha)=\mc{O}_{C_b\cup C_{b+1}}(1,-1)$.  
Hence we obtain $l(\beta)>l(\Psi(\beta))$ and $l(\Psi(\alpha))=1$, 
where $\Psi= T_{\mc{O}_{C_b}}\circ T_{\mc{O}_{C_{b+1}}(-1)}$.

\paragraph{Case (ii.2): $a=-1$.}
By the argument in the beginning of Case (ii), we can write
\begin{equation}\label{equation:decompose}
\bigoplus_p\mathcal{H}^p(\beta)=
\bigoplus_j^{r_0}\mathcal{R}_{0,j}\oplus
\bigoplus_j^{r_1}\mathcal{R}_{1,j}\oplus\bigoplus_j^{r_2}\mathcal{R}_{2,j}
\oplus\bigoplus_j^{r_3}\mathcal{R}_{3,j}\oplus\mathcal{R}_4,
\end{equation}
where $\mathcal{R}_{k,j}$ and $\mathcal{R}_4$ are sheaves of the forms in
the following figure.
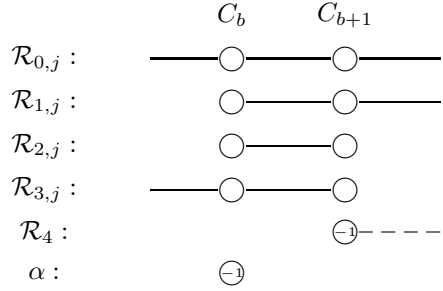
\[ 
\xymatrix@R=1ex@M=.3ex{       &  & C_b & C_{b+1} & \\
    \mc{R}_{0,j} : &\ar@{-}[r] &{\no}\ar@{-}[r] &{\no}\ar@{-}[r]& \\
          \mathcal{R}_{1,j}  :&           &{\no}\ar@{-}[r]  &{\no}\ar@{-}[r]& \\
            \mathcal{R}_{2,j}:&           &{\no}\ar@{-}[r]  &{\no}        & \\
         \mathcal{R}_{3,j}   :& \ar@{-}[r]&{\no}\ar@{-}[r]  &{\no}        &   \\
            \mathcal{R}_4    :&           &                   &{\mone}\ar@{--}[r]& \\
            \alpha           :&           &{\mone}          &                 &}
\]
Here, noting that $\mc{R}_4 \subset \cohom^1(\alpha)$ is unique,
we normalize the degrees on $C_{b+1}$ by the condition
$$
\deg_{C_{b+1}}\mathcal{R}_4=-1.
$$
In addition, we can see
\begin{equation}\label{equation:deg0}
\deg_{C_b}\mc{R}_{0,j}=\deg_{C_b}\mc{R}_{3,j}=0
\end{equation}
as follows.
If $\mc{O}_{C_b\cup\cdots}(0,*)\in\Sigma (\beta)$, then \eqref{Ext^1} and Condition \ref{condition:(1)} imply \eqref{equation:deg0}.
Thereby assume $\mc{O}_{C_b\cup\cdots}(0,*)\not\in\Sigma (\beta)$.
Then we get  
$$
E^{2,-2}_2\cong \Hom_X^0(\mc{O}_{C_b}(-1),\mc{H}^2(\beta))^\vee=0
$$
in Claim \ref{claim:useful} and therefore $\deg_{C_b}\mc{R}_{0,j}$
(or $\deg_{C_{b}}\mc{R}_{3,j}$) is zero
if it is a direct summand of $\cohom^1(\beta)$.
>From this and Condition \ref{condition:(2)}, we conclude that \eqref{equation:deg0} holds
for all $j$.

As a consequence of \eqref{equation:deg0} and the uniqueness of $\mc{R}_4$,
we have $\dim E_2^{1,-1}=1$ in \eqref{equation:spectral}.
Thus \eqref{equation:addCondition(3)} becomes
\begin{equation}\label{equation:cok+ker}
\dim \Coker d_2^{0, -1}+\dim \Ker d_2^{0,0}=1.
\end{equation}
Now we divide the proof of Case (ii.2) into the two cases: (a) $l(\mc{R}_4)>1$ and (b) $l(\mc{R}_4)=1$
\subparagraph{Case (ii.2.a): $l(\mc{R}_4)>1$.}
In this case, \eqref{Ext^1} implies
that $\deg_{C_{b+1}}\mathcal{R}_{2,j}=\deg_{C_{b+1}}\mathcal{R}_{3,j}=-1$
and that $\deg_{C_{b+1}}\mathcal{R}_{0,j}$, $\deg_{C_{b+1}}\mathcal{R}_{1,j}
\in \{0, -1\}$.
Thus, specifying degrees in \eqref{equation:decompose}, we write
 
$$
\bigoplus_p\mathcal{H}^p(\beta)=
\bigoplus_j^{r_0}\mathcal{R}_{0,j}
\oplus\bigoplus_j^{s_1}\mathcal{S}_{1,j} \oplus \cdots \oplus \bigoplus_j^{s_6}\mathcal{S}_{6,j}
\oplus\bigoplus_j^{r_3}\mathcal{R}_{3,j}\oplus\mathcal{R}_4,
$$
where $\mathcal{S}_{k,j}$'s are sheaves of the forms in the following figure.
\[ 
\xymatrix@R=1ex@M=.3ex{     &  & C_b & C_{b+1} & \\
\mc{R}_{0,j}:& \ar@{-}[r]&{\zero}\ar@{-}[r]&{\no}\ar@{-}[r]& \\
          \mathcal{S}_{1,j}:&          &{\zero}\ar@{-}[r]&{\zero}\ar@{-}[r]& \\
          \mathcal{S}_{2,j}:&          &{\mone}\ar@{-}[r]&{\zero}\ar@{-}[r]& \\ 
          \mathcal{S}_{3,j}:&          &{\zero}\ar@{-}[r]&{\mone}\ar@{-}[r]&\\
          \mathcal{S}_{4,j}:&          &{\mone}\ar@{-}[r]&{\mone}\ar@{-}[r]& \\
          \mathcal{S}_{5,j}:&          &{\zero}\ar@{-}[r]&{\mone}& \\
          \mathcal{S}_{6,j}:&         &{\mone}\ar@{-}[r]&{\mone}& \\ 
          \mathcal{R}_{3,j}:& \ar@{-}[r]&{\zero}\ar@{-}[r]&{\mone}& \\
          \mathcal{R}_4:&           &                 &{\mone}\ar@{-}[r]& }
\]
Then \eqref{equation:cok+ker} and Condition \ref{condition:(3)} imply that
\begin{equation}\label{C_b}
|s_1+s_3+s_5-s_2-s_4-s_6|= 1.
\end{equation}
We first consider $\Psi'=T_{\mc{O}_{C_b}(-1)}\circ T_{\mc{O}_{C_{b+1}}(-2)}$ and note that
$\Psi '(\alpha)\cong \mc{O}_{C_{b+1}}(-2)$. Moreover we obtain
$$
\sum_p l(\Psi '(\cohom ^p(\beta)))-l(\beta)= s_1-s_2-s_3+s_4-s_5-s_6-2r_3-1
$$ 
from direct computation.
Then by Lemma \ref{strategy}, we have 
\begin{equation}\label{beta}
l(\Psi '(\beta))-l(\beta)\le s_1-s_2-s_3+s_4-s_5-s_6-2r_3-1.
\end{equation}
>From (\ref{C_b}) and \eqref{beta}, we get  
\begin{equation}\label{RHS}
l(\Psi '(\beta))-l(\beta) \le 2s_4-2s_3-2s_5-2r_3
\end{equation}
and 
\begin{equation}\label{RHS'}
l(\Psi '(\beta))-l(\beta) \le 2s_1-2s_2-2s_6-2r_3.
\end{equation}
If $l(\beta)> l(\Psi '(\beta))$, then we have nothing to do any more. Hence 
let us consider the case
$$
l(\beta)\le l(\Psi '(\beta)).
$$
Now note that $s_1\cdot s_4=0$ by \eqref{Ext^1}.
If $s_1=0$, \eqref{RHS'} implies $s_2=s_6=r_3=0$ and $l(\Psi '(\beta))=l(\beta)$.
Then \eqref{beta} means
\begin{equation*}
s_4\ge s_3+s_5+1.
\end{equation*}
It follows from this that $s_4\ne 0$, which implies $s_5=0$ and
$\deg_{C_{b+1}} \mc{R}_{0,j}=-1$ by \eqref{Ext^1}.
Hence in this case, we have 
$$
l(T_{\mc{O}_{C_b\cup C_{b+1}(-1,-2)}}(\beta))-l(\beta) \le 2s_3-2s_4+1 \le -1
$$ 
and
$T_{\mc{O}_{C_b\cup C_{b+1}(-1,-2)}}(\alpha)\cong\mc{O}_{b+1}(-3)[1]$ as desired.
If $s_4=0$, by a similar argument, we see $s_3=s_5=s_6=r_3=0$,
$\deg_{C_{b+1}} \mc{R}_{0,j}=0$ and $s_1\ge s_2+1$. Then we obtain
$$
l(T'_{\mc{O}_{C_b\cup C_{b+1}}}(\beta))-l(\beta) \le 2s_2-2s_1+1 \le -1
$$
and $T'_{\mc{O}_{C_b\cup C_{b+1}}}(\alpha)\cong\mc{O}_{b+1}[-1]$, which finishes the proof.
  
\paragraph{Case(ii.2.b): $a=-1$ and $l(\mc{R}_4)=1$.}
In this case, \eqref{Ext^1} implies that $\deg_{C_{b+1}} \mc{R}_{0,j}
= \deg_{C_{b+1}} \mc{R}_{1,j}=0$.
Noting \eqref{equation:deg0}, we specify the degrees in \eqref{equation:decompose} and write
$$
\bigoplus_p\mathcal{H}^p(\beta)=
\bigoplus_j^{r_0} \mc{R}_{0,j} \oplus
\bigoplus_j^{s_1}\mathcal{S}_{1,j}\oplus
\dots
\oplus\bigoplus_j^{s_8}\mc{S}_{8,j}
\oplus\mathcal{R}_4,
$$
where $\mathcal{S}_{k,j}$'s are sheaves of the following forms.
\[ 
\xymatrix@R=1ex@M=.3ex{     &  & C_b & C_{b+1} & \\
\mc{R}_{0,j}: & \ar@{-}[r] & {\zero} \ar@{-}[r] & {\zero} \ar@{-}[r] & \\
          \mathcal{S}_{1,j}:&          &{\zero}\ar@{-}[r]&{\zero}\ar@{-}[r]& \\
          \mathcal{S}_{2,j}:&          &{\mone}\ar@{-}[r]&{\zero}\ar@{-}[r]& \\ 
          \mathcal{S}_{3,j}:&          &{\zero}\ar@{-}[r]&{\zero}&\\
          \mathcal{S}_{4,j}:&          &{\zero}\ar@{-}[r]&{\mone}& \\
          \mathcal{S}_{5,j}:&          &{\mone}\ar@{-}[r]&{\zero}& \\
          \mathcal{S}_{6,j}:&         &{\mone}\ar@{-}[r]&{\mone}& \\ 
          \mathcal{S}_{7,j}:& \ar@{-}[r]&{\zero}\ar@{-}[r]&{\mone}& \\
          \mathcal{S}_{8,j}:& \ar@{-}[r]&{\zero}\ar@{-}[r]&{\zero}& \\
          \mathcal{R}_4:&           &                 &{\mone}& }
\]
\begin{claim}\label{claim:ii-2-b}
Under the above assumption, we have the following.
\begin{enumerate}
\item
$|s_1+s_3+s_4-(s_2+s_5+s_6)|=1$.
\item
If $s_3=s_6=s_7=s_8=0$, then we have $s_1 \ne s_2$.
\item
If $s_1=s_2=s_3=s_6=0$, then we have $s_7 \ne s_8$.
\end{enumerate}
\end{claim}
\begin{proof}
(i) follows from Condition \ref{condition:(3)} and \eqref{equation:cok+ker}.
To show (ii), assume
\begin{equation}\label{equation:ii.2.b.1}
s_3=s_6=s_7=s_8=0
\end{equation}
and
\begin{equation}\label{equation:ii.2.b.2}
s_1=s_2.
\end{equation}
(i) means that $|s_4-s_5|=1$ in this case.
Write $s_k=\sum_p s_k^p$ where $s_k^p$ counts the number of
direct summands $\mc{S}_{k,j}$ in $\cohom^p(\beta)$.
By \eqref{equation:ii.2.b.2} we can apply Lemma \ref{lemma:equality}
to deduce that $s_1^p=s_2^{p-1}$ for all $p$.
On the other hand, Condition \ref{condition:(3)} under the assumption \eqref{equation:ii.2.b.1}
gives rise to equalities and inequalities 
$s_5^p+ s_2^p=s_4^{p+1}+s_1^{p+1}$ for $p \ne 0, 1$,
$s_5^0+ s_2^0 \ge s_4^1+s_1^{1}$ and $s_5^{1}+ s_2^1 \le s_4^2+s_1^{2}$.
Thus we obtain
$s_5^p=s_4^{p+1}$ for $p \ne 0, 1$,
$s_5^0 \ge s_4^1$ and $s_5^{1} \le s_4^2$.
Moreover, \eqref{equation:cok+ker} says either
$$
\left\{
\begin{aligned}
s_5^0 &= s_4^1+1 \\
s_5^{1} &= s_4^2
\end{aligned}
\right.
\qquad
\text{or}
\qquad
\left\{
\begin{aligned}
s_5^0 &= s_4^1 \\
s_5^{1} &= s_4^2 -1
\end{aligned}
\right.
$$
holds.
We consider only the first case because the second case is similar.
In this case, Lemma \ref{lemma:equality} applied to $C=C_{b+1}$
yields $ s_4^{p-1}\le s_5^p$ for all $p$.
Then we have
$$
s_5^2 = s_4^3 \le s_5^4 = s_4^5 \le \cdots \quad \textrm{and} \quad s_4^0 \le s_5^1 = s_4^2 \le s_5^3 = s_4^4\le \cdots.
$$
Because $\beta$ is a bounded complex, we have $s_4^p=s_5^p=0$ for $p\gg 0$, and consequently $s_5^2=s_4^0=0$.
It follows from this and \eqref{equation:ii.2.b.1}
that
$$\Ext^2_X(\cohom^2(\beta), \mc{R}_4)=\Ext^2_X(\mc{R}_4, \cohom^0(\beta))=0.$$
Recall $\mc{R}_4$ is a direct summand of $\cohom^1(\beta)$ by Condition \ref{condition:(2)}.
Then Lemma \ref{decomposable} implies that $\mc{R}_4[-1]$ is a direct summand of $\beta$.
Since $\beta$ is spherical, this means that $\beta = \mc{R}_4[-1]$
and hence that $c_1(\alpha) \ne c_1(\beta)$.
This is a contradiction to Condition \ref{condition:chern class}.
(iii) can be shown in a similar way.
\end{proof}
Since $c_1(\beta)=C_b$ holds by Condition \ref{condition:chern class},
we see that $r_0+s_1+s_2$ is even
and $r_0+s_1+\dots+s_8$ is odd.
Therefore, $s_3+\dots+s_8$ is odd and especially we have $s_3+s_5+s_8 \ne s_4+s_6+s_7$.
Since $l(T_{\owe_{C_{b}}(-1)}\circ T_{\owe_{C_{b+1}}(k)}(\alpha))=1$ for all $k$,
the following completes the proof for the case (ii.2.b).

\begin{claim}
\begin{enumerate}
\item
If $s_3+s_5+s_8 > s_4 + s_6 + s_7$, then
$l(T_{\owe_{C_{b}}(-1)}\circ T_{\owe_{C_{b+1}}(-1)}(\beta)) < l(\beta)$.
\item
If $s_3+s_5+s_8 < s_4 + s_6 + s_7$, then
$l(T_{\owe_{C_{b}}(-1)}\circ T_{\owe_{C_{b+1}}(-2)}(\beta)) < l(\beta)$.
\end{enumerate}
\end{claim}
\begin{proof}
To prove (i), suppose that the inequality
\begin{equation}\label{equation:assumption of (i)}
s_3+s_5+s_8 > s_4 + s_6 + s_7
\end{equation}
holds.
If we further assume $s_6 \ne 0$,
then \eqref{Ext^1} implies $r_0=s_1=s_3=s_8=0$ and \eqref{equation:assumption of (i)} becomes
$s_5 > s_4 + s_6 + s_7$.
This contradicts $|s_4-(s_2+s_5+s_6)|=1$ from Claim \ref{claim:ii-2-b}
and thus we obtain
$$s_6=0.$$
Then, putting $\Psi=T_{\owe_{C_{b}}(-1)}\circ T_{\owe_{C_{b+1}}(-1)}$,
we have
\begin{equation}\label{equation:ii-2-b}
l(\Psi(\beta))- l(\beta) \le (s_2 + s_4 + 2s_7+1) -(s_1 + s_3 + s_5 + 2s_8)
\end{equation}
by Lemma \ref{strategy}.

We first consider the case $s_3 = 0$. By contradiction,
assume that $l(\Psi(\beta))- l(\beta) \ge 0$.
Then, combining \eqref{equation:ii-2-b} with $s_5-s_4=s_1-s_2\pm 1$ from Claim \ref{claim:ii-2-b} (i)
and $s_5-s_4>s_7-s_8$ from \eqref{equation:assumption of (i)},
we see $s_5-s_4=s_1-s_2+1$ and $s_1 - s_2 = s_7 - s_8$.
Now we have $s_1s_7 = s_2s_8 = 0$ by \eqref{Ext^1} and therefore we obtain
$s_1 = s_2$ and $s_7 = s_8$.
Since either of these is zero, this contradicts Claim \ref{claim:ii-2-b} (ii) and (iii).

Next consider the case $s_3 \ne 0$.
In this case, we have $r_0 = s_2 = s_6 = s_7 = 0$ by \eqref{Ext^1}.
Then \eqref{equation:ii-2-b} and \eqref{equation:assumption of (i)} imply 
\begin{equation}\label{equation:s_4+1}
l(\Psi(\beta))- l(\beta) \le (s_4+1)-(s_1+s_3+s_5+2s_8)\le -s_1-s_8\le 0.
\end{equation}
Assume $l(\Psi(\beta)) = l(\beta)$.
Then the equalities hold in \eqref{equation:s_4+1} and
it follows that $s_1=s_8=0$ and $s_4+1=s_3+s_5$.
Combining it
with $|s_3+s_4-s_5|=1$ from Claim \ref{claim:ii-2-b} (i),
we also see $s_3=1$ and $s_4=s_5$.
Moreover, since the equality holds in \eqref{equation:ii-2-b}, the spectral sequence in Lemma \ref{strategy} must be $E_2$-degenerate.
Namely, for the class $e^p(\beta) \in \Hom_X(\cohom^p(\beta), \cohom^{p-1}(\beta)[2])$,
the map
$$
\cohom^{-1}(\Psi(e^p(\beta))): \cohom^{-1}(\Psi(\cohom^{p}(\beta))) \to \cohom^1(\Psi(\cohom^{p-1}(\beta)))
$$
is zero (see Proposition \ref{differential map}).
Note
\begin{itemize}
\item $\Ext^2_X(\mc{F}, \mc{S}_{3,1})$ for $\mc{F}=\mc{S}_{4, j}, \mc{S}_{5,j},
\mc{R}_4$.
\item The map $\Ext^2_X(\mc{S}_{3,1}, \mc{F}) \to \Hom_X(\cohom ^{-1}(\Psi(\mc{S}_{3,1})), \cohom ^1(\Psi(\mc{F})))$
induced by $\Psi$
is isomorphic for $\mc{F}=\mc{S}_{5,j}, \mc{R}_4$ and of rank $1$ for $\mc{F}=\mc{S}_{4, j}$.
\end{itemize}
Hence, for $\mc{F}$ as above, if an entry of $e^p(\beta)$ in $\Ext^2_X(\mc{F}, \mc{S}_{3,1})$
or $\Ext^2_X(\mc{S}_{3,1}, \mc{F})$ is non-zero,
then it must be in the kernel of
$\Ext^2_X(\mc{S}_{3, 1}, \mc{S}_{4, j}) \to \Ext^2_X(\mc{S}_{3, 1}, \mc{S}_{4, j}|_{C_{b+1}})$.
This contradicts the surjectivity of $d_2^{0, 0}$ in Condition \ref{condition:(3)}.
Thus we obtain (i).
The proof of (ii) is similar.
\end{proof}


\subsection{Case (iii)
} 

Condition \ref{condition:(3)} implies that $\mc{R}$ and $\mc{R}'$ above must be in $\cohom^1(\beta)$.
Moreover, 
they are unique in a decomposition of $\cohom^1(\beta)$,
by virtue of the inequality $\dim E_2^{1, -1} \le 2$ from \eqref{equation:addCondition(3)}.
Thus \eqref{Ext^1} allows us to write
$$
\bigoplus_p\mathcal{H}^p(\beta)=
\bigoplus_j^{r_1}\mathcal{R}_{1,j}\oplus\bigoplus_j^{r_2}\mathcal{R}_{2,j}
\oplus\mathcal{R}_3\oplus\mathcal{R}_4,
$$
where $\mathcal{R}_{k,j}$'s, $\mathcal{R}_3$ and $\mathcal{R}_4$ are sheaves of the following forms.
\[ 
\xymatrix@R=1ex@M=.3ex{     &               & C_{b-1} & C_b &  C_{b+1}& \\
          \mathcal{R}_{1,j}: & \ar@{-}[r]&{\no}\ar@{-}[r]&{\no}\ar@{-}[r]&{\no}\ar@{--}[r]&   &\\         
          \mathcal{R}_{2,j}: &           &{\no}\ar@{-}[r]&{\no}\ar@{-}[r]&{\no}\ar@{--}[r]&   &\\
          \mathcal{R}_3    : & \ar@{--}[r]&{\mone}          &                 &{\no}\ar@{--}[r]& & :\mathcal{R}_4\\
          \alpha:           &          &                  &{\cia}          &                 &}
\]
Here we assume that 
$\deg_{C_{b-1}}\mc{R}_3=-1$ by tensoring a suitable line bundle.

\begin{claim}
We have $a=-1$.
\end{claim}
\begin{proof}
Claim \ref{claim:a>=-1} says $a \ge -1$.
If $a\ge 0$, then we have $\Ext_X^1(\owe_{C_{b}}(a), \mc{R}) \ne 0$ for any $\mc{R} \in \Sigma(\beta)$.
It follows from Condition \ref{condition:(2)} that $\mc{H}^q(\beta)=0$ for $q\ne 1$. This is absurd, since $c_1(\beta)=C_b$ by Condition \ref{condition:chern class}.
\end{proof}

The inequality $\dim E_2^{1, -1} \le 2$ from \eqref{equation:addCondition(3)}
also implies that $\Ext_X^1(\mathcal{R}_{k,j},\mc{O}_{C_b}(-1))=0$ for $k=1, 2$ and for all $j$. 
In particular we get
$$
\deg_{C_b}\mc{R}_{1,j}=\deg_{C_b}\mc{R}_{2,j}=0.
$$

Now we give a proof for Case (iii) by induction on $l(\mc{R}_3)$.
First suppose $l(\mc{R}_3)=1$.
We write
$$
\bigoplus_j^{r_1}\mathcal{R}_{2,j}=\bigoplus_j^{s_1}\mathcal{S}_{1,j}
\oplus\bigoplus_j^{s_2}\mathcal{S}_{2,j},
$$
where $\mathcal{S}_{k,j}$'s are sheaves of the following forms.
\[ 
\xymatrix@R=1ex@M=.3ex{ &               & C_{b-1} & C_b &  C_{b+1}& \\
          \mathcal{R}_{1,j}:& \ar@{-}[r]&{\zero}\ar@{-}[r]&{\zero}\ar@{-}[r]&{\no}\ar@{--}[r]& \\
          \mathcal{S}_{1,j}:&          &{\zero}\ar@{-}[r] &{\zero}\ar@{-}[r]&{\no}\ar@{--}[r]& \\
          \mathcal{S}_{2,j}:&          &{\mone}\ar@{-}[r]&{\zero}\ar@{-}[r]&{\no}\ar@{--}[r]& \\
          \mathcal{R}_3:&              &{\mone}         &                 &{\no}\ar@{--}[r]& &:\mathcal{R}_4& \\
          \alpha:           &          &                  &{\mone}          &                 &}
\]
Because of the existence of $\mc{R}_3$, we have $s_1\ne s_2$ by Lemma \ref{lemma:r_4=0}.
Define
$$
\Psi _0= \begin{cases}
T_{\mc{O}_{C_{b-1}\cup C_b}(-1,-1)} & \text{ if } s_1<s_2, \\
T'_{\mc{O}_{C_{b-1}\cup C_b}} & \text{ if } s_2<s_1.
\end{cases}
$$
Then $(\Psi_0(\alpha), \Psi _0(\beta))$ fits in Case (ii) and
$\Psi_0(\beta)$ satisfies $l(\Psi_0 (\beta))\le l(\beta)$.
Since we have proved Case (ii), we finish the case $l(\mc{R}_3)=1$.

Next suppose $l(\mc{R}_3)>1$.
In this case, \eqref{Ext^1} implies
$$
\deg_{C_{b-1}} \mc{R}_{2,j}=-1.
$$
Define
$$
\Psi'= T_{\owe_{C_b}(-1)} \circ T_{\owe_{C_{b-1}}(-2)}.
$$
Then we have $\Psi'(\alpha) \cong \owe_{C_{b-1}}(-2)$
and $l(\Psi'(\beta)) \le l(\beta)$.
Moreover,
we can see that $\Psi'(\beta)$ satisfies the induction hypothesis (on $l(\mc{R}_3))$.
This finishes the proof of Case (iii) and we get the assertion of 
Proposition \ref{proposition:step -2 of A_n}.
\qed\\

\noindent \textsc{Department of Mathematics, Graduate School of Science, Hiroshima University\\
1-3-1 Kagamiyama, Higashi-Hiroshima 739-8526, JAPAN\\
E-mail:} \texttt{akira@math.sci.hiroshima-u.ac.jp}
 
\medskip
\noindent \textsc{Department of Mathematics, Kyoto University\\
Kyoto 606-8502, Japan\\
E-mail:}
\texttt{hokuto@math.kyoto-u.ac.jp}

\begin{thebibliography}{Hum72}

\bibitem[Art66]{Artin:ira}
M.~Artin.
\newblock On isolated rational singularities of surfaces.
\newblock {\em Amer.\ J.\ {M}ath. \ \emph{\textbf{88}}}, pages 129--136,
  (1966), MR0199191, Zbl 0142.18602.

\bibitem[BO01]{Bondal:rvfd}
A.~Bondal and D.~Orlov.
\newblock Reconstruction of a variety from the derived category and groups of
  autoequivalences.
\newblock {\em Compositio Math.\ \emph{\textbf{125}}}, pages 327--344, (2001), MR1818984, Zbl 0994.18007.

\bibitem[Bri99]{Bridgeland:etc}
T.~Bridgeland.
\newblock Equivalences of triangulated categories and {F}ourier-{M}ukai
  transforms.
\newblock {\em Bull.\ {L}ondon {M}ath.\ {S}oc.\ \emph{\textbf{31}}}, pages
  25--34, (1999), MR1651025, Zbl 0937.18012.

\bibitem[BM98]{Bridgeland:fmqv}
T.~Bridgeland and A.~Maciocia.
\newblock  \emph{Fourier-Mukai transforms for quotient varieties}.
\newblock Preprint math.AG/9811101.

\bibitem[BM01]{Bridgeland:csed}
T.~Bridgeland and A.~Maciocia.
\newblock Complex surfaces with equivalent derived categories.
\newblock {\em Math.\ {Z}.\ \emph{\textbf{236}}}, pages 677--697, (2001), MR1827500, Zbl pre01665938.

\bibitem[Esn85]{Esnault:rmqs}
H.~Esnault.
\newblock Reflexive modules on quotient surface singularities.
\newblock {\em J.\ {R}eine\ {A}ngew.\ {M}ath. \emph{\textbf{362}}}, (1985), MR0809966, Zbl 0553.14016.

\bibitem[GM96]{Gelfand:mha}
S.~I. Gelfand and Yu.~I. Manin.
\newblock {\em Methods of homological algebra}.
\newblock Springer, (1996), MR0323842, Zbl 0855.18001.

\bibitem[Hum72]{Humphreys:LART}
J.~E. Humphreys.
\newblock {\em Introduction to {L}ie algebras and representation theory}.
\newblock Springer, (1972), MR0323842, Zbl 0254.17004.

\bibitem[KV00]{kapranov:kdh}
M.~Kapranov and E.~Vasserot.
\newblock Kleinian singularities, derived categories and {H}all algebras.
\newblock {\em Math.\ {A}nn.\ \emph{\textbf{316}}}, pages 565--576, (2000), MR1752785, Zbl 0997.14001.

\bibitem[KS90]{Kashiwara:som}
M.~Kashiwara and P.~Schapira.
\newblock {\em Sheaves on Manifolds}.
\newblock Springer, (1990), MR1074006, Zbl 0709.18001.

\bibitem[Kaw02]{Kawamata:deke}
Y.~Kawamata.
\newblock ${D}$-equivalence and ${K}$-equivalence.
\newblock {\em J.\ {D}ifferential\ {G}eom.\ \emph{\textbf{61}}}, pages
  147--171, (2002), MR1949787, Zbl pre02052875.

\bibitem[Orl97]{Orlov:edck3}
D.~O. Orlov.
\newblock Equivalences of derived categories and ${K}3$ surfaces.
\newblock In {\em Algebraic geometry, 7.\ {J}.\ {M}ath.\ {S}ci. (New York)
  \emph{\textbf{84}}}. Consultants {B}ureau, {N}ew {Y}ork, pages 1361--1381,
  (1997), MR1465519, Zbl 0938.14019.

\bibitem[Orl02]{Orlov:dcav}
D.~O. Orlov.
\newblock Derived categories of coherent sheaves on abelian varieties and
  equivalences between them.
\newblock {\em Izv.\ {R}oss.\ {A}kad.\ {N}auk\ {S}er.\ {M}at.
  \emph{\textbf{66}}}, pages 131--158, (2002).
\newblock translation in {I}zv.\ {M}ath.\ \emph{\textbf{66}} (2002),\ 569--594, MR1921811, Zbl 1031.18007.

\bibitem[Rie03]{Riemenschneider:srtm}
O.~Riemenschneider.
\newblock Special representations and the two-dimensional {M}cKay
  correspondence.
\newblock {\em Hokkaido {M}ath.\ {J}.\ \emph{\textbf{32}}}, pages 317--333,
  (2003), MR1996281, Zbl 1046.14002.

\bibitem[ST01]{Seidel:bga}
P.~Seidel and R.~Thomas.
\newblock Braid group actions on derived categories of coherent sheaves.
\newblock {\em Duke {M}ath.\ {J}.\ \emph{\textbf{108}}}, pages 37--108, (2001), MR1831820, Zbl pre01820814.

\bibitem[Tod03]{To03}
Y.~Toda.
\newblock \emph{Fourier-Mukai transforms and canonical divisors}.
\newblock Preprint math.AG/0312015, (2003).

\bibitem[Ver96]{Verdier:cdca}
J.~L. Verdier.
\newblock Des cat{\' e}gories d{\' e}riv{\' e}es des cat{\' e}gories ab{\'
  e}liennes.
\newblock {\em Ast{\' e}resque \emph{23}}, (1996), MR1453167, Zbl 0882.18010.

\end{thebibliography}
\end{document}